\documentclass[psamsfonts]{amsart}
\input{diagrams.tex}
\diagramstyle[scriptlabels,height=8mm,width=8mm]
\newcommand\query[1]{}
\theoremstyle{plain}

\newtheorem{thm}{Theorem}[section]
\newtheorem{cor}[thm]{Corollary}
\newtheorem{lem}[thm]{Lemma}
\newtheorem{prop}[thm]{Proposition}

\theoremstyle{definition}
\newtheorem{defi}[thm]{Definition}
\newtheorem{conj}[thm]{Conjecture}
\newtheorem{conv}[thm]{Convention}
\newtheorem{nota}[thm]{Notation}
\newtheorem{rem}[thm]{Remark}
\newtheorem{exa}[thm]{Example}
\newtheorem{sit}[thm]{}

\newcommand{\brem}{\begin{rem}}
\newcommand{\erem}{\end{rem}}
\newcommand{\bexa}{\begin{exa}}
\newcommand{\eexa}{\end{exa}}
\newcommand{\bdefi}{\begin{defi}}
\newcommand{\edefi}{\end{defi}}
\newcommand{\bcor}{\begin{cor}}
\newcommand{\ecor}{\end{cor}}
\newcommand{\blem}{\begin{lem}}
\newcommand{\elem}{\end{lem}}
\newcommand{\bconv}{\begin{conv}}
\newcommand{\econv}{\end{conv}}
\newcommand{\bconj}{\begin{conj}}
\newcommand{\econj}{\end{conj}}
\newcommand{\bprop}{\begin{prop}}
\newcommand{\eprop}{\end{prop}}
\newcommand{\bthm}{\begin{thm}}
\newcommand{\ethm}{\end{thm}}
\newcommand{\bnota}{\begin{nota}}
\newcommand{\enota}{\end{nota}}
\newcommand{\bsit}{\begin{sit}}
\newcommand{\esit}{\end{sit}}

\newcommand{\new}{\newcommand}

\new{\tr}{,\dots,}
\new{\si}{\sigma}
\new{\Si}{\Sigma}
\new{\h}[1]{\hat {#1}}
\new{\bari}[1]{{#1}^{(i)}}
\new{\barun}[1]{{#1}^{(1)}}
\new{\barn}[1]{{#1}^{(n)}}
\new{\quot}[2]{\raise 2pt\hbox{#1}/\raise-4pt\hbox{#2}}
\new{\xb}{\bar x}
\new{\yb}{\bar y}
\new{\zb}{\bar z}
\newenvironment{eg}{\begin{array}{l}}{\end{array}}
\new{\ga}{\gamma}
\new{\ol}{\bar}
\new{\iif}{\Longleftrightarrow}
\new{\lp}{\hbox{\large (}}\new{\pl}{\hbox{\large )}}
\new{\Lp}{\hbox{\Large (}}\new{\pL}{\hbox{\Large )}}
\new{\LP}{\hbox{\LARGE (}}\new{\PL}{\hbox{\LARGE )}}
\new{\hp}{\hbox{\huge  (}}\new{\ph}{\hbox{\huge  )}}

\newcommand{\la}{\label}

\newcommand{\ms}{\medskip}

\newcommand{\sm}{\setminus}
\newcommand{\s}{\sigma}
\newcommand{\re}{resp., }

\newcommand{\hra}{\hookrightarrow}

\newcommand{\PP}{{\mathbb P}}
\newcommand{\R}{{\mathbb R}}
\newcommand{\C}{{\mathbb C}}
\newcommand{\A}{{\mathbb A}}
\newcommand{\Q}{{\mathbb Q}}
\newcommand{\Z}{{\mathbb Z}}
\newcommand{\N}{{\mathbb N}}
\newcommand{\LND}{{\rm LND}}

\newcommand{\nF}{{F^{\rm norm}}}
\newcommand{\nG}{{\G^{\rm norm}}}

\newcommand{\np}{{\pi^{\rm norm}}}

\newcommand{\M}{M_{\sigma}}

\newcommand{\G}{{\Gamma}}

\newcommand{\tH}{{\widetilde H}}
\newcommand{\hE}{{\hat E}}
\newcommand{\hD}{{\hat D}}
\newcommand{\CV}{{C_{\rm vert}}}
\newcommand{\CNV}{{C_{\rm non-vert}}}
\newcommand{\CH}{{C_{\rm horiz}}}
\newcommand{\GH}{{\Gamma_{\rm horiz}}}
\newcommand{\CNH}{{C_{\rm non-horiz}}}
\newcommand{\GNH}{{\Gamma_{\rm non-horiz}}}

\newcommand{\GV}{{\Gamma_{\rm vert}}}
\newcommand{\GS}{{\Gamma_{\rm slant}}}

\newcommand{\CS}{{C_{\rm slant}}}

\newcommand{\GNVs}{{\Gamma^*_{\rm non-vert}}}

\newcommand{\fnh}{{f_{\rm non-horiz}}}

\newcommand{\fh}{{f_{\rm horiz}}}
\newcommand{\fv}{{f_{\rm vert}}}
\newcommand{\fs}{{f_{\rm slant}}}
\newcommand{\SG}{{S_{\Gamma}}}
\newcommand{\SC}{{S_C}}
\newcommand{\HG}{{H_{\Gamma}}}
\newcommand{\HC}{{H_C}}
\newcommand{\hX}{{\hat X}}
\newcommand{\hY}{{\hat Y}}
\newcommand{\p}{{\partial}}
\newcommand{\Aut}{{\operatorname{Aut}}}
\newcommand{\Var}{{\operatorname{Var}}}

\newcommand{\be}{\begin{eqnarray}}
\newcommand{\ee}{\end{eqnarray}}
\newcommand{\no}{\noindent}

\title {Simple birational extensions 
of the  polynomial ring $\C^{[3]}$}

\author{Shulim Kaliman}
\address{Department of Mathematics,
University of Miami,
Coral Gables, FL  33124, U.S.A.}
\email{kaliman@@math.miami.edu}

\author{St\'ephane V\'en\'ereau}
\address{Universit\'e
Grenoble I, Institut Fourier, UMR 5582 CNRS-UJF, BP 74,
38402 St.\ Martin
d'H\`eres c\'edex, France}
\email{venereau@@ujf-grenoble.fr}

\author{Mikhail Zaidenberg}
\address{Universit\'e
Grenoble I, Institut Fourier, UMR 5582 CNRS-UJF, BP 74,
38402 St.\ Martin
d'H\`eres c\'edex, France}
\email{zaidenbe@@ujf-grenoble.fr}

\thanks{Research of the  first author
partially supported by the NSA grant  MDA904-00-1-0016.\\
The third author is grateful to the IHES and to the 
MPI at Bonn (where a part of the work was done) 
for their hospitality and excellent working conditions.\\
\mbox{\hspace{11pt}}{\it 1991 Mathematics Subject Classification}:
14R10, 14R25.\\
\mbox{\hspace{11pt}}{\it Key words}: affine space,
polynomial ring, variable, affine modification, birational extension.}

\begin{document}

\maketitle



\begin{abstract}
The Abhyankar-Sathaye Problem asks whether
any biregular embedding $\varphi:\C^k \hra \C^n$
can be rectified,
that is, whether  there exists an automorphism
$\alpha\in \Aut \,\C^n$ such that
$\alpha\circ\varphi$ is a linear embedding.
  Here we study this problem for the  embeddings
$\varphi:\C^3 \hra \C^4$ whose image $X=\varphi(\C^3)$
is given
in $\C^4$ by an equation  $p=f(x,y)u+g(x,y,z)=0$,
where
$f\in\C[x,y]\backslash\{0\}$ and $g\in\C[x,y,z]$. 
Under certain additional assumptions
we show that, indeed, the polynomial $p$
is a variable of the polynomial ring $\C^{[4]}=\C[x,y,z,u]$
(i.e., a coordinate of a polynomial automorphism of $\C^4$).
This is an analog of a theorem due to Sathaye \cite{Sat1}
which concerns the  case of embeddings $\C^2\hra \C^3$.
Besides, we
generalize a theorem of Miyanishi \cite[Thm. 2]{Miy1} giving, 
for a polynomial $p$ as above, 
a criterion for as when $X=p^{-1}(0)\simeq\C^3$.
\end{abstract}

\vspace{-5pt}

\tableofcontents

\vskip 0.2in

\section*{Introduction}
Generalizing a theorem of A. Sathaye \cite{Sat1} 
it is proven in \cite{KaZa1}
that if a surface $X=p^{-1}(0)\subseteq\C^3$ 
with $p=fu+g\in\C^{[3]}$
and $f,g\in \C[x,y]$ is acyclic 
(that is, $\tH_*(X;\,\Z)=0$)
then $p$ is a variable of the polynomial ring $\C^{[3]}$
i.e., a coordinate of an automorphism 
$\alpha\in\Aut\,\C^3$. Thus $X$ can be rectified,
and so is isomorphic to $\C^2$. 
This does not hold any more in $\C^4$ 
(even with $f\in\C[x]$).
Indeed \cite{Ru, ML1} the Russell cubic $3$-fold 
$$X:=\{p=x^2u+x+y^2+z^3=0\}\subseteq\C^4$$
is an exotic $\C^3$ i.e., 
is diffeomorphic to $\R^6$ 
and non-isomorphic to $\C^3$.
In \ref{c3g}, \ref{mthsec5} and \ref{tra} below 
we give a criterion for as when a $3$-fold 
$X=p^{-1}(0)\subset\C^4$ with 
\query{(label=maineq)}\be\la{maineq}
p=f(x,y)u+g(x,y,z)\in \C^{[4]}\qquad
(f\in \C^{[2]}\backslash\{0\},\,\,g\in\C^{[3]}) \ee
is acyclic (resp., is isomorphic to $\C^3$ resp., 
is an exotic $\C^3$). In particular, we show in \ref{3f} that 
if $X$ is acyclic then actually it is 
diffeomorphic to $\R^6$. If furthermore (\ref{mthsec5})
$X$ is isomorphic to $\C^3$ then any fiber 
$X_{\lambda}:=p^{-1}(\lambda)\,\,\,(\lambda\in\C)$
of the polynomial $p$ is isomorphic to $\C^3$ as well,
and moreover (with an appropriate choice of coordinates 
$(x,y)$) all fibers of the morphism
$$\rho:=(x,p) : \C^4\to\C^2$$
are reduced and isomorphic to $\C^2$. 
We do not know whether in that case a polynomial
$p$ in (\ref{maineq}) must be a variable of
the polynomial ring $\C^{[4]}$, and $\rho$
must be a trivial family. However,
in section \ref{VAR}
in many cases we provide affirmative answers to
these questions and give simple concrete 
examples where the answers remain unknown. 
Due to the Quillen-Suslin Theorem,
the latter question would be answered in positive if the
following conjecture   
\cite[(3.8.5)]{DolVei} (cf. \cite{Sat2, KaZa2}) 
were true for $n=2=\dim S$:  

\smallskip \no {\bf Dolgachev-Weisfeiler Conjecture.} {\it Let $f:X\to S$ 
be a flat affine 
morphism of smooth schemes
with every fiber isomorphic (over the residue field) 
to an affine space $\A^n_k$. 
Then $f$ is locally trivial in the Zariski topology
(i.e., is a fiber bundle).} 

\smallskip 

\no Whereas the former question 
(as whether $p$ is a variable of the polynomial ring $\C^4$)
is a particular case (with $n=4$ and $k=3$) of the famous 

\smallskip \no {\bf 
Abhyankar-Sathaye Embedding Problem:} {\it Is it true that any biregular embedding 
$\C^k\hookrightarrow \C^n$ is rectifiable i.e., 
is equivalent to a linear one under the action 
of the group $\Aut\, \C^n$ on $\C^n$?}

\smallskip

Geometrically, the situation 
can be regarded as follows
(cf. \cite{KaZa1}). The morphism 
$\si: (x,y,z,u)\longmapsto (x,y,z)$ 
represents the $3$-fold $X$ as 
a birational modification of $Y:=\C^3$. 
The latter essentially consists 
in replacing the divisor
$D:=D(f)\subseteq Y$ by the exceptional one
$E:=\{f=g=0\}\subseteq X$.
All the important properties of 
$X$ can be recovered in terms 
of the restriction
$\si|_E : E \to D$. 
In our setting both $E$ and $D$ are
cylindrical surfaces, namely
$E=C\times\C$ and $D=\G\times\C$, where
$C:=D(f)\cap D(g)\subset Y$ is the center
of the blowing up $\si$ and 
$\G:=f^{-1}(0)\subseteq\C^2_{x,y}$.
This makes it possible to formulate 
the criteria mentioned above 
in terms of the natural projection 
$\pi: C \to \G,\,\,(x,y,z)\longmapsto (x,y)$,
and enables us in concrete examples
to verify these criteria. 
 
Let us briefly describe the content of the paper. 
Section 
\ref{sec-prel} contains preliminaries; 
it can be omitted 
at the first reading and consulted 
when necessary. 
However, some results obtained here 
(and used later on in the proofs) 
are of independent interest.
For instance, this concerns \ref{nn'} where we 
treat the question for as when a birational 
extension of a UFD is again a UFD. Furthermore, 
generalizing an observation due to V.
Shpilrain and J.-T. Yu \cite{SY1, SY2} 
we claim
in \ref{isom}-\ref{isomm} 
that for arbitrary polynomials
$p,\,q\in\C[x_1,\ldots,x_n][y]$ 
the hypersurfaces in $\C^{n+1}$ 
given respectively by
the equations
$$y=q(p(y))\qquad\mbox{and}\qquad y=p(q(y))$$
are isomorphic and moreover, 
$1$-stably equivalent
(see \ref{1stab}). 
We also use the fact (see \ref{Alex}) that 
a one-point compactification of 
an acyclic smooth affine variety 
is a homology manifold which is 
a homology sphere and 
satisfies the Alexander duality. 

In section \ref{S2} we study the topology of the 
$3$-folds $X$ as above.
More generally, we work with a 
$3$-fold $X=\{p=fu+g=0\}\subset Y\times\C$,
where $f,g\in\C[Y]$ with $Y$ being 
a smooth acyclic affine $3$-fold
and $D:=f^{-1}(0)\subseteq Y$ being a cylinder
$D=\G\times\C$ over an affine curve $\G$
(whereas in subsection \ref{pracs} 
$Y$ itself is supposed to be
a cylinder over an acyclic 
affine surface $Z$ i.e., 
$Y=Z\times\C$, with $\G\subseteq Z$). 
The main results
of section \ref{S2} (see \ref{3f}, \ref{3g} and \ref{c3g})
provide a criterion for as when such a $3$-fold $X$ is 
diffeomorphic to $\R^6$. 

In subsection \ref{subEXO} we determine when 
$X=p^{-1}(0)\subseteq\C^4$ 
with $p=fu+g$ as in (\ref{maineq}) is 
an exotic $\C^3$. 
The main tool used here is Derksen's 
version of the Makar-Limanov invariant
\cite{ML1, De} described in subsection \ref{subMLD}. 

Subsection \ref{ssec-c3} is devoted to 
a study of embeddings $\C^3\hookrightarrow\C^4$ given
by an equation  $p=fu+g=0$ with 
$f\in\C^{[2]}\backslash\C$ and 
$g\in \C^{[3]}$. In \ref{mthsec5} we show that 
in appropriate new coordinates 
in the $(x,y)$-plane,
the $x$-coordinate restricted to any fiber of $p$
gives a $\C^2$-fibration. On the other hand,
the restriction of $p$ to any hyperplane 
$x= {\rm const}$ is a variable of the polynomial ring
$\C[y,z,u]$ 
(in the latter case 
we say in brief that $p$ is a {\it residual $x$-variable}).

A complete analog of Theorem 7.2 in \cite{KaZa1} 
cited at the beginning holds if the polynomial 
$p\in \C^{[4]}$ is linear with respect to two 
(and not just one) variables. Indeed (\ref{mthsec6})
if the $3$-fold $X$: $$p=a(x,y)u+b(x,y)v+c(x,y)=0$$
in $\C^4$ is smooth and acyclic then $p\in \C^{[4]}$ 
is a variable. We give a simple criterion (in terms of
the coefficients $a,b,c\in\C[x,y]$) for 
as when this is the case. 

In section \ref{VAR} we concentrate on the 
Abhyankar-Sathaye Problem for our particular 
class of embeddings.
The main results \ref{partialres}, \ref{adpr} 
of subsection \ref{ssecPART}
provide sufficient conditions 
for as when a residual
$x$-variable $p=fu+g\in\C^{[4]}$ 
as in (\ref{maineq}) is 
indeed a variable. In \ref{factor} 
we define a canonical factorization 
$f=\fh\fs\fv$, and we show (\ref{partialres})
that an embedding 
$\C^3\stackrel{\simeq}{\longrightarrow} 
X:=p^{-1}(0)\subseteq\C^4$
can be rectified if either the polynomial 
$\fh$ is reduced, or
$\fv=1$. Moreover, in appropriate new coordinates 
$(x,y)$
we have $\fnh:=\fs\fv\in\C[x]$, 
and $p$ is an $x$-variable 
i.e.,
a variable in $\C[x]^{[3]}$ over $\C[x]$.
For instance 
(see \ref{partialres}, \ref{cova}) 
this is the case
if $\deg_z g\le 1$, or $f$ 
is a power of an irreducible polynomial,
or else $f\in\C[x]$ (the latter strengthens
a result of M. Miyanishi \cite[Thm. 2]{Miy1}, 
where it is supposed in addition that $g\in\C[y,z]$).
As another examples, we show (see \ref{x5})
that the polynomials 
$$p_1:=xy^2u+y+xz+xyz^2\qquad {\rm and}\qquad 
p_3:=xy^2u+y+x^2z+x^3yz^2$$ are variables of $\C^{[4]}$.
However, we do not know whether or not so is 
$$p_2:=xy^2u+y+x^2z+xyz^2\,$$ 
(whereas $p_2^{-1}(\mu)\simeq\C^3\,\,\,\forall\mu\in\C$
and
$p_2$ is  a residual
$x$-variable and a $\C(x)$-variable, 
see  \ref{cform} and \ref{QN}). 

In subsection \ref{1stably}  
we establish (see \ref{x-planes}) that every embedding
$\C^3\hookrightarrow\C^4$ given by an equation 
$p=fu+g=0$ as in (\ref{maineq}) can be rectified in $\C^5$. 

In the last subsection \ref{SW} 
(attributed to the second author)
we generalize a theorem of D. Wright
\cite{Wr} which says that Sathaye's Theorem holds for
the embeddings $\C^2\hookrightarrow\C^3$ 
given by an equation $p=fu^n+g=0$ with $f,g\in\C[x,y]$
and $n\in\N$. Namely, 
it is shown in \ref{Wr} that a residual
$x$-variable of the form $p=fu^n+g\in\C^{[4]}$, 
where $f,g\in\C[x,y,z]$ and $n\ge 2$, 
actually is an $x$-variable.


\query{(label=sec-prel)}\section{Preliminaries}\la{sec-prel}
\query{(label=sec-UFD)}\subsection{Affine modifications of
UFD's}\la{sec-UFD}

{We} start by recalling the notion of affine modification
\cite{KaZa1}; at the same time, we introduce the notation 
that will be used throughout the  paper.

\query{(label=mm)}\no \bnota\la{mm} {\rm  
Let $Y$ be a reduced, irreducible affine variety
over $\C$, $A=\C[Y]$ be the  algebra of regular functions on $Y$,
$I\subseteq A$ be a non-trivial ideal and $f\in I$ be a non-zero
element of $I$. {\it The affine modification of the  variety $Y$
along the  divisor $D_f=f^*(0)$ with center} $I$ is the  affine
variety $X=\rm{spec}\, A'$, where $$A':=A[I/f]=\{a'={a_k /f^k}|a_k \in I^k\}\subseteq {\rm Frac}\,A\,.$$
The inclusion
$A\hra A'$ corresponds to a birational morphism $\s:X\to Y$
with the  exceptional divisor 
$E=\s^{-1}(D)=(f\circ\s)^{-1}(0)\subseteq X$,
where $D:={\rm supp}\, D_f$. 
The restriction $\s|(X\backslash E): X\backslash E\to Y\backslash D$ is an
isomorphism.

Let $D=\bigcup_{i=1}^{n} D_i $
resp., $E=\bigcup_{i=1}^{n'} E_i $ be the decomposition into
irreducible components, which we assume to be Cartier divisors. 
Letting
\begin{eqnarray}   \s^*(D_i) = \sum_{j=1}^{n'} 
m_{ij} E_j,\,\,\,\,i=1,\dots,n\,,
\end{eqnarray}
we consider the $n\times n'$ {\it multiplicity matrix} $\M=(m_{ij})$
with non-negative integer entries. Clearly, $m_{ij}>0\Leftrightarrow
\si(E_j)\subseteq D_i$.} \enota

The following simple observation will be useful in \ref{pi1} below. 
We denote
${\rm reg}\,E_j=E_j\backslash  {\rm sing}\, E_j$ and 
${\rm reg}\,D_i=D_i\backslash  {\rm sing}\, D_i$.

\query{(label=gone)}\no\blem\la{gone} In the notation 
as above, suppose that the affine varieties $X$ and $Y$ are smooth. 
If $m_{ij}=1$ then $\si(E_j)\not\subseteq {\rm sing}\, D_i$. 
Moreover $m_{ij}=1$ if and only if
$\si({\rm reg}\,E_j)\subseteq {\rm reg}\,D_i$ and $\si$ 
sends the analytic discs in $X$ transversal to $E_j$ at a point
$Q\in {\rm reg}\,E_j$ biholomorphically onto analytic discs in 
$Y$ transversal to $D_i$ at the point $P:=\si(Q)\in {\rm reg}\,D_i$. 
\elem

\no\proof We may assume that $m_{ij}>0$ that is, $\si(E_j)\subseteq D_i$. 
For a point 
$Q\in E_j$ with the image $P=\si(Q)\in D_i$,
we let $U\ni P$ resp., $V\ni Q$ 
be a neighborhood such that $D_i\cap U=f_i^*(0)$ resp., 
$E_j\cap V=h_j^*(0)$, where $f_i$ resp., $h_j$ 
is a holomorphic function in $U$ resp., $V$. 
Then we have
$$f_i\circ\si|_V=h_j^{m_{ij}}\cdot h\qquad 
{\rm with}\quad h(Q)\neq 0\,,$$
which gives the equality of $1$-forms: 
\query{(label=cir)}\be\la{cir} 
d(f\circ\si)(Q)=m_{ij}h_j^{m_{ij}-1}(Q)h(Q)\cdot dh_j(Q)\,.\ee
Assuming that $Q\in {\rm reg}\,E_j$ 
($\Leftrightarrow dh_j(Q)\neq 0$) and $m_{ij}=1$, from (\ref{cir})
we obtain:
$$d(f\circ\si)(Q)=h(Q)\cdot dh_j(Q)\neq 0\,,$$
whence $df_i(P)\neq 0$ and $d\si(T_QX)\not\subseteq T_P D_i$.
This yields the implication "$\Longrightarrow$". On the other hand,
if $P\in {\rm reg}\, D_i$ ($\Leftrightarrow df_i(P)\neq 0$) and $d\si(T_QX)\not\subseteq T_P D_i$ then
$d(f\circ\si)(Q)\neq 0$, and so by (\ref{cir}) $m_{ij}=1$,
which gives "$\Longleftarrow$". \qed\medskip

 For an algebra $B$, denote by $B^*$ its group of invertible elements.

\query{(label=nn')}\no \bprop \la{nn'}
In the notation as in \ref{mm} above, 
assume moreover that the algebras $A=\C[Y]$
 and $A'=\C[X]$ are UFD's, and $ A^*={A'}^*$.
 Then $n=n'$ and the  multiplicity matrix $\M$ is unimodular.
\eprop

 \noindent \proof Since the  algebras $A$ and $A'$ are UFD's there exist
 irreducible elements $f_1,\dots,f_n\in A$ resp.,
$h_1,\dots, h_{n'}\in A'$ such that
 $D_i=f_i^* (0),\,\,\,i=1,\ldots,n,$ resp.,
$E_j=h_j^*(0),\,\,\,j=1,\ldots, n'.$
 Clearly,   $h_j=(a_j/b_j)\circ\s$
 with coprime $a_j,\,b_j\in A, \,\,j=1,\ldots,n'$.
 The regular functions $a_j,\,b_j$ do not vanish in
 $Y\backslash D$, and hence their irreducible factors
 are proportional to some of the  elements
 $ f_1, \ldots, f_n$.  Thus
 for each $j=1,\dots,n'$
 there exists $\gamma_j\in A^*$ such that
 \begin{eqnarray} \la{mij} h_j=(\gamma_j\prod_{i=1}^n
 f_i^{m_{ji}'})\circ\s\,
 \end{eqnarray}\query{(label=mij)}
 with $m_{ji}'\in \Z$. On the  other hand,
 for each $i=1,\dots,n$ there exists  $\delta_i\in {A'}^{*}=A^*$
 such that
 \begin{eqnarray} \la{mjk} f_i\circ\s=\delta_i\prod_{k=1}^{n'} h_k^{m_{ik}}\,.
 \end{eqnarray}\query{(label=mjk)}
   Plugging successively (\ref{mij}) and (\ref{mjk}) one into another
  and taking into account the  assumption
 that the algebras $A'$ and $A$ are UFD's and $\s$ is a birational morphism,
 we obtain the  equalities
 \begin{eqnarray}
   \la{id}
   \sum_{k=1}^{n'}m_{ik}m'_{kj}=\delta_{ij},\,\,\,i,j=1,\dots,n\,,\\
   \sum_{i=1}^nm_{ji}'m_{ik}=\delta_{jk},\,\,\,j,k=1,\dots,n'\,,
 \end{eqnarray}\query{(label=id)}
 that is, $\M'\cdot\M=I_{n'}$ and  $\M\cdot\M'=I_{n}$,
 where  $I_{k}$ denotes
 the  identity matrix of order $k$.
 Hence $n=n'$ and  $\M'=\M^{-1}$, so $\M$ is unimodular. 
 \qed\medskip

\query{(label=pi1)}\no \bprop\la{pi1} 
Suppose that the varieties $X$ and $Y$ as in \ref{mm} are smooth,
and the  multiplicity matrix $\M$ is an upper triangular unipotent square
matrix: $$n=n',\quad m_{ii}=1,\,\,i=1,\ldots,n,
\qquad m_{ij}=0\,\,\,\forall i>j\,.$$ Then 
$\s_*:\pi_1(X)\to\pi_1(Y)$ is an isomorphism.\eprop

\no\proof Letting $G=\pi_1(Y\backslash  D),\,\,\,\,
G'=\pi_1(X\backslash  E)$ and $\si'=\si|_{X\backslash  E}$, we have that
$\si'_* : G'\to G$ is an isomorphism which sends the subgroup $$H':=\big\langle\big\langle\alpha_{E_1},\ldots,\alpha_{E_n}\big\rangle\big\rangle\,\,\subseteq G'$$
into the subgroup
$$H:=\big\langle\big\langle\alpha_{D_1},
\ldots,\alpha_{D_n}\big\rangle\big\rangle\,\,\subseteq G\,,$$
where for a hypersurface $Z$ in a complex manifold  $X$,
$\alpha_Z$ denotes a vanishing loop of $Z$ in $X\setminus Z$, 
whereas for a group $G$ and elements $a_1,\ldots,a_k\in G$, 
$\big\langle\big\langle a_1,\ldots,a_k\big\rangle\big\rangle$ denotes the minimal 
normal subgroup of $G$ generated by 
$a_1,\ldots,a_k$. Moreover, 
$\si_* :\pi_1(X)\to\pi_1(Y)$ is a surjection
with kernel ${\rm ker}\,\si_*\cong H/\si'_*(H')$ 
(see e.g., \cite[the proof of Prop. 3.1]{KaZa1}). 
Thus we must show that $\si'_*(H')=H$.

For $i=1,\ldots,n$ denote $$H_i:=\big\langle\big\langle\alpha_{D_1},
\ldots,\alpha_{D_i}\big\rangle\big\rangle
\,\,\subseteq G\quad
{\rm resp.,}\quad H'_i:=\big\langle\big\langle\alpha_{E_1},
\ldots,\alpha_{E_i}\big\rangle\big\rangle\,\,\subseteq G'\,,$$
so that $H_n=H$ and $H'_n=H'$.
Clearly, $\si'_*(\alpha_{E_1})\sim \alpha_{D_1}$ 
(where $\sim$ 
stands for conjugation), whence $\si'_*(H'_1)=H_1$. 
We show by induction that $\si'_*(H'_i)=H_i$ 
for all $i=1,\ldots,n$.

Assume that this equality holds for $i=k-1<n$. 
As $m_{kk}=1$, by 
\ref{gone} we may conclude that for a general point 
$Q\in E_k$, $P:=\si(Q)$ is a smooth point 
of $D_k$, and $\si$ sends biholomorphically
a smooth analytic disc transversal to the divisor $E_k$ at $Q$ 
onto a transversal disc
to $D_k$ at $P$. As the matrix $\M$ is upper 
triangular
we obtain $\si'_*(H'_i)\subseteq H_i\,\,\,(i=1,\ldots,n)$ 
and furthermore, $\si'_*(\alpha_{E_k})\sim 
\alpha_{D_k}\,\mod\,H_{k-1}$.
As $\si'_*(H'_{k-1})=H_{k-1}$ it follows that also 
$\si'_*(H'_{k})=H_{k}$, therefore $\si'_*(H')=H$, 
as desired. \qed\medskip

\query{label=acvar}\subsection{Acyclic varieties}\label{acvar}
Recall that the acyclicity of a topological space 
$X$ means that
 its reduced homology vanishes: 
${\widetilde H}_* (X):={\widetilde H}_* (X; \Z)=0$. 
The following proposition is an immediate corollary of 
the above results. 

\query{(label=new)} \no \bprop\la{new}
Assume that the affine varieties $X$ and $Y$ as in 
\ref{mm} above are smooth and acyclic. Then
\begin{enumerate} 
\item[(a)]
$n=n'$ and the  multiplicity matrix $\M$ is unimodular.
\end{enumerate}
If moreover, $\M$ is an upper triangular unipotent
matrix then
\begin{enumerate}
\item[(b)] $\s_*:\pi_1(X)\to\pi_1(Y)$ is an isomorphism.
\item[(c)] $X$ is contractible if and only if so is $Y$.
\item[(d)] In the latter case the varieties $X$ and $Y$
are both diffeomorphic to the affine space
$\R^{2m}$ provided that $m:=\dim_{\C} Y\ge  3$.
\end{enumerate}
 \eprop

 \noindent \proof By \cite[1.18-1.20]{Fuj} 
(see also \cite[3.2]{Ka1})
 for any smooth acyclic affine variety $Z$ 
the algebra $\C[Z]$ is UFD and  its
 invertible elements are constants. Thus
(a) follows from \ref{nn'}.
(b) directly follows from \ref{pi1}. 
In virtue of the Hurewicz and Whitehead theorems
\cite[Ch.2, \S 11.5, \S 14.2]{FoFu},
an acyclic manifold is contractible 
if and only if it is simply connected.
Hence (c) follows from (b). In turn, 
(d) follows from the
Dimca-Ramanujam  theorem \cite{ChDi, Za2}.
\qed\medskip

\query{(label=ssy)}\no\bsit\la{ssy}{\rm
In section \ref{ssec-c3} we will apply the following 
corollary (see \ref{miy} below) of Miyanishi's characterization of 
$\C^3$ \cite{Miy1, Miy2}. On the other hand, 
this corollary also follows from \cite{Sat2} and \cite{KaZa2},
as stated in \cite[Cor. 0.2]{KaZa2}. Moreover, 
by \cite[L. III]{Ka2}
it would be enough to suppose in \ref{miy} that only general 
(rather than all) fibers of $x$ were isomorphic to $\C^2$. 
}\esit

\query{(label=miy)}\no\bthm\label{miy} Let $X$ 
be a smooth acyclic affine 3-fold. Then 
$X\simeq \C^3$ if and only if there exists 
a regular function 
$x\in \C[X]$ with all fibers isomorphic to $\C^2$. \ethm

\query{(label=bimo)}\bsit\la{bimo} 
{\rm It is known \cite[Thm.1.1]{KaZa1} that 
any birational morphism $\si:X\to Y$ 
of affine varieties 
is an affine modification. The divisor $D\subseteq Y$ 
of modification and the exceptional divisor 
$E \subseteq X$ of $\si$ can be defined 
as minimal reduced divisors
such that the restriction 
$\si|_{X\backslash  E} : X\backslash  
E\to Y\backslash  D$ 
is an isomorphism. In the next proposition 
and its corollary 
we provide
conditions 
(more general than those in \cite[Thm. 3.1]{KaZa1}) 
which guarantee preservation of the 
homology group under modification. 
(Note that the divisors $\hE$ and $\hD$ 
below do not need to satisfy 
the assumption of minimality.)} \esit

\query{(label=prhom)}\no\bprop\label{prhom}
Given affine varieties ${\hat X},\,{\hat Y}$ 
and decompositions
\query{partitions?} 
$${\hat X}=({\hat X}\backslash  \hE)\cup \hE\quad {\rm and}
\quad {\hat Y}=({\hat Y}\backslash  \hD)\cup \hD,$$
where $\hE \subseteq {\hat X}$ and $\hD\subseteq {\hat Y}$ 
are reduced divisors, let $\si:{\hat X}\to {\hat Y}$ be a birational 
\query{Do we really need birational?}
morphism  
which respects these decompositions. 
Suppose that the following hold.
\begin{enumerate}
\item[(i)] $\hE,\hD$ are topological manifolds, 
$\hE \subseteq {\rm reg}\, {\hat X},\,\,\hD\subseteq {\rm reg}\,{\hat Y}$, 
and $\si^*(\hD)=\hE$.
\item[(ii)] The induced homomorphisms
$$(\si|_{\hE})_* : H_*(\hE)\to 
H_*(\hD)\quad {\rm
and}\quad (\si|_{{\hat X}\backslash  \hE})_* : H_*({\hat X}\backslash  \hE)
\to H_*({\hat Y}\backslash  \hD)$$ are
isomorphisms. 
\end{enumerate}
Then $\si_* : H_*({\hat X})
\to H_*({\hat Y})$ is
an isomorphism as well; in particular, ${\hat X}$ 
is acyclic if and only if ${\hat Y}$ is so. \eprop

\no\proof We apply the Thom isomorphism  
\cite[\query{p. 297,} 7.15]{Do}
to the pairs $({\hat X},\hE)$ resp., $({\hat Y},\hD)$
(notice that locally near $\hE$
resp., $\hD$ these are pairs of topological manifolds). 
As $\si^*(\hD)=\hE$ and 
$(\si|_{\hE})_* : H_0(\hE)\stackrel{\cong}{\longrightarrow} H_0(\hD)$, 
$\si$ maps the irreducible components of $\hE$
into those of $\hD$ providing a one-to-one correspondence, 
and (as in \ref{gone}) 
sends their transverse classes 
\cite[Ch. VIII\query{, p.315}]{Do}  
to the corresponding transverse classes. 
By  functoriality of the cap-product \cite[\query{p.239,} VII.12.6]{Do}
for every $i\geq 0$
the following diagram is commutative:\\
\begin{diagram}[height=2.3em,width=3em]
H_i({\hat X},\,{\hat X}\backslash \hE)&\rTo^{\cong}_{t_{\hat X}}&\tH_{i-2}(\hE)\\
\dTo_{\si_*} & &\dTo_{\si_*}& & &                           \\
H_i({\hat Y},\,{\hat Y}\backslash \hD)&\rTo^{\cong}_{t_{\hat Y}}& \tH_{i-2}(\hD)\\
\end{diagram}

\no (where $t$ stands for the Thom isomorphism, 
and the homology groups
in negative dimensions are zero). 
This allows to replace the relative homology groups in the 
exact sequences of pairs as to obtain the following 
commutative diagram: 

\begin{diagram}[height=2.2em,width=3.2em]
\ldots  \to & \tH_{i-1}(\hE) & \longrightarrow & \tH_i({\hat X}\backslash \hE) & \longrightarrow & 
\tH_{i}({\hat X}) & \longrightarrow & \tH_{i-2}(\hE) &\longrightarrow & \tH_{i-1}({\hat X}\backslash \hE) & \to \,\ldots \\
 & {\wr}\,\dTo_{\si_*} & & {\wr}\,\dTo_{\si_*} & & \dTo_{\si_*} & & {\wr}\,\dTo_{\si_*} & & {\wr}\,\dTo_{\si_*} & & \\
\ldots  \to & \tH_{i-1}(\hD) & \longrightarrow & \tH_i({\hat Y}\backslash \hD) & \longrightarrow &
\tH_{i}({\hat Y}) & \longrightarrow & \tH_{i-2}(\hD) &\longrightarrow & \tH_{i-1}({\hat Y}\backslash \hD) & \to \,\ldots \\
\end{diagram}\\

\no By (ii) the four vertical arrows (as shown at the diagram) 
are isomorphisms, whence
by the 5-lemma, the middle one is so as well, as stated.
\qed\medskip

Actually, the divisors $E$ and $D$ we 
deal with in section \ref{S2} below are not always topological manifolds. 
However, in our setting we can apply \ref{prhom} by 
decomposing further as follows.

\query{label=prdec}\no\bcor\label{prdec}
The same conclusion as in \ref{prhom} 
holds if (instead of (i), (ii)) we assume that
there are decompositions $\hE=E'+E''$ and $\hD=D'+D''$
satisfying the following conditions:
\begin{enumerate}
\item[(i$'$)] $E'\backslash  E'',\,D'\backslash  D'',\,E''$ and
$ D''$
are topological manifolds, 
$\hE \subseteq {\rm reg}\, {\hat X},\,\,\hD\subseteq {\rm reg}\,{\hat Y}$, 
and 
$\si^*(D'\backslash  D'')=E'\backslash  E'', \,\,\si^*(D'')=E''$.
\item[(ii$'$)] The induced homomorphisms 
$\,\,\,(\si_{{\hat X}\backslash  \hE})_*: H_*({\hat X}\backslash  \hE)
\to H_*({\hat Y}\backslash  \hD)$,
$$
(\si|_{E'\backslash  E''})_* : 
H_*(E'\backslash  E'')\to 
H_*(D'\backslash  D'')
\quad {\rm and}\quad 
(\si|_{E''})_* : H_*(E'')
\to H_*(D'')
$$ 
are isomorphisms. 
\end{enumerate}\ecor

\no\proof Indeed, \ref{prhom} implies that 
(under our assumptions)
$$(\si|_{{\hat X}\backslash  E''})_* : 
H_*({\hat X}\backslash  E'')\to 
H_*({\hat Y}\backslash  D'')$$ is an isomorphism. Now (with
$\hE$ and $\hD$ in \ref{prhom} replaced by $E''$ resp., $D''$),
\ref{prhom} implies that  
$\si_*: H_*({\hat X})\to 
H_*({\hat Y})$ is an isomorphism as well. \qed\medskip

\query{(label=homa)}\bsit\la{homa} 
{\rm Recall \cite[Ch. 5, \S 5.3]{Mau}, 
\cite[Ch. 2, \S 8.1]{FuVi}\footnote{We are grateful 
L. Guillou for useful discussions on homology manifolds
and references.} 
that a simplicial polyhedron $P$ is called a
{\it homology $n$-manifold} if for any point $p\in P$ we have 
$H_*(P,P\backslash  \{p\})\cong \tH_*(S^n)$ \footnote{Or equivalently, 
$H_*(lk(p))\cong H_*(S^{n-1})$, where $lk(p)$ denotes the link of 
$p$ in $P$, or else, for any $q$-simplex $\si$ in $P$, 
$H_*(lk(\si))\cong H_*(S^{n-q-1})$.}. }\esit

\query{label=Alex}\no\bprop\label{Alex} 
Let $X$ be a smooth acyclic affine variety of
complex dimension $n$. Then the one-point compactification ${\dot
X}$ of $X$ is a homology $2n$-manifold which is a homology
$2n$-sphere: $H_*({\dot X})\cong H_*(S^{2n})$.
In particular, the Alexander duality holds for ${\dot X}$. \eprop

\no\proof Notice first that $X$ is diffeomorphic to the
interior of a compact manifold, say, $X_R$ with boundary
$\partial X_R$ (indeed, one can take for $X_R$ the
intersection of $X$ with a ball of a large enough radius $R$ in
an affine space $\C^N\supseteq X$). Hence ${\dot X}\simeq
X_R/\partial X_R\simeq X_R\cup_{\partial X_R} C(\partial X_R)$ (with $CY$ denoting
the cone over $Y$), and any triangulation of $X_R$ 
naturally extends to those of $\dot X$. 

By the Poincar\'e-Lefschetz duality for a
manifold with boundary \cite[5.4.13]{Mau}, 
$$H^{2n-i}(X_R)\cong H_i (X_R,\partial
X_R)\cong {\widetilde H}_i (X_R/\partial X_R)\cong 
{\widetilde H}_i({\dot X})\,,$$ 
whence
${\dot X}$ is a homology $2n$-sphere. Using the acyclicity of
$X_R$ and of $C(\partial X_R)$ and applying the Mayer-Vietoris
sequence to the decomposition 
${\dot X}= X_R \cup_{\partial X_R} C(\partial X_R)$,
we see that $H_i({\dot X})\cong H_{i-1}(\partial X_R)$. Thus the
smooth manifold $\partial X_R$ is a homology ($2n-1$)-sphere.

Clearly,  
$H_*({\dot X},{\dot X}\backslash  \{x\})\cong \tH_*(S^{2n})$
for any point $x\in X={\dot X}\backslash \{\infty\}$. 
For the vertex $x=\infty$ of the cone $C(\partial X_R)$ and for 
any $i\in \N$, 
by excision and from the exact homology sequence of a pair
we obtain: 
$$H_i({\dot X},{\dot X}\backslash  \{\infty\})\cong 
H_i(C(\partial X_R),C(\partial X_R)\backslash  \{\infty\})$$ $$
\cong \tH_{i-1}(\partial X_R)\cong \tH_{i-1}(S^{2n-1})
\cong \tH_i(S^{2n})\,.$$
Thereby ${\dot X}$ is a homology manifold 
and is a homology $2n$-sphere.
For any point 
$x\in {\dot X}$, 
from the exact 
homology sequence of the pair 
$({\dot X}, {\dot X}\backslash  \{x\})$
it follows that ${\dot X}\backslash  \{x\}$ is acyclic. 
Now the proof of  
the Alexander duality for the usual sphere \cite[5.3.19]{Mau}
goes {\it mutatis mutandis} for ${\dot X}$ (cf.
\cite[Thm. 6.4]{Be}, \cite[p. 176]{Wi}). \qed\medskip

\query{(label=subMLD)}
\subsection{Digest on Makar-Limanov and 
Derksen invariants}
\la{subMLD}
These invariants (introduced in \cite{ML1, KaML1, De}) allow
in certain cases to distinguish space-like
affine varieties from the affine spaces. On this purpose, 
we use them in subsection \ref{subEXO} (to establish \ref{tra}). 
Let us first recall the following notions and facts. 

\query{(label=LND)}\no\bsit\la{LND} 
{\it Locally nilpotent derivations.}
{\rm 
Let $A$ be an affine domain over $\C$.
A derivation $\partial\in {\rm Der}\,A$
of $A$ is called {\it locally nilpotent}
({\it LND} for short) if for each $a\in A$ there exists
$n=n(a,\,\p)\in\N$ such that $\p^{(n)}(a)=0$;
the set of all non-zero locally nilpotent
derivations of the algebra $A$ is denoted by LND$(A)$.
Given $\p\in {\rm LND}\,(A)$, 
the function $\deg_{\p}(a):=\min \{n(a,\,\p)-1\}$ 
is a degree function
on $A$. 
The kernel $A^{\p}={\rm ker}\,\p$ of $\p$ is a
$\p$-invariant subalgebra of $A$; 
its elements are called {\it $\p$-constants}. 
For the proof of the following lemma see e.g., 
\cite{ML1, ML2, KaML1, De}, 
\cite[\S 7]{Za2} or \cite[5.1(6)]{Ka2}.
}\esit

\query{(label=trdeg)}\no\blem\la{trdeg}    
The following statements hold:

\begin{enumerate}
\item[(a)] tr.deg $[A:A^{\p}]=1$.
\item[(b)] The subalgebra $A^{\p}$ is algebraically closed.
\item[(c)] It is factorially closed i.e.,
$uv\in A^{\p}\backslash \{0\}\Longrightarrow u,\,v\in A^{\p}$.
\item[(d)] If $u^k+v^l\in A^{\p}\backslash \{0\}$
for some $k,\,l\ge 2$ then $u,\,v\in A^{\p}$.\\ 
Moreover,
$p(u,v)\in A^{\p}\backslash \C\Longrightarrow u,\,v\in A^{\p}$ 
for any polynomial $p\in\C^{[2]}$ 
with general fibers being irreducible and 
non-isomorphic to $\C$.
\end{enumerate}
\elem

\query{(label=INV)}\no\bsit\la{INV}  {\it  Invariants.} {\rm
The {\it Makar-Limanov invariant}
of the algebra $A$ is the subalgebra
$${\rm ML}\,(A)=\bigcap_{\p\in {\rm LND}(A)} A^{\p}\subset A\,,$$
whereas the  {\it Derksen invariant}:
$${\rm Dk}\,(A)=\C\Big\lbrack \bigcup_{\p\in {\rm LND}(A)}
A^{\p}\Big\rbrack\subset A$$
is the subalgebra of $A$ generated by the $\p$-constants of
all locally nilpotent derivations on $A$. If 
ML$(A)=\C$ (\re Dk$(A)=A$) then we
say that
the corresponding invariant is trivial. 
This is, indeed, the case
for a polynomial algebra $A=\C^{[n]}\,\,(n\in\N)$.
}\esit

\query{(label=SPE)}\no\bsit\la{SPE}  {\it  Specializations.} {\rm
Let $X$ be an affine variety, and set $A=\C[X]$.
To study the locally nilpotent derivations
on $A$, it is possible
to proceed by induction on the dimension of $X$. 
Namely, let $\p\in {\rm LND}\,(A)$, and
let  $u\in {\rm ker}\,\p$
be non-zero. As 
$\,\p (u-c)=0\,\,\forall c\in\C$, 
the principal ideal $(u-c)$ 
of the algebra
$A$ is invariant under $\p$, and so $\p$ descends to
the quotient $B_c=A/(u-c)=\C[S_c]$, where 
$S_c=u^{-1}(c)$ is 
a fiber of $u$. For a general $c\in\C$,
this specialization $\p_c$ is a non-zero locally nilpotent
derivation of the algebra $B_c$ (see \ref{ACT} below).
Clearly, the restriction to $S_c$ of any $\p$-constant
$v\in {\rm ker}\,\p$ is a $\p_c$-constant.
}\esit

\query{(label=ACT)}\no\bsit\la{ACT}  {\it $\C_+$-actions.} {\rm
Otherwise, the above specialization can be described 
via the natural
correspondence between the locally nilpotent derivations
of the algebra $A$ and the regular actions of the additive group
$\C_+$ on the variety $X= {\rm spec}\,A$ (e.g., see \cite{Re, Za2}).
Indeed, the subalgebra ${\rm ker}\,\p$ 
coincides with the algebra of invariants $A^{\varphi}$
of the associated $\C_+$-action ${\varphi}={\varphi}_{\p}$.
If $u$ is a ${\varphi}$-invariant then
clearly, the $\C_+$-action ${\varphi}|_{S_c}$ 
on a fiber $S_c={\rm spec}\,B_c$
of $u$
is associated with the above specialization $\p_c$. 
Hence $\p_c\in {\rm LND}\,(B_c)$ 
if and only if the $\C_+$-action ${\varphi}|_{S_c}$
is non-trivial. 
}\esit

\query{(label=JD)}\no\bsit\la{JD}  {\it 
Jacobian derivations \cite{KaML1, KaML2}.} 
{\rm For an $n$-tuple of polynomials
$p_1,\dots,p_{n-1},\,q\in\C^{[n]}$, the Jacobian
$$\p(q):={\rm jac}\,(p_1,\dots,p_{n-1},\,q)=
{\p(p_1,\dots, p_{n-1},\,q)\over \p(x_1,\dots,x_n)}$$
(regarded as a function of $q$, whereas the 
polynomials $p_1,\dots,p_{n-1}$ are fixed) gives 
a derivation on the polynomial algebra
$\C^{[n]}$. This derivation is non-zero provided that  the
polynomials $p_1,\dots,p_{n-1}$ are algebraically independent.
For $p:=p_1$, the principal ideal $(p)\subseteq \C^{[n]}$ 
is invariant under $\p$, whence $\p$ descents to
a derivation of the quotient algebra $A:=\C^{[n]}/(p)=\C[X]$,
where $X:=p^*(0)\subset \C^{n}$ \footnote{To 
simplify the exposition we suppose that $X$ is a hypersurface,
whereas the results of \cite{KaML2}
hold in a more general setting.}.
Denote  by JLND$(A)$ the set of all locally nilpotent
Jacobian derivations
of the algebra $A$. Two derivations
$\p,\,\p'\in {\rm Der}\,(A)$ are called
{\it equivalent} if $A^{\p}=A^{\p'}$. Actually, two derivations
are equivalent iff they generate the same degree function 
\cite{KaML2}: 
$\deg_{\p}=\deg_{\p'}$. We have the following
theorem.
}\esit

\query{(label=KML1)}\no\bthm\la{KML1} 
{\rm (\cite{KaML1, KaML2})}
If $p\in \C^{[n]}$ is a non-constant irreducible
polynomial then
every locally nilpotent derivation 
$\p\in {\rm LND}\,(A)$ (where as above, 
$A=\C^{[n]}/(p)$) is equivalent
to a Jacobian locally nilpotent
derivation $\p'\in {\rm JLND}\,(A)$.
Henceforth, we have \ethm
$${\rm ML}\,(A)=\bigcap_{\p\in
{\rm JLND}(A)} A^{\p}\qquad
{\rm and}\qquad {\rm Dk}\,(A)=
\C\Big\lbrack\bigcup_{\p\in 
{\rm JLND}(A)} A^{\p}\Big\rbrack\,.$$

\query{(label=WDF)}\no\bsit\la{WDF}  {\it 
Weight degree functions and the associated graded algebras} 
{\rm (see e.g., \cite{KaML1, KaML2, Za2}). 
A {\it weight degree function}
on the algebra $A=\C[X]=\C^{[n]}/(p)$ is defined by 
assigning weights $d_i=d(x_i)\in \R$ to the variables
$x_1,\dots,x_n\in\C^{[n]}$ and letting for $a\in A$:
$$d_A(a):=\inf\,\{d(q)\,\big|\,q\in \C^{[n]},\,q|_X=a\}\,.$$
Here as usual
$$d(q)=\max\,\{d(m)\,|\,m\in M(q)\}\,$$
with 
$M(q)$ denoting the set 
of the monomials $m=c_m\prod_{i=1}^n x_i^{\alpha_i}$
of $q$, where $d(m):=\sum_{i=1}^n \alpha_i d_i$.
The  weight degree function $d_A:A\to \R\cup\{-\infty\}$
defines an ascending filtration 
$\{A_t\}_{t\in \R}$ of the algebra $A$
with
$A_t:=\{a\in A\,|\,d_A(a)\le t\}$.
We also let  $A_t':=\{a\in A\,|\,d_A(a)< t\}\subset A_t$, and
we consider the {\it associated graded algebra}
$$\h A=\bigoplus_{t\in \R} A_t/A_t'$$
(actually, the set of non-zero
homogeneous components $ A_t/A_t'$ of the algebra $\h A$
is at most countable).
}\esit

\query{(label=AHD)}\no\bsit\la{AHD}  {\it 
Associated homogeneous derivations. }
{\rm 
For a polynomial $q\in \C^{[n]}$, we consider its
{\it $d$-principal part} (in other words,
the {\it principal $d$-quasihomogeneous part})
$\h q:=\sum_{m\in M(q),\,d(m)=d(q)} m$. For an element
$a\in A$, we let $\h a$ to be its image in the graded
algebra $\h A$ (clearly, $\h a\in A_t/A_t'$ with $t=d_A(a)$).
Notice that $\h a=\h q\,|_{\h X}$
for a polynomial $q\in\C^{[n]}$ such that
$q\,|_{X}=a$ and $d(q)=d_A(a)$; the latter equality 
holds if and only if
$\h q|_{\h X}\neq 0$.

If $\p\in {\rm Der}\,(A)$ is a derivation then the degree
$$d_A(\p) := \inf\,\{d_A(a)-d_A(\p a)\,|\,a\in A\}$$
is finite \cite{KaML2}. Letting
$$\h\p\h a:={\widehat {\p a}} \quad {\rm if} 
\quad d_A(a)-d_A(\p a)=d_A(\p)
\quad {\rm and} \quad \h\p\h a=0 \quad {\rm otherwise,}\,$$
and extending $\h\p$ in a natural way to a homogeneous
derivation of the graded algebra $\h A$, 
we obtain a correspondence
Der$(A)\to $Der$_{\rm gr} (\h A)$. 
It has the following properties.
}\esit

\query{(label=KML2)}\no\bthm\la{KML2}
{\rm (\cite{KaML1, KaML2})}
Let  
$p\in \C^{[n]}$ be a polynomial with a non-constant, 
irreducible $d$-principal part $\h p$.
Set $\h X=\h p^{-1}(0)\subset\C^n$. Then the following hold.

\begin{enumerate}
\item[(a)] $\h A\simeq\C[\h X]$.
\item[(b)] If $\p\in {\rm LND}\,(A)$ then
$\h\p \in {\rm LND}\,(\h A)$, and for any $a\in A^{\p}$,
 $\h a\in\h A^{\h\p}$. 
\item[(c)] \footnote{A similar fact remains true 
for any affine domain \cite{ML3}.}
For any non-zero Jacobian derivation
$$\p={\p(p,\,p_2,\dots, p_{n-1},\,*)\over\p(x_1,\dots,x_n)}$$
of the algebra $A=\C[X]$ there exists an equivalent one
$$\p'={\p(p,\,p'_2,\dots,p'_{n-1},\,*)\over\p(x_1,\dots,x_n)}$$
such that the $d$-principal parts $\h p,\,\h p'_2,\dots,\h p'_{n-1}$
are algebraically independent.
\item[(d)] 
If the latter condition holds then the associated derivation
$\h{\p'}\in$Der$(\h A)$ of the associated graded algebra 
$\h A$ is
also a Jacobian one, namely,
$$\h{\p'}={\p(\h p,\,\h {p'}_2,\dots,\h {p'}_{n-1},\,*)
\over\p(x_1,\dots,x_n)}\,\,.$$
\end{enumerate}
\ethm

\query{(label=GDI)}\no\bsit\la{GDI}  {\it Graded invariants.}
{\rm For a graded algebra
$\h A$, we denote by ${\rm Dk_{gr}}\,(\h A)$ the following
 `graded' version of the Derksen invariant:
$${\rm Dk_{gr}}\,(\h A)= \C\Big\lbrack\bigcup_{\p\in {\rm LND_{gr}}(\h
  A)} \h A^{\p}\Big\rbrack\,.$$
The way we use in the next section the Derksen invariant
(similar to that of \cite{De, KaML1, ML1, Za2})
is based on the next simple lemma.
}\esit

\query{(label=DML)}\no\no\blem\la{DML}
Given an irreducible hypersurface $X\subseteq\C^n$,
suppose that the algebra $A=\C[X]$ is equipped with 
a weight degree function $d_A$
such that
the hypersurface $\h X\subset \C^n$ is also irreducible. If
$${\rm Dk_{gr}}\,(\h A)\subset {\h A}_{\le 0}:= 
\bigoplus_{t\le 0} A_t/A'_t\subset \h A$$ 
then $${\rm Dk}\,(A)\subset A_0=\{a\in A\,|\,d_A(a)\le 0\}\,.$$
Henceforth,
$X\not\simeq \C^{n-1}$ unless $A=A_0$. \elem

\no\proof By \ref{KML2} (a),(b)
for every $\p\in {\rm LND}\,(A)$ and for every
$\p$-constant $a\in A^{\p}$, we have ${\h a} \in {\rm ker}\,{\h \p}$,
where ${\h \p}\in {\rm LND_{gr}}({\h A})$, whence
${\h a} \in {\h A}_{\le 0}$. Therefore, $d_A(a)\le 0$ i.e.,
$a\in A_0$, as stated. \qed

\query{(label=ssec-var)}\subsection{Variables 
in polynomial rings}

\query{(label=var)}\no\bsit\la{var} {\rm 
For a commutative ring $B$, we let 
$B^{[n]}=B[y_1,\ldots,y_n]$.
A polynomial $p\in B^{[n]}$ is called a 
{\it $B$-variable} (or simply
a {\it variable} if no ambiguity occurs) 
if it is a coordinate ($p=p_1$)
of an automorphism $\beta=(p_1,p_2,\ldots,p_n)
\in\Aut_B B^{[n]}$ 
(thus $\beta|_B= {\rm id}_B$). 
The set of all $B$-variables
of $B^{[n]}$ is denoted as $\Var_B  B^{[n]}$.
For $B=\C[x_1,\cdots,x_n]$, a $B$-automorphism 
is called in brief an {\it
$(x_1,\cdots,x_n)$-automorphism}, and a $B$-variable is
also called an  $(x_1,\cdots,x_n)$-variable. 
}\esit

\query{(label=slice)} \no\brem\la{slice}{\rm
It is easily seen that a polynomial $p\in \C^{[n]}$ is an $(x_1,\dots,x_k)$-variable 
(where $0\le k\le n-1$)
if and only if $X:=p^{-1}(0)\simeq \C^{n-1}$ and
there exists $\p\in \LND\,( \C^{[n]})$ with 
$x_1,\dots,x_k\in {\rm ker}\,\p$ such 
that $\p p = 1$ (i.e.,
$p$ is a {\it slice of} $\p$).
Indeed, $\p p = 1$ implies that
$X$ can be identified with the orbit space 
$\C^n/\varphi_{\p}={\rm spec}\,(\C^{[n]})^{\varphi_{\p}}=
{\rm spec}\,({\rm ker}\,\p)$ of the associated $\C_+$-action
$\varphi_{\p}$ on $\C^n$ (see \ref{ACT}), and so
$\C^{[n]}=({\rm ker}\,\p)[p]\cong \C^{[n-1]}[p]$.

Notice that as $\C^n\simeq X \times\C$, 
the assumption
$X\simeq \C^{n-1}$ above is
superfluous provided that 
the Zariski Cancellation Conjecture
holds.
}\erem

\query{(label=spe)}\no\bsit\la{spe} {\rm
For a polynomial $q=q(x,y,z,\ldots)\in\C^{[n]}$, 
we denote by 
$q_{\lambda}\,\,\,(\lambda\in\C)$ the specialization
$q(\lambda,y,z,\ldots)\in\C^{[n-1]}$. 

We say that a polynomial
$p(x,y_1,\cdots,y_n)\in \C[x][y_1,\cdots,y_n]$ is a 
{\it residual $x$-variable} if for every 
$\lambda\in \C$,
the specialization $p_{\lambda}$ is a variable of 
$\C[y_1,\cdots,y_n]$. It is easily seen
that any $x$-variable is a residual $x$-variable.
}\esit

\query{(label=nil)}\no\bsit\la{nil} 
{\rm Let $f\in B$. An element 
$b\in B$ is said to be {\it invertible} 
(resp., {\it nilpotent})
$\mod f$ if its image $[b]\in B/(f)$ is so. 
Thus $b$ is nilpotent $\mod f$ if and only if 
$b\in {\rm rad} \,(f)$. 
Observe that if $B$ is UFD and 
$f=\prod_{i=1}^n f_i^{a_i}$ 
is a canonical factorization then 
$b\in {\rm rad} \,(f)\Leftrightarrow 
b\in (f^{\rm red})$, where $f^{\rm red}:=
\prod_{i=1}^n f_i\in B$. 
}\esit

To detect variables in polynomial rings, 
the following results will be useful.
The statement (a) below is due to P. Russell \cite[Prop. 2.2]{Ru}, 
and is reproved in \cite{Ve} basing on (b) 
(see \cite[Prop. 1.4]{Ve}).

\query{(label=RuVe)}\no\bprop\la{RuVe} For a
commutative ring $B$, the following hold. 

\begin{enumerate}
\item[(a)] For any $f$, $b_0$, $b_1\in B$ and $q\in B[z]$
with $b_1$ invertible
$\mod f$ and $q$ nilpotent $\mod f$, we have
$p:=fu+b_0+b_1z+q(z)\in \Var_B B[z,u]$.
\item[(b)] For $f\in B$, we let $B_f:=B/(f)$, 
and we denote by $\rho : B \to B_f$ 
the canonical surjection.
If $p\in \Var_B B[z,u]$ is such that 
$\rho(p(z,0))\in \Var_{B_f} B_f[z]$ 
then also
$p_f(z,u):=p(z,fu)\in \Var_B B[z,u]$.
\item[(c)] Let $B=\C[x]$ and $f\in B$. 
If $p\in \Var_B B[y,z,u]$ is such that 
for every root $x_0$ of $f$, 
the specialization $p_{x_0}(y,z,0)$ is a variable of 
$\C[y,z]$ \footnote{Or equivalently,  
$\rho(p(y,z,0))\in \Var_{B_f} B_f[y,z]$.}
then also
$p_f(y,z,u):=p(y,z,fu)\in \Var_B B[y,z,u]$.
\end{enumerate}
\eprop

\no {\it Proof of (c).} It suffices to prove (c)
for $f=x$ and then conclude by induction on the degree of $f$. Let
$\gamma\in \Aut_B B[y,z,u]$ be such that
$\gamma(y)=p(x,y,z,u)$. Denote $\gamma_0$ the specialization of
$\gamma$ at $x=0$. Clearly, $\gamma_0(y)=p(0,y,z,u)$ is a variable of
$\C[y,z,u]$; in particular, 
$$
\C[y,z,u]/(p(0,y,z,u))\simeq \C^{[2]}\,.
$$
If $\phi$ denotes this isomorphism then 
$$
\C[y,z,u]/(p(0,y,z,u),u)\simeq \C^{[2]}/(\phi(u))\,.
$$
By our assumption, 
$p(0,y,z,0)$ is
a variable of $\C[y,z]$, and so \footnote{Hereafter
under $(f,g)$ we mean the ideal 
generated by $f$ and $g$.}
$$
\C[y,z,u]/(p(0,y,z,u),u)\simeq \C[y,z]/(p(0,y,z,0))\simeq \C^{[1]}\,.
$$
Hence by the Abhyankar-Moh-Suzuki Theorem, $\phi(u)$ is a variable
in $\C^{[2]}$. Therefore we may assume that $\phi(u)=u$ i.e., 
there is a
$\C[u]$-isomorphism 
$$
 \C[u][y,z]/(p(0,y,z,u))\simeq \C[u]^{[1]}\,.
$$
By the Abhyankar-Moh-Suzuki theorem as generalized by Russell and
Sathaye \cite[Thm. 2.6.2]{RuSat}, $p(0,y,z,u)$ is a
$\C[u]$-variable of $\C[u][y,z]$. Let $\alpha_0$ 
be the corresponding
$\C[u]$-automorphism of $\C[u][y,z]$
such that $\alpha_0(y)=p(0,y,z,u)$, and denote by 
$\bar \gamma$
the composition 
$$  
\bar \gamma:=\gamma\gamma_0^{-1}\alpha_0\in \Aut_B B[y,z,u]\,.
$$
Then we get
$$
\bar \gamma(y) = \gamma\gamma_0^{-1}\alpha_0(y)=
\gamma\gamma_0^{-1}(p(0,y,z,u))= \gamma(y)= p(x,y,z,u)
$$
and 
$$
\bar \gamma(u)=\gamma\gamma_0^{-1}\alpha_0(u)=
\gamma\gamma_0^{-1}(u)
\equiv u \mod x\,.
$$
If $\sigma_u$ is the $(x,y,z)$-automorphism of 
$\C(x)[y,z][u]$ defined
by $\sigma_u(u)=xu$ then the composition
$\sigma_u\bar \gamma \sigma_u^{-1}$ (which {\it \`a priori} is an $x$-automorphism of $\C(x)[y,z,u]$) 
actually is a $B$-automorphism of
$B[y,z,u]$. Thus 
$\sigma_u\bar \gamma \sigma_u^{-1}(y)=p(x,y,z,xu)$ is
an $x$-variable of $\C[x][y,z,u]$.\qed \medskip

\begin{cor}\la{corRuVe}\query{(label=corRuVe)}
Let $p\in \C[x][y,z,u]$ be an $x$-variable. 

\begin{enumerate}
\item[(a)] If for
$q(x)\in \C[x]$, $p_q:=p(x,y,z,q(x)u)$ 
is a residual $x$-variable 
then it actually is an $x$-variable. 
\item[(b)] Consequently, $p(x,y,z,q(x)u)$ is an $x$-variable 
if and only if $p(x,y,z, q_{\rm red}(x)u)$ is so.
\end{enumerate}
\end{cor}

\proof (a) As $p_q$ 
is a residual $x$-variable,
for every root $x_0$ of $q$, 
$p_{x_0, 0}:=p(x_0,y,z,0)$ is a variable 
of $\C[y,z,u]$. 
In particular,
$\C[y,z,u]/(p_{x_0, 0})=\C[y,z]/(p_{x_0, 0})
\otimes\C[u]\cong \C^{[2]}$. 
It follows that $\C[y,z]/(p_{x_0, 0})\cong \C^{[1]}$, 
and then by
the Abhyankar-Moh-Suzuki Theorem   
$p_{x_0, 0}$ is a variable of $\C[y,z]$. 
It follows from 
\ref{RuVe}(c) that $p(x,y,z,q(x)u)$ 
is an $x$-variable, 
as stated. 

The proof of (b) relies on (a) and 
on the fact that 
$p(x,y,z,q(x)u)$ and 
$p(x,y,z,q_{\rm red}(x)u)$ are simultaneously 
residual $x$-variables. 
\qed\medskip

The proof of the following results is inspired 
by those of 
\cite[Thm. 1.1]{SY1} and \cite[Thm. 1.4]{SY2}.

\query{(label=isom)}\bprop\la{isom}
For any commutative ring $B$, the following hold. 

\begin{enumerate}
\item[(a)] For arbitrary pair of polynomials 
$p, q\in B[y]$ 
there exists an automorphism 
$\gamma=\gamma_{p,q}\in\Aut_B B[y,v]$  
such that 
$$\gamma\bigl ((y-q(p(y)),v)\bigl )
=(y-p(q(y)),v)\,.$$
\item[(b)] In particular, for $q(y)=-ay$ 
with $a\in B$  
we have: 
$$\gamma\bigl ((y+ap(y),v)\bigr )=(y+p(ay),v)\,.$$
Moreover, if $a^k|p(0)$ then there
exists an automorphism $\gamma_k\in\Aut_B B[y,v]$ 
such that
$$\gamma_k\bigl ((y+ap(y),v)\bigl )=(y+
\frac{p(a^{k+1}y)}{a^k},v)\,.$$ 
\end{enumerate}
\eprop

\no \proof
It is not difficult to verify 
that the desired automorphism 
$\gamma$ and its
inverse $\gamma^{-1}$ can be defined as follows :
$$
\left\{
           \begin{eg}
             \gamma(y)=v+q(y)\\
             \gamma(v)=y-p(v+q(y))    
           \end{eg}
         \right. 
\left\{
           \begin{eg}
             \gamma^{-1}(y)=v+p(y)\\
             \gamma^{-1}(v)=y-q(v+p(y)) \,\,.   
           \end{eg}
         \right. 
$$
Iterating $\gamma$ yields 
an appropriate $\gamma_k$ as needed in (b).
\qed\medskip

\query{(label=isomm)}\bcor\la{isomm} 
\begin{enumerate}
\item[(a)] 
For any 
pair $p,q\in B[y]$ the homomorphism 
$$ 
\begin{array}{rccl}
   \phi \;: &  B[y]/(y-q(p(y))) & 
\longrightarrow & B[y]/(y-p(q(y))) \\
        &        y        & \longmapsto     & q(y)
\end{array} 
$$
is an isomorphism. 
\item[(b)] In particular, for any $a\in B$
and $p\in B[y]$ the homomorphism
$$ 
\begin{array}{rccl}
   \phi \;: &  B[y]/(y+ap(y)) & 
\longrightarrow & B[y]/(y+p(ay)) \\
        &        y        & \longmapsto     & ay
\end{array} 
$$
is an isomorphism. Moreover, if $a^k|p(0)$ then also
$$
\begin{array}{rccl}
   \phi_k \;: &  B[y]/(y+ap(y)) & \longrightarrow & B[y]/(y+\frac{p(a^{k+1}y)}{a^{k}}) \\
        &        y        & \longmapsto     & a^{k+1}y
\end{array} 
$$
is an isomorphism.
\end{enumerate}
\end{cor}

\no \proof
We just apply \ref{isom} 
and the obvious isomorphism $B[y]/(f(y))\simeq
B[y,v]/(f(y),v)$.
\qed \medskip

\query{(label=jac)}\bsit\la{jac}{\rm
We denote by S$\Aut\,\C^{[n]}\subset\Aut\,\C^{[n]}$ 
the subgroup 
which consists of the
automorphisms $\alpha$ with 
the Jacobian jac$(\alpha)=1$. 
}\esit

For the polynomial extension $B^{[2]}$ 
of the algebra $B:=\C[x]$, 
we have the following fact. 

\query{(label=sp2)}\no\bprop\la{sp2}
\begin{enumerate}
\item[(a)]
  For any set of distinct points 
$\lambda_1,\ldots,\lambda_n\in\C$,
the homomorphism of multi-specialization
$$
  \begin{array}{rccl}
    \Phi: & \mbox{\rm SAut}_{B} B[y,z] & \longrightarrow &
    ({\rm SAut}\,\C[y,z])^n \\
          & \ga & \longmapsto & ({\ga}_{\lambda_1} \tr {\ga}_{\lambda_n})
  \end{array}
$$
is surjective.
\item[(b)] The multi-specialization map 
  $$
    \begin{array}{rccl}
      \Phi': & \mbox{\rm Var}_{B}\,B[y,z] & \longrightarrow &
      (\mbox{\rm Var}\,\C[y,z])^n \\
            & p & \longmapsto & ({p}_{\lambda_1} \tr {p}_{\lambda_n})
    \end{array}
  $$
  is surjective.
\end{enumerate}
\eprop

\no\proof
By the Jung-van der Kulk decomposition theorem 
it is enough to show that the image of $\Phi$ contains 
the elements $g_{\alpha}\in$(SAut$\,\C[y,z])^n$, where 
  $$
   g_{\alpha}:= (\mbox{Id} \tr \mbox{Id},\alpha,\mbox{Id} \tr \mbox{Id})\,,
  $$
  and $\alpha\in$ SAut$\,\C[y,z]$ is a triangular automorphism of the form
  $$
    \begin{array}{cccc}
      \mbox{either} &
         \left\{
           \begin{eg}
             \alpha(y)=y+r(z)\\
             \alpha(z)=z 
           \end{eg}
         \right.                     & \mbox{ or } & \left\{
                                                            \begin{eg}
                                                              \alpha(y)=y\\
                                                              \alpha(z)=z+s(y).
                                                             \end{eg}
                                                          \right.
    \end{array} 
  $$
Indeed, any automorphism $\beta_c : (y,z)\longmapsto (cy+r(z), c^{-1}z)$
with $c\in\C^*$ can be decomposed as $\beta_c = \delta_c\circ\alpha$ with
$\alpha : (y,z)\longmapsto (y+c^{-1}r(z), z)$ as above and 
$\delta_c : (y,z)\longmapsto (cy, c^{-1}z)$. In turn, we have
$\delta_c=\gamma_1\circ \gamma_{-c}$, where
$\gamma_{c} : (y,z)\longmapsto (cz, -c^{-1}y)$ can be written as
$\gamma_{c}=\alpha'_{c}\circ\alpha''_{c}\circ\alpha'_{c}$
with
$$\alpha'_{c}: (y,z)\longmapsto (y, z-c^{-1}y)\qquad {\rm and}\qquad 
\alpha''_{c}: (y,z)\longmapsto (y+cz, z)\,.$$

  Without loss of generality we may assume that
$g_{\alpha}=(\alpha,\mbox{Id} \tr \mbox{Id})$. 
Let $\varphi\in B=\C[x]$ be such
  that $\varphi(\lambda_1)=1$, 
$\varphi(\lambda_2)=\cdots=\varphi(\lambda_n)=0$, 
and consider the triangular  automorphism 
   $\ga\in\Aut_{B}\,B[y,z]$, where
  $$  
    \begin{array}{ccc}
      \left\{
        \begin{eg}
          \ga(y)=y+\varphi(x)r(z)\\
          \ga(z)=z
        \end{eg}
      \right.                & \mbox{resp., } & \left\{
                                                       \begin{eg}
                                                         \ga(y)=y\\
                                                         \ga(z)=z+\varphi(x)s(z).
                                                       \end{eg}
                                                     \right. 
    \end{array} 
  $$   
  It is easily seen that, indeed, 
  $\Phi(\ga)=g_{\alpha}$. This proves (a).

(b)
  By (a) it suffices to show that every variable 
$p=\alpha(y)\in\Var_{\C}\,\C[y,z]$ with $\alpha\in\Aut_{\C} \C[y,z]$
can be regarded
  as a coordinate function of an automorphism 
$\alpha'\in $SAut$\,\C[y,z]$. Indeed, let $\alpha':=\alpha\circ \si$
with $ \si : (y,z) \longmapsto (y, j^{-1}z)$, 
where $j:=$ jac$(\alpha)$. Then 
$\alpha'(y)=\alpha(y)=p$ and 
$\alpha' \in $S$\Aut\,\C[x,y]$. Thus (b) follows.   
\qed\ms

\query{(label=ft)}\no\brem\la{ft} We would like to use this 
opportunity to indicate a flow in the proof of Theorem 7.2 
in \cite{KaZa1}
(which generalizes the Sathaye Theorem, as mentioned in the Introduction).
Namely, proving Claim 3 of Proposition 7.1 in \cite{KaZa1}
and carrying induction 
by the multiplicity of a root $x=0$ of the polynomial $p$, 
it has been forgotten to extend it over all the roots. Instead, 
one might apply 
\ref{sp2} (or \ref{RuVe}(b); cf. \cite[Thm. 3.6]{Ve}) 
which would simplify the proof considerably. \erem


\query{(label=S2)}\section{Simple 
modifications of acyclic 3-folds
along cylindrical divisors}\la{S2}

We focus below on a special 
case of affine modifications called
simple birational extensions (see \ref{sibi} below), 
applied to acyclic affine 3-folds. 
Our aim in this sections is to give a criterion 
for as when the acyclicity 
is preserved under such a modification
(cf. \cite[Thm. 3.1]{KaZa1}). In the special 
case when the divisor of modification is a cylinder,
we obtain in Theorem \ref{3f} below necessary conditions
for preserving the acyclicity. 
In Theorem \ref{3g} (cf. also \ref{c3g})
we show that these conditions are also 
sufficient, provided that the 
given acyclic 3-fold is a cylinder as well. 

\query{(label=simo)} \subsection{Simple affine modifications}
\la{simo}
\query{(label=sibi)} \bdefi\la{sibi} {\rm
A {\it simple birational extension}
of a domain $A$ (over $\C$) is an algebra $A':=A[g/f]$, 
where $f, g \in A$ are such that the ideal $I=(f,g) \subset A$
is of height 2 (i.e., the center of modification
$C=V(I)=D_f^{\rm red}\cap D_g^{\rm red}$ 
is of codimension 2 in $Y:={\rm spec}\,A$).
We also call $X:={\rm spec}\,A'$ a 
{\it simple affine modification} of the 
affine variety $Y$. 

In the sequel $A$ is UFD, and 
the above condition simply means that
$f, g\in A$ are coprime. More generally 
\cite[Prop. 1.1]{KaZa1},
any affine modification can be obtained as 
$A'=A[g_1/f_1,\ldots,g_n/f_n]$
with $f_i,g_i\in A$; here again, 
if $A$ is UFD then we may suppose 
$f_i,g_i$ being coprime $(i=1,\ldots,n)$.
}\edefi

\query{(label=OBS)}\bsit\la{OBS} {\rm Observe that the variety
 $X={\rm spec}\,A'$   can be realized
 as  the hypersurface in $Y\times \C$  with
 the  equation $fu+g=0$
 (where $u$ is a coordinate in $\C$), 
and the blowup morphism is just the first projection: 
$\si={\rm pr}_1|_X:X\to Y$.
 The exceptional divisor
 $E=\s^{-1}(C)=\s^{-1}(D)=\{f=0\}\subseteq X$ 
is cylindrical: $E\simeq C\times\C$
 (cf. \cite[Prop. 1.1]{KaZa1}).
}\esit

We have the following simple but useful lemma.

\query{(label=smooth)} \no\blem \la{smooth}
 In the notation as above, 
assume that the affine variety $Y:={\rm spec}\,A$ is smooth. 
Then the
simple affine modification $X={\rm spec}\,A[g/f]$ is also
smooth if and only if 

\begin{enumerate}
\item[(i)] the divisor $D_g=g^*(0)$ is smooth
and reduced at each point of the  center 
$C=D_f^{\rm red}\cap D_g^{\rm red}$, and 
\item[(ii)] $D_f$ 
and $D_g$ meet transversally at those points 
of $C$ where the divisor $D_f$ is also smooth and reduced. 
\end{enumerate}

\no Furthermore, if $X$ is smooth then $df=0$ 
at each singular point of 
the center $C$. 
\elem

\no\proof The second assertion 
easily follows from the first one.
To prove the first one, let $F=0$ be the equation 
of the hypersurface $X\subseteq Y\times\C$
with $F:=uf-g\in\C[Y][u]$ being irreducible.
As $X\backslash E\simeq Y\backslash D$ is smooth, 
we should only control the smoothness of $X$ at
the points of the exceptional divisor 
$$E=\{f=g=0\}=C\times\C\subseteq Y\times\C\,.$$ 
At a point $Q=(P,u)\in E$, we have:
$$dF=udf+dg+fdu=udf+dg=0$$ if and only if
$dg=-udf$ is proportional to $df$. 
Now the assertion easily follows. \qed\medskip

\query{(label=prac)}
\subsection{Preserving acyclicity: 
necessary conditions}\la{prac}
In this subsection we adopt the following
convention and notation. 

\no \bconv \la{conv0}\query{(label=conv0)} {\rm
\begin{enumerate}
\item[(i)] $X$ and $Y$
denote smooth, acyclic
  affine 3-folds \footnote{Hence the algebras 
$A$ and $A'$ are UFD's, see the proof of \ref{new}.} such that
 \item[(ii)] $ A'=\C[X]$ is
 a simple birational extension of the algebra
 $A=\C[Y]$ i.e., $A'=A[g/f]$  with  coprime elements $f,\,g\in A$, 
and
\item[(iii)]   $D:=D_f^{\rm red}\simeq\G\times\C$ (where $\G$ is
 an affine curve) is a cylindrical divisor.
\end{enumerate}}\econv

\query{(label=Ci)} \no\bsit \la{Ci}
{\rm  Denote by $\pi:D\to \G$
 the morphism
 induced by the canonical projection
 $\G\times\C\to\G$.  By abuse of notation,
 we equally denote by $\pi$ the restriction
 $\pi|_C:C\to \G$.
 Thus we have the following commutative diagram:

\begin{picture}(150,75)
\unitlength0.2em
\put(45,25){$E\simeq C\times\C\longrightarrow C$}
\put(45,5){$D\simeq \G\times \C\longrightarrow \G$}
\put(47,22){$\vector(0,-1){11}$}
\put(50,16){{$\sigma$}}
\put(82,22){$\vector(0,-1){11}$}
\put(85,16){{$\pi$}}
\end{picture}

\no If $ C=\bigcup_{i=1}^{n'} C_i$ resp.,   $\G=\bigcup_{j=1}^n \G_i$
 is the irreducible decomposition then
 the irreducible components of the divisor
 $E$ resp.,   $D$ are $E_i\simeq C_i\times \C$,
 resp.,   $D_j\simeq \G_j\times\C$. As in \ref{mm} we 
denote by $\M=(m_{ij})$ the
 multiplicity matrix of $\si$ . Recall that by \ref{mm}, \ref{gone}
and \ref{smooth},
$m_{ij}\geq 1$ if and only if $C_j\subset
    D_i^{\rm red}$, and  $m_{ij}= 1$ if and only if
    the surfaces $D_i^{\rm red},\,D_g^{\rm red}\subseteq Y$
    meet transversally at general points of the curve $C_j$.
}
\esit
 
The terminology in the following definitions comes from the real picture
corresponding to our situation.

\no \bsit \la{def}\query{(label=def)}
{\rm Notice that for each $j=1,\dots,n$ there exists
$i\in\{1,\dots,n\}$ such that $\pi(C_j)\subseteq \G_i$,  and this index
$i=i(j)$ is unique unless $\pi|_{C_j}={\rm const}$.

An irreducible component $C_j$ of the curve $C$ is called
{\it vertical}
if $\pi|C_j={\rm const}$  (i.e. deg$\,(\pi|_{C_j})=0$) and {\it
  non-vertical} otherwise (thus the vertical components
of $C$ are disjoint and each of them is isomorphic to $\C$). The
  uniqueness of the index $i=i(j)$ 
for a non-vertical component $C_j$ and
the unimodularity of $\M$ (see \ref{new}) 
imply that the $j$-th column of the matrix
$\M$ is the $i$-th vector of the standard basis
$(\bar{e}_1,\dots\bar{e}_n)$ in $\R^n$, and two different non-vertical
components $C_j$ and $C_{j'}$ of $C$ project into two different
irreducible components $\G_i$ resp., $\G_{i'}$ of $\G$.
Hence up to reordering, we may assume that
$C_1,\dots,C_k$ are the non-vertical components of $C$ and
$\pi(C_i)\subseteq\G_i$, $i=1,\dots, k$. Then we have
\begin{displaymath}
\mbox{$\M=$}
\left( \begin{array}{ccc}
I_k & \vert & B \\
-\!\!-\!\!- &\vert& -\!\!-\!\!- \\
0 &  \vert & B' \\
\end{array} \right)
\end{displaymath}
with a unimodular matrix $B'$. Consequently, for every $i=1,\ldots,k$,
the surfaces $D_i^{\rm red}$ and $D_g^{\rm red}$ meet transversally
at general points of the curve $C_i$.
}\esit

\query{(label=non-smth)} \no\bsit\la{non-smth}
{\rm The irreducible components $\G_1,\dots,\G_k$  are also called
 {\it non-vertical} resp.,
 $\G_{k+1},$$\dots,\G_n$  are called {\it vertical}. Among the
 non-vertical components
 $C_i$ resp., $\G_i$ we distinguish those with
 deg$\,(\pi|_{C_j})=1$ which we call {\it horizontal} and those with
 deg$\,(\pi|_{C_j})\geq 2$ which we call  {\it slanted}. We
 reorder again to obtain that $C_1,\dots,C_h$ resp., $\G_1,\dots,\G_h$
 are the horizontal components of 
$C$ resp., $\G$, and
 $C_{h+1},\dots,C_k$ resp., 
$\G_{h+1},\dots,\G_k$ are the slanted ones.
An irreducible component of $C$ resp., 
$\G$ which is not horizontal
is refered to as a {\it non-horizontal} 
component; {\it $\CH$}
denotes  the union
of all horizontal components of $C$. 
In the same way we define
$\CNH,\,\,\CV,\,\,\CNV,\,\,\CS$, 
$\GH,\,\,\GNH,\,\,D_{\rm horiz}$, etc. 
Thus we have:
$$\GH=\bigcup_{i=1}^h \G_i,\quad 
\GS=\bigcup_{i=h+1}^k \G_i,\quad
\GV=\bigcup_{i=k+1}^n \G_i\,,$$
and similarly
$$\CH=\bigcup_{i=1}^h C_i,\quad 
\CS=\bigcup_{i=h+1}^k C_i,\quad
\CV=\bigcup_{i=k+1}^n C_i\,.$$
}\esit

\query{(label=factor)} \no\bsit\la{factor}
{\rm 
Let furthermore $f=\prod_{j=1}^n f_j^{a_j}
=\fh\cdot \fnh$
be  factorizations such that
$$f_j^*(0)=D_j^{\rm red}\,\,\,(j=1,\ldots,n),
\quad (f_{\rm horiz})^{-1}(0)=D_{\rm horiz}\simeq
\G_{\rm horiz}\times\C$$
and
$$(f_{\rm non-horiz})^{-1}(0)=D_{\rm non-horiz}
\simeq \G_{\rm non-horiz}\times\C\,.$$
}\esit

\query{(label=isol)} \no\bsit\la{isol}
{\rm An irreducible component is called 
{\it isolated} if it is a connected
component.
}\esit

Let us give a typical example which illustrates our definitions. 

\query{(label=exmp)} \no\bexa\la{exmp} {\rm Letting $Y=\C^3$ 
with coordinates $x,y,z$, set
$f=xy, g=y+xz$.  Then $X=\{xyu+y+xz=0\}\subset\C^4=Y\times\C$, $D=\{xy=0\}\subset\C^3,\,\,\,\si(x,y,z,u)=(x,y,z)$ and
$E=\{xy=y+xz=0\}\subseteq \C^4$. Therefore, 
$\G = \{ xy=0 \} \subset \C^2$,
$\GH = \{ y=0 \}$ and $ \GV=\{x=0\}$ whereas $C=\CV\cup\CH$
with $\CH=\{y=z=0\}$ and $\GV=\{x=y=0\}$. 
}\eexa

 The main result of this subsection 
 is the following theorem.

\query{(label=3f)} 
\no\bthm\la{3f} Let $X$ and $Y$ be
 smooth acyclic  affine 3-folds satisfying
 the conditions (ii) and (iii) of \ref{conv0}. 
 Then  (in the notation of \ref{def}-\ref{isol}) 
 the following hold.

\begin{enumerate}
\item[($\alpha$)]
$\pi|(\CH\setminus \CNH):\CH\setminus \CNH\to\GH\setminus \GNH$
is an isomorphism.
\item[($\beta$)]
  The slanted  components $C_{h+1},\dots,C_k$ and  
  $\G_{h+1},\dots,\G_k$
  are isolated and homeomorphic to $\C$.
\item[($\gamma$)]
  $\fnh=p\circ f_{h+1}$ with $p\in\C[z]$,
  and every non-horizontal component of $\G$ is homeomorphic to
  the affine line $\C$. In other words,
  $$
   \fnh=c\cdot \prod_{i=h+1}^n(f_{h+1}-\lambda_i)^{a_i}
  \qquad (c\in\C^*,\,\,\,\lambda_{h+1}=0)$$
   with $\pi(f_{h+1}^{-1}(\lambda_i))=\G_i$
  homeomorphic to $\C\,\,\,(i=h+1,\ldots,n)$.
\item[($\delta$)]
  Up to a further reordering, 
 the multiplicity matrix has the following form:
  \begin{displaymath}
  \mbox{$\M=$}
  \left( \begin{array}{ccccc}
  I_h & \vert & 0 & \vert & B_0 \\
  -\!\!-\!\!- &\vert& -\!\!-\!\!-&\vert& -\!\!-\!\!- \\
  0 & \vert & I_{k-h} & \vert & 0 \\
  -\!\!-\!\!- &\vert& -\!\!-\!\!-&\vert& -\!\!-\!\!- \\
  0 & \vert & 0 & \vert & I_{n-k} \\
  \end{array} \right)
  \end{displaymath}
\end{enumerate}

\no Consequently (by \ref{new}), 
$\si_*:\pi_1(X)\to\pi_1(Y)$ is an isomorphism. Moreover,
$X$ is contractible (and hence diffeomorphic to $\R^6$)
if and only if $Y$ is so. 
\ethm

 The  rest of this subsection is  devoted to the proof of
 Theorem \ref{3f}.
 It is convenient to introduce the following terminology and  notation.

\query{(label=nm)}\noindent \bsit\la{nm}
{\rm Let $F$ be a curve. We say that a
  point  $P\in F$
is {\it multibranch} (resp.,   {\it unibranch})
if it is a center of $\mu_P=\mu_P(F)>1$ (resp.,   $\mu_P=1$)
local analytic branches of $F$. We denote by $\nF$
the normalization of $F$ and by ${\ol F}$
its smooth complete model.
The points of
${\ol F}\sm \nF$ are called the {\it punctures} of $F$.
A morphism  of curves $\rho:F\to G$ can be  lifted
to the normalizations resp.,   the completions; we denote the lift
by ${\rho^{\rm norm}}:\nF\to\nG$
resp.,   ${\ol \rho}:{\ol F}\to{\ol G}$.  }\esit

\query{(label=bs)}\noindent \bsit\la{bs}
{\rm Consider further an irreducible smooth curve  $F$
of genus $g$ with $n$ punctures,
and let the
1-cycles $a_1,\dots,a_g,\,b_1,\dots,b_g$
on $F$  provide a
symplectic basis of the group
$H_1({\ol F})=H_1({\ol F};\,\Z)$.
Then there is
an injection $H_1({\ol F})\hra H_1(F)$
onto the subgroup
generated by the classes
$[a_1],\dots,[a_g],\,[b_1],\dots,[b_g]$;
we  may identify the group $H_1({\ol F})$
with its image in $H_1(F)$.
The group $H_1(F)$ being freely generated
by the classes
$[a_1],\dots,[a_g],\,[b_1],\dots,[b_g]$, $[c_1],\dots,[c_{n-1}]$,
where $c_1,\dots,c_n$ are simple 1-cycles
around the punctures of $F$,
we have a (non-canonical) decomposition
\query{(label=di)}\be\la{di} H_1(F)\simeq H_1({\ol F})\oplus G(F)\quad{\rm with}\quad
G(F):=\big\langle[c_1],\dots,[c_{n}]\big\rangle\simeq\Z^{n-1}\,\ee
($\sum_{i=1}^n [c_i]=0$ being the
only relation between the generators
$[c_1],\dots,[c_n]$ of the group $G(F)$).}\esit

\query{(label=pg)}\noindent \bsit\la{pg}
{\rm We denote by $S_{\G}=S_{\G}^{(0)}\cup S_{\G}^{(1)}$ (where $S_{\G}^{(0)}\cap S_{\G}^{(1)}=\emptyset$) 
the finite subset
of the curve $\G$ such that a point
$P$ belongs to $S_{\G}$ if and only if it  satisfies at least
one of the following three
conditions:

\begin{enumerate}
\item[(i)] $P=\pi(C_i)$ for a
vertical component $C_i$ of $C$ ($\Leftrightarrow P\in S_{\G}^{(1)}$);
\item[(ii)] $P=\pi(Q)$ for a
multibranch point $Q$ of $C$;
\item[(iii)] $P$ is a
multibranch point of $\G$.
\end{enumerate}}\esit

\query{(label=ps)} \noindent \bsit\la{ps}
{\rm Set
 $S_C =\pi^{-1}(S_{\G})\subseteq C$ and $S_C^{(i)} 
=\pi^{-1}(S_{\G}^{(i)})$,
$i=0,1$. Thus the analytic set $S_C$ contains
the union $S_C^{(1)} =\CV$ of vertical components
of $C$, whereas the 
residue set $S_C^{(0)}=S_C\sm S_C^{(1)}$ is finite. 

It is easily seen that if
$Q\in S_C^{(0)}$ then over 
any local analytic branch $B_i$ of $\G$ at
the point $P:=\pi(Q)$ there is 
at least one local analytic branch 
$A_j\subseteq\pi^{-1}(B_i)\cap D_g^{\rm red}$ 
of $C$ at $Q$. Indeed, the surfaces 
$\pi^{-1}(B_i)\simeq B_i\times\C$
and $D_g$ meet at $Q$, and so the polynomial
$g|_{\pi^{-1}(B_i)}\in {\mathcal O}[B_i][z]$ vanishes 
at $Q=:(P,z_0)$,
but its specialization at $P$ is nonzero, as $Q\in S_C^{(0)}$
means that no vertical component of $C$ passes through $Q$. 
Consequently $\mu_Q(C)\ge \mu_P(\G)$ and $P\not\in\GV$, whence
$\pi^{-1}(\GV)\subseteq S^{(1)}_C=\CV$.
}\esit

\query{(label=str)}\noindent \bsit\la{str} {\rm
We let $\G^*=\G\sm S_{\G}$ 
(more generally, $\G^*_{\rm
  something}=\G_{\rm something} \cap \G^*$) and
$C^*=C\sm S_C=\pi^{-1}(\G^*)=\pi^{-1}(\GNVs)\subseteq C$.}
\esit

\query{(label=li)} \noindent \bsit\la{li}
{\rm For a complex hermitian manifold $M$
and a closed analytic subset $T$ of $M$,
by a {\it link} $lk_P(T)$ of $T$ at a
point $P\in T$
we mean the
intersection $T\cap S_{\varepsilon}$ of
$T$ with a small enough sphere $S_{\varepsilon}$
in $M$ centered at $P$. We also call {\it link}
the corresponding homology class
$[lk_P(T)]\in H_*(T\sm\{P\})=H_*(T\sm\{P\};\,\Z)$,
and we still denote it simply by $lk_P(T)$.  }\esit

\query{(label=hc)} \noindent \bsit\la{hc}
{\rm We denote by $\HG$ resp.,
$\HC$ the subgroup of
the group $H_1(\G^*)=H_1(\G^*;\,\Z)$
resp.,   $H_1(C^*)$ generated by the links\footnote{Here 
we are in the
  special case when  $M=X$ resp., $Y$ and
$T=C$ resp. $\G$.} $lk_P(\G)$ resp., $lk_Q(C)$ of
the points $P\in\SG$
resp., $Q\in S_C^{(0)}$.  Notice that
$ lk_P(\G)=\sum_{i=1}^{\mu} lk_P(B_i)$, 
where $B_1,\dots,B_{\mu}$ are the
local branches of $\G$ at $P$ and $\mu:=\mu_P (\G)$.  Clearly,
$H_{\G}\subseteq G(\G^*)$ resp., $H_C\subseteq G(C^*)$, 
and so we obtain (non-canonical) isomorphisms
\query{(label=nlo)}\be\la{nlo}
H_1(C^*)/H_C\simeq H_1(\overline {C^*})\oplus [G(C^*)/H_C]\ee 
resp.,
\query{(label=nlo*)}\be\la{nlo*} H_1(\G^*)/H_{\G}
\simeq H_1({\ol \G})\oplus [G(\G^*)/H_{\G}]\,.\ee} \esit

 The next proposition is our main technical
tool in the proof of Theorem \ref{3f}.

\query{(label=mainprop)} \no \bprop \la{mainprop}
Under the assumptions as in 
\ref{conv0}, consider the restriction $\pi=\pi_{|C^*} : C^*\to\G^*$.
Then $\pi_*(H_C)\subseteq H_{\G}$, and $\pi$ induces
the following isomorphisms:
\begin{eqnarray}\la{iso}
  {\widehat \pi_*}:H_1(C^*)/\HC\stackrel{\cong}{\longrightarrow}
  H_1(\G^*)/\HG\,,
\end{eqnarray}\query{(label=iso)}
\begin{eqnarray}\la{isocomp}
  {\ol \pi}_*:H_1(\overline {C^*})\stackrel{\cong}{\longrightarrow}
  H_1({\ol \G})\,,
\end{eqnarray}\query{(label=isocomp)}
\begin{eqnarray}\la{isoG}
  {\tilde\pi}_*:G(C^*)/H_C \stackrel{\cong}{\longrightarrow}
  G(\G^*)/H_{\G}\,.
\end{eqnarray}\query{(label=isoG)}
\eprop

The proof is based on Lemmas \ref{nle1} and \ref{nle2} below.
Let us introduce the following notation.

\no\bsit\la{notaq}\query{(label=notaq)}
{\rm
Denote  $$S_E=\si^{-1}(S_C)=\SC\times \C \subseteq E\quad {\rm resp., }\quad S_E^{(i)}=\si^{-1}(S_C^{(i)})\,\,\,(i=0,1)$$ and
$$ S_D=\pi^{-1}(\SG)\subseteq D\,,\qquad {\rm so\,\,\, that }\quad
S_D\simeq \SG\times \C \subseteq \G\times\C\simeq D\,.$$
Furthermore, set
$$E^*=\si^{-1}(C^*)=E\sm S_E = C^*\times\C\qquad {\rm resp., }
\qquad D^*= \pi^{-1}(\G^*)=D\sm S_D \simeq \G^*\times \C\,$$
and
$$X^*=X\sm S_E\qquad {\rm resp., }\qquad Y^*=Y\sm S_D.$$
 Thus
$$X\setminus E= X^*\sm E^*\qquad {\rm resp., }\qquad
Y\sm D=Y^*\sm D^*\,,$$
and so we have an isomorphism
$$\s|_{X^*\sm E^*}: X^*\sm E^*\to Y^*\sm
D^*\,.$$}\esit

\no\blem\la{nle1}\query{(label=nle1)} There are monomorphisms
\query{(label=ae)}\be\la{ae}
\rho_X : H_3(X^*)\to H_1(C^*)\qquad {\rm resp.,}\qquad \rho_Y :
H_3(Y^*)\to H_1(\G^*)\ee
such that $\pi_*\rho_X=\rho_Y\si_*$ and
$${\widehat \pi}_* : H_1(C^*)/\rho_X(H_3(X^*))\to
H_1(\G^*)/\rho_Y(H_3(Y^*))$$
is an isomorphism.\elem

\no\proof

The map of pairs $\s|_{X^*}: (X^*,X^*\sm E^*)\to (Y^*,Y^*\sm D^*)$
induces the following commutative diagram (where the horizontal
lines are exact homology sequences of pairs):

\begin{diagram}[height=2.3em,width=2em]
\dots& \rTo& H_i(X^*\sm E^*)               & \rTo& H_i(X^*)     &
\rTo^{i_*} & H_i(X^*,\,X^*\sm E^*)& \rTo{\partial_*} &H_{i-1}(X^*\sm
E^*)                            & \rTo& \dots\\
     &     &\dTo_{\wr\, \si_*}&     & \dTo_{\si_*} &
           & \dTo_{\si_*}        &                  &
\dTo_{\wr\, \si_*} &     &                                      &(*_1)    \\
\dots&\rTo & H_i(Y^*\sm D^*)               & \rTo& H_i(Y^*)     &
\rTo^{i_*}&H_i(Y^*,\,Y^*\sm
D^*)&\rTo^{\partial_*}&H_{i-1}(Y^*\sm D^*)&\rTo&\dots\\
\end{diagram}\\

\no {\bf Claim.} (a) $H_3( X^* \sm E^*)=H_3(Y^* \sm D^*)=0$
{\it and}
(b) $H_2(X^*)=H_2(Y^*)=0$.

\ms \no {\it Proof of the claim.} (a)
In virtue of the isomorphism
$$\s_*:H_3(X^*\sm E^*)=
H_3(X\sm E)\stackrel{\cong}{\longrightarrow}H_3(Y\sm D)=H_3(Y^*\sm D^*)$$
it is enough to show the vanishing of one of these groups, say,
of $H_3(X\sm E)$.
Let as before, ${\dot F}$ denote the one-point
compactification of a topological space $F$. As
$X$ is acyclic (whence by  \ref{Alex}, 
the Alexander duality can be applied to 
$\dot X$)
and moreover $E$ is closed in $X$, the Alexander duality
\query{\cite[p.301, 8.18]{Do}} gives
an isomorphism
\query{(label=ad)}\be\la{ad}
  H_3(X\sm E)\cong H^2({\dot E},\,\{\infty\})\,.
\ee
Consider
the homeomorphisms
$E\approx C\times \R^2,\,\,\,{\dot E}\approx ({\dot C}
\times{\dot \R}^2)/({\dot C} \vee {\dot \R}^2)$ and replace here
${\dot \R}^2$
by $S^2$. Since $({\dot C}\times S^2,\,{\dot C} \vee S^2)$ has the homotopy
type of a pair of cell complexes, by \cite[\query{p.102,} 4.4]{Do} we get
$$
H^2({\dot E},\,\{\infty\}) \cong {\widetilde H}^2 ({\dot E})
\cong H^2({\dot C}\times S^2,\,{\dot C} \vee S^2)\,.
$$
 The K\"unneth  formula for cohomology  \cite[\query{p.195,} (11.2)]{Mas}
yields a monomorphism
$$
 \mu: \sum_{p+q=2}H^p(\dot C,\{\infty\})\otimes
H^q(S^2,\{\infty\}) =:H\to
H^2({\dot C}\times S^2,\,{\dot C} \vee S^2)
$$
with the cokernel  
$${\rm coker}\,\mu=\sum_{p+q=3} {\rm Tor}\,(H^p(\dot
C,\{\infty\}), H^q(S^2,\{\infty\}))\,.$$ As $$H^0(\dot
C,\{\infty\})=H^0(S^2,\{\infty\})=H^1(S^2,\{\infty\})=0$$ we have
$ H=0$. The group $H^*(S^2,\{\infty\})$ being torsion free,
we also have coker$\,\mu=0$, and so
$H^2({\dot C}\times S^2,\,{\dot C} \vee S^2)=0$ as well.
Thus in view of (\ref{ad}),
$H_3(X\sm E)=0$.  This proves (a).

 (b) By the Alexander duality we obtain
$$H_2(X^*)= H_2(X \sm S_E)\cong
  H^3({\dot S}_E,\{\infty\})=H^3({\dot S}_E)\,.
$$
The topological space ${\dot S}_E$ is
homeomorphic to a bouquet of
  4-spheres $S^4$ and 2-spheres $S^2$.
The 4-spheres  are
provided by the one-point compactification of the product
$ S_C^{(1)}\times\C=\CV\times\C$ (recall (see \ref{def}) 
that  the components
of the curve $\CV$ are
disjoint and each one is isomorphic
to $\C$), whereas the 2-spheres $S^2$ are
provided by the one-point compactification
of the product
$S_C^{(0)}\times\C$.

Hence $ H^3({\dot S_E})\cong H_2(X^*)=0$.
Similarly, we have $H_2(Y^*)=0$.
This proves the claim. \qed\medskip

In virtue of the above claim, ($*_1$)
(with $i=3$) leads to the
commutative diagram\\
\begin{diagram}[height=2.3em,width=3em]
0 &\rTo & H_3(X^*) & \rTo^{i_*}& H_3(X^*,X^*\sm
E^*)&\rTo^{\partial_*}&
H_{2}(X^*\sm E^*)&\rTo & 0\\
 & &\dTo_{\si_*}& &\dTo_{\si_*}& &\dTo_{\wr\,\si_*}& & &             (*_2)\\
0&\rTo &H_3(Y^*)&\rTo^{i_*}& H_3(Y^*,Y^*\sm
D^*)&\rTo{\partial_*}&
H_{2}(Y^*\sm D^*)&\rTo& 0\\
\end{diagram}\\

\no Next we apply the Thom isomorphism  \cite[\query{p. 297,} 7.15]{Do}
to the pairs of  manifolds $(X^*,E^*)$ resp., $(Y^*,D^*)$
(cf. the proof of \ref{prhom}). 
Indeed, the curve $C^*$ resp., $\G^*$ being locally
unibranch (whence homeomorphic to its normalization),
$E^*=C^*\times\C$ resp., $D^*=\G^*\times\C$ is a topological manifold.
Notice that (in virtue of \ref{gone} and \ref{def}) 
$\si^*$ maps the transverse classes 
\cite[Ch. VIII\query{, p.315}]{Do} of $E^*$
into those of $D^*$.
By  functoriality of the cap-product \cite[\query{p.239,} VII.12.6]{Do},
for every $i\geq 2$
the following diagram is commutative:\\
\begin{diagram}[height=2.3em,width=3em]
H_i(X^*,\,X^*\sm E^*)&\rTo^{\cong}_{t_X}&\tH_{i-2}(E^*)\\
\dTo_{\si_*} & &\dTo_{\si_*}& & &                           (*_3)\\
H_i(Y^*,\,Y^*\sm D^*)&\rTo^{\cong}_{t_Y}& \tH_{i-2}(D^*)\\
\end{diagram}\\

\no  By making use of   ($*_3$) together with the following commutative
diagram (see  \ref{Ci}):\\
 \begin{diagram}[height=2.3em,width=3em]
H_*(E^*)     & \rTo^{\cong} & H_*(C^*)\\
\dTo_{\si_*} &               & \dTo_{\pi_*} & &              (*_4)\\
H_*(D^*)     & \rTo^{\cong} & H_*(\G^*)\\
\end{diagram}\\

\no we may replace in  ($*_2$)
the group $H_3(X^*,\,X^*\sm E^*)$
resp.,   $H_3(Y^*,\,Y^*\sm D^*)$
by the group $H_1(C^*)$ resp.,  $H_1(\G^*)$
to obtain a commutative diagram 

\begin{diagram}[height=2.3em,width=3em]
0 & \rTo & H_3(X^*) & \rTo^{\rho_X} & H_1(C^*) & \rTo^{\partial_*} &
H_{2}(X^*\sm E^*) & \rTo& 0\\
& & \dTo_{\si_*} & & \dTo_{\pi_*} & &\dTo_{\wr}&            & & (*_5)\\
0 &\rTo& H_3(Y^*) &\rTo^{\rho_Y} &
H_1(\G^*)& \rTo^{\partial_*}&
H_{2}(Y^*\sm D^*)& \rTo & 0\\
\end{diagram}\\

\no where $ \rho_X:=\si_*\circ t_X\circ i_*$ resp.,
$ \rho_Y:=\pi_*\circ t_Y\circ i_*$.
The diagram ($*_5$) yields the assertions of the lemma. \qed\medskip

\no\blem\la{nle2}\query{(label=nle2)}
$\rho_X (H_3(X^*))=H_C$ {\it and} $\rho_Y(H_3(Y^*))=H_{\G}$. \elem

\no  \proof  We start by constructing
an appropriate free base of the $\Z$-module
$H_3(X^*)$ (resp., $H_3(Y^*)$).
 By the Alexander duality we have isomorphisms
$$\tH_3(X^*)\cong H^2({\dot S_E})\cong H^2({\dot S^{(0)}_E})
\cong
\Z^{b_0(S_C^{(0)})}\,,$$
where $S^{(0)}_E:=S_C^{(0)}\times\C$
(similarly, ${\tH}_3(Y^*)\cong \Z^{b_0(\SG)}$).
Thus $H_3(X^*)=\tH_3(X^*)\cong H^2({\dot S^{(0)}_E})$ is a
free $\Z$-module,
and so   the universal coefficient formula
provides
yet another form of the Alexander duality:
$$H_3(X^*)\cong H^2({\dot S^{(0)}_E})\cong H_2({\dot S^{(0)}_E})\qquad
{\rm (resp.,}\quad H_3(Y^*)\cong H_2({\dot S}_D){\rm )}\,$$
(see e.g., \cite[\query{p.130,} Ch. 2, \S 15.5]{FoFu}).
Now a free base of $H_3(X^*)$ can be reconstructed as follows.
On each component $T_Q:=\si^{-1}(Q)=\{Q\}\times\C\simeq\C$
of $S_E^{(0)}$ (where
$Q\in S_C^{(0)}$) fix a point
$Q'\in T_Q$, and fix a complex 2-ball $B^2_{Q'}$ in $X$
transversal to
$T_Q$ at $Q'$. Then the 3-cycle $[s_{Q'}]\in H_3(X^*)$
which corresponds to the 2-cycle
$[{\dot T_Q}]\in H_2({\dot S^{(0)}_E})$ by
the Alexander duality, can be represented by
a small 3-sphere (say) $s_{Q'}$ in $B^2_{Q'}$
centered at $Q'$ \cite[Ch. 2, \S 17\query{, p.159}]{FoFu}.
These cycles
$[s_{Q'}]\in H_3(X^*)$ with $\pi(Q')=Q$, where $Q$
runs over $S^{(0)}_C$ form the desired base.

It is convenient to fix the choice
of a point $Q'\in T_Q$ as follows.
Consider the curve $ C'\subseteq X$ given by the equations
$ u=f\circ\s=0$  (thus $C'=E\cap D_u^{\rm red}$).
It is easily seen that the restriction $\s|C':C'\to C$
is an isomorphism.
Its inverse $C\to C'\subseteq E$ is
a section of the projection $\s|E:E=C\times\C\to C$.
 Letting $Q':=T_Q\cap C'$ and
placing the transversal ball $B^2_{Q'}$ into the
divisor $D_u^{\rm red}=\{u=0\}\subseteq X$
(notice that by \ref{smooth},
$Q'$ is a smooth point of
$D_u^{\rm red}\simeq D_g^{\rm red}$),
we let $\delta_{Q'}:=s_{Q'}\cap C'$.
Then $ \s(\delta_{Q'})=:\delta_{Q} =lk_{Q}(C)$
is a link of the curve $C$ at
the point $Q\in C$, and so $\delta_{Q}\in H_C$.
We claim that $\rho_X([s_{Q'}])=[\delta_{Q}]\in
 H_1(C^*)$,  or in other words, that
$i_*([s_{Q'}])=t_X^{-1}([\delta_{Q'}])$.

Indeed,  consider the Zariski open dense subsets
$X^{**}:=X^*\setminus {\rm sing}\,E^*\subseteq X^*$ resp.,
${\rm reg}\,E^*:=E^*\setminus {\rm sing}\,E^*\subseteq E^*$.
Let $U$  be a small tubular neighborhood
of the closed submanifold ${\rm reg}\,E^*\subseteq X^{**}$
with a retraction $\tau:U\to {\rm reg}\, E^*$.
We have $\delta_{Q'}\subseteq {\rm reg}\,E^*$,
and the class
$[\tau^{-1}(\delta_{Q'})]\in H_3(X^{**},\,X^{**}
\setminus {\rm reg}\,E^*)=
H_3(U,\,U\sm {\rm reg}\,E^*)$
corresponds to the class $[\delta_{Q'}]$ under the Thom
isomorphism $H_3(X^{**},\,X^{**}\setminus {\rm reg}\,E^*)
\stackrel{\cong}{\longrightarrow} H_1({\rm reg}\,E^*)$
\cite[Ch. 4, \S 30\query{, p. 262}]{FoFu}.
The Thom
isomorphism being functorial under open embeddings
\cite[VIII.11.5]{Do},
we have
$t_X^{-1}([\delta_{Q'}])=[\tau^{-1}(\delta_{Q'})]
\in H_3(X^{*},\,X^{*}\setminus E^*)$.

The divisors
$D_u^{\rm red}$ and $E^*$ in $X^*$ are transversal,
and so the 3-sphere $s_{Q'}$ is
transversal to the divisor $E^*$ along the 1-cycle
$\delta_{Q'}\subseteq E^*$.
Hence $s_{Q'}\cap U$ and $\tau^{-1}(\delta_{Q'})$
represent the same relative homology class in
$H_3(U,\,U\sm {\rm reg}\,E^*)=
H_3(X^{**},\,X^{**}\setminus {\rm reg}\,E^*)$.
This proves the equality
$\rho_X([s_{Q'}])=[\delta_{Q}]$.

Therefore, we have shown that $\rho_X(H_3(X^*))\subseteq H_C$,
and every element
$[\delta_Q]$ of
a free base of the
$\Z$-module $H_C$ is the image of an element
$[s_{Q'}]\in H_3(X^*)$
with $\si(Q')=Q$. This proves the first equality of the lemma;
the proof of the second one is similar. \qed\medskip

\no{\it Proof of Proposition \ref{mainprop}.}
In virtue of \ref{nle1} and \ref{nle2} we have 
$\pi_*(H_C)\subseteq H_{\G}$ and
$\hat\pi_*$ in (\ref{iso}) is, indeed, an isomorphism.
Furthermore, as $H_1 ({\ol \G}) =H_1(\G )/ G( \G^* )$ and
$H_1 (\overline {C^*}) =H_1(C )/ G( C^* )$ (see \ref{bs}),
the homomorphism ${\ol \pi}_* : H_1 (\overline {C^*}) \to
H_1 ({\ol \G})$ factors through ${\widehat\pi}_*$, whence
${\ol \pi}_*$ is a surjection. To show
that it is also injective, consider a nonzero element
$[ \alpha ] \in H_1 (\overline {C^*})$ 
generated by a 1-cycle
$\alpha$ in $C_i^*$ (where $i\le k$) such that
the 1-cycle $\beta := \pi ( \alpha )$ in $\G_i^*$
gives a zero element of $H_1 ({\ol \G }_i)$.
Hence the class $[\beta ]\in H_1 (\G_i^*)$
is contained in the subgroup $G (\G_i^*)$. 
It is easily seen that $\pi_* (G (C_i^* ) )$
contains a subgroup of finite index in $G (\G_i^*)$. 
Thus for some
$m>0$ we have $m [\beta ]= \pi_* ([\gamma ])$ 
with $[\gamma ]
\in G (C_i^*)$. As $m [\alpha ] -[\gamma ]$ 
is still a
nonzero element of $H_1 (\overline {C^*})$ 
whose image in
$H_1 (\G_i^*)$ is zero, we have
${\widehat\pi}_*(m [\alpha ] -[\gamma ])=0$, 
which is wrong since
${\widehat\pi}_*$ is an isomorphism.  
This proves that ${\ol \pi}_*$
in (\ref{isocomp}) is an isomorphism.

It follows that the decompositions (\ref{nlo}) 
and (\ref{nlo*}) 
in \ref{hc} can be chosen in such 
a way that $\pi_*$ respects them. 
Then by (\ref{iso}) and (\ref{isocomp}), $\tilde\pi_*$ in 
 (\ref{isoG}) is also an isomorphism. 
The proof is completed.
\qed\medskip

\query{(label=prim)} \no\bsit\label{prim} {\rm 
The proof of the next 
lemma is based on the following simple observation.
Let $G$ is a free abelian group of finite rank. 
For an element $a\in G$, 
the following conditions are equivalent:

\begin{enumerate} 
\item[(i)] $a$ is {\it primitive} (that is, if $a=kb$ with
$b\in G,\,\,k\in\Z$ then $k=\pm 1$);
\item[(ii)] given a free base $a_1,\ldots,a_n$ of a $\Z$-module
$G$, 
the coordinates in
the presentation $a=\sum_{i=1}^n k_ia_i$ are relatively prime.
\end{enumerate} 
}
\esit

\query{(label=lpri)} \no\blem\la{lpri} 
Let $F$ be a curve with the irreducible components 
$F_1,\ldots,F_L$, which are all non-compact, 
and with the multiple points 
$P_1,\ldots,P_M$.
Denote by $p_1,\ldots, p_N$ the punctures of the curve 
$F^*:=F\backslash \{P_1,\ldots,P_M\}$, 
and let $B_i$ be a local branch of $\bar F$ centered at $p_i$.
Consider the group $G:=G(F^*)/H_F$, where 
$$G(F^*):=\big\langle lk_{p_i}(B_i)\,\big|\,
i=1,\ldots,N\big\rangle\subseteq H_1(F^*)$$ and 
$$
H_F:=\big\langle lk_{P_j}(F)\,\big|\,j=1,\ldots,M\big\rangle\subseteq G(F^*)\,.$$ 
Then $G\cong \Z^{N-M-L}$. 
Furthermore, for each $i=1,\ldots,N$ the class
$[lk_{p_i}(B_i)]\in G$ is primitive unless it is zero.
The latter holds if and only if $p_i$ is the only puncture
of the corresponding component $F_l$, 
and $F_l$ is an isolated component of $F$.\elem

\no\proof Denoting $\nu: \nF\to F$ 
a normalization, we let
$\nu^{-1}(P_m)=\{p_j\}_{j\in J_m}$ with 
$J:=\coprod_{m=1}^M J_m\subseteq\{1,\ldots,N\}$ 
\footnote{Here 
$\coprod$ stands for the disjoint union.}. Thus 
$lk_{P_m}(F)=\sum_{j\in J_m} lk_{p_j}(B_j)$.
We also let $\{p_i\}_{i\in I_l}$ 
be the set of punctures of a component 
$F_l^*:=F_l\cap F^*$ 
of the curve $F^*$. For every $m=1,\ldots,M$ 
(resp., $l=1,\ldots,L$)
we pick an index $j_m\in J_m$ 
(resp., $i_l\in I_l\backslash J$)
(notice that the index set $I_l\backslash J$ 
of the punctures of $F_l$ is non-empty, as
$F_l$ is assumed to be non-compact).
Then we have:
$$ 
G(F^*)=\big\langle lk_{p_i}(B_i)\,\big|
\,i\not\in\{i_1,\ldots,i_L\}\big\rangle$$
$$=
\big\langle lk_{p_i}(B_i), \,\,lk_{P_m}(F)\,\big | \,m=1,,\ldots,M, \,\,\,i\not\in\{i_1,\ldots,i_L,j_1,\ldots,j_M\}\big\rangle\,.
$$
Thus $G(F^*)\cong \Z^{N-L}\cong H_F\oplus \Z^{N-L-M}$ with 
$H_F\cong \Z^M$, and so
\query{(label=bas)}\be\la{bas}\,\quad
G=G(F^*)/H_F= \big\langle [lk_{p_i}(B_i)]\,\big |  \,
i\not\in\{i_1,\ldots,i_L,j_1,\ldots,j_M\}
\big\rangle\cong \Z^{N-L-M}\,.
\ee
The elements $[lk_{p_i}(B_i)]\in G$ 
of the free base (\ref{bas})
satisfy the condition (ii) of 
\ref{prim}, whence are primitive. 
The indices $j_m\in J_m$ 
being arbitrary, as card$J_m\ge 2$ the classes 
$[lk_{p_i}(B_i)]\in G$ with
$i\in J$ are all primitive as well. By the same reason, 
a class $[lk_{p_i}(B_i)]\in G$
with $i\in I_l \backslash J$ is primitive unless 
$\{i\}=I_l \backslash J$
(i.e., unless $p_i$ is the only puncture of $F_l$).
If such a component $F_l$ 
meets another one $F_{l'}$ at a multibranch point $P_m$,
then choosing $j_m$ from $J_m\cap I_{l'}$ and 
decomposing 
$[lk_{p_i}(B_i)]\in G$ in the base (\ref{bas})
we obtain that 
at least one of the coordinates equals $-1$, hence 
\ref{prim}(ii) is fulfilled, and so $[lk_{p_i}(B_i)]$ 
is primitive too. Finally, $[lk_{p_i}(B_i)]\in G$ 
is not primitive 
if and only if $p_i$ is the only puncture
of $F_l$ and $F_l$ is isolated (in that case 
$lk_{p_i}(B_i)+\sum_{P_m\in F_l} lk_{P_m}(F)=0$ 
in $G(F^*)$, whence 
$[lk_{p_i}(B_i)]=0$ in $G$).
\qed\medskip

In the proof of Theorem \ref{3f} below
(based on \ref{mainprop} and \ref{lpri})
the role of $F$ in \ref{lpri} is played
respectively by the curves $\G$ and 
$\CNV\backslash \CV$. 
As in the proof of \ref{lpri}, 
it will be important to bear in mind 
the freedom of choice when selecting 
the indices 
$\{i_1,\ldots,i_L,j_1,\ldots,j_M\}$ 
as in (\ref{bas}). From \ref{lpri} 
we obtain such a corollary.
\medskip

\query{(label=nonram)} \no\bcor\la{nonram}
In the notation as in \ref{mainprop}, 
let $Q\in \bar C_i\,\,\,(i\in\{1,\ldots,k\})$
be a puncture 
of $C_i^*$ such that $P:=\bar\pi(Q)\in\bar\G_i$ is not 
a puncture at infinity of the affine curve $\G_i$. 
Then the following hold.
\begin{enumerate}
\item[(a)] ${\rm mult}_Q \bar\pi=1$ 
(that is, $\bar\pi$ is non-ramified at $Q$);
\item[(b)] $P$ is a puncture of $\G_i^*$;
\item[(c)] the components $C_i$ and $\G_i$ are horizontal
(that is, $1\le i\le h$).
\end{enumerate}
\ecor

\no\proof The curve $C^*_i$ possesses 
yet another puncture over a place at infinity of $\G_i$,
whence by (\ref{lpri}), $[lk_Q(C^*_i)]$ 
is a primitive element of the group 
$G(C^*)/H_C$. As 
${\tilde\pi}_*$ in 
(\ref{isoG}) is an isomorphism,  
${\tilde\pi}_*([lk_Q(C^*_i)])=\kappa 
[lk_P(\G^*_i)]\in G(\G^*)/H_{\G}\backslash \{0\}$ 
is primitive as well, whence
$\kappa={\rm mult}_Q \bar\pi=1$, which yields (a) and (b).  

Suppose that there are two different 
local branches $A$ and $A'$ of $\bar C$ 
over a local branch $B$ of 
$\bar\G$ at $P$, and let $Q,\,Q'\in \bar C$ 
be their centers.  The same argument as above shows that
${\tilde\pi}_*([lk_Q(A)])={\tilde\pi}_*([lk_{Q'}(A')])=[lk_P(B)]$,
where both $[lk_Q(A)],\,[lk_{Q'}(A')]\in G(C^*)/H_C$ are primitive. 
Hence
\query{(label=ecyc)}\no\be\la{ecyc} [lk_Q(A)]=[lk_{Q'}(A')]\neq 0\,.\ee
It is not difficult to verify that 
(in contradiction with (\ref{ecyc}))
there exists a base (\ref{bas})
which includes both $[lk_Q(A)]$ and $[lk_{Q'}(A')]$,
unless $A$ and $A'$ are 
local branches of the curve $C$ at a multibranch point 
$Q\in S_C^{(0)}$ with $\mu_Q(C)=2$. But in the latter case 
$[lk_Q(A)]=-[lk_{Q'}(A')]$, which again contradicts (\ref{ecyc}).
Therefore, there is only one local branch of $\bar C$ over $B$,
which yields (c). 
\qed\medskip

Now we are ready to prove Theorem \ref{3f}.\medskip

\no{\it Proof of Theorem \ref{3f}.} 
It is convenient to proceed first with the proof of $\,(\beta)$.
From \ref{nonram} we get: 
$$\GS\cap \SG=\emptyset=\CS\cap\SC$$ 
(indeed, $\G^*_{\rm slant}$ resp., $C^*_{\rm slant}$
cannot have punctures other 
than the places at infinity). 
Thus the slanted components $C_i$  and $\G_i$ 
($h+1\le i \le k$)
of $C$ resp., $\G$ are isolated and 
do not contain multibranch points. 
In particular, $C_i=C_i^*$ resp., $\G_i=\G_i^*$, and 
these curves are homeomorphic to their normalizations. 
Furthermore, the group $G(C_i^*)$ (resp., $G(\G_i^*)$)
is a direct summand of $G(C^*)/H_C$ (resp., 
$G(\G^*)/H_{\G}$). Hence by \ref{mainprop},
\query{(label=sta)}
\be\la{sta} {\tilde\pi}_*|_{G(C_i^*)} : 
G(C_i^*)\stackrel{\cong}{\longrightarrow} 
G(\G_i^*)\ee
is an isomorphism, as well as 
\query{(label=stasta)}
\be\la{stasta} {\bar\pi}_*|_{H_1(\bar C_i)} : 
H_1(\bar C_i)\stackrel{\cong}{\longrightarrow} H_1(\bar \G_i)\,.\ee
Since $\deg\,{\bar\pi}|_{\bar C_i}>1$, 
(\ref{stasta}) shows that 
$H_1(\bar C_i)=0=H_1(\bar \G_i)$
i.e., both $\bar C_i$ and $\bar\G_i$ are rational curves. 
Therefore by 
(\ref{di}) in \ref{bs},
$$H_1(C_i^{\rm norm})\cong H_1(C_i)=H_1(C_i^*)=G(C_i^*)$$ and 
$$H_1(\G_i^{\rm norm})\cong H_1(\G_i)=H_1(\G_i^*)=G(\G_i^*)\,.$$
Now (\ref{sta}) yields that 
$$\np_* : H_1(C_i^{\rm norm}) 
\to H_1(\G_i^{\rm norm})$$
is an isomorphism. By a theorem of Hurwitz
(see e.g., \cite[I.2.1]{LiZa}), a morphism
$\rho\,:\,F\to G$
of smooth irreducible affine curves is an
isomorphism once the induced homomorphism
$\rho_*\,:\,H_1(F)\to H_1(G)$ is an isomorphism,
unless $H_1(F)= H_1(G)=0$
i.e., $F\simeq G\simeq \C$. 
Thus $C_i^{\rm norm}\simeq \G_i^{\rm norm}\simeq \C$, and so
the curves $C_i$
and $\G_i$ are homeomorphic to $\C$, which proves $(\beta)$.

$\,(\alpha)$ By \ref{ps} $\pi^{-1}(\GV)\subseteq\CV$, whence 
$\pi(\CH\backslash\CV)\subseteq \GH\backslash\GV$.
To show that actually,
the latter inclusion is an equality, 
suppose on the contrary that
for a point $P\in \GH\backslash\GV$,
$\pi^{-1}(P)=\emptyset$. Then all $\mu:=\mu_P(\G)$ 
branches $A_1,\ldots,A_{\mu}$
of $\bar C$ over the branches $B_1,\ldots,B_{\mu}$
of $\G$ at $P$ have centers at infinity. 
By \ref{nonram}(b)
the curve $\G^*$ has a puncture over $P$. Thus
$P\in S_{\G}^{(0)}$ is a multibranch point of 
$\GH\backslash\GV$, and so  
the primitive classes (see \ref{lpri})
$[lk_{P}(B_j)]\in G(\G^*)/H_{\G}\,\,\,(j=1\ldots,\mu)$ 
are subjected to the (only) relation
$\sum_{j=1}^{\mu} [lk_{P}(B_j)]=0$. 
Therefore, the primitive classes $[lk_{Q_j}(A_j)]\in G(C^*)/H_{C}\,\,\,(j=1\ldots,\mu)$ are also related by 
\query{(label=3sta)}\be\la{3sta}
\sum_{j=1}^{\mu} [lk_{Q_j}(A_j)]=0\,,\ee
where $Q_j\in \bar C$ is the center of the branch $A_j$. 
Constructing a free base (\ref{bas})
of the free $\Z$-module 
$G(C^*)/H_{C}$, we may
suppose that for any $j=1,\ldots,\mu$, 
$j\not\in \{i_1,\ldots,i_L\}$ (and definitely, 
$j\not\in \{j_1,\ldots,j_M\}$ as $Q_j$ 
is a puncture at infinity
of $C$). Thus the classes 
$[lk_{Q_j}(A_j)]\,\,\,(j=1\ldots,\mu)$
make a part of a free base (\ref{bas}), 
and so cannot satisfy (\ref{3sta}), 
a contradiction. 

Denote $F:=\CH\backslash \CNH=\CH\backslash\CV$ and 
$G:=\GH\backslash\GNH=\GH\backslash\GV$.
We have shown that the morphism $\pi_{|F} : F\to G$ 
of degree 1
is bijective. It follows that it is an isomorphism.
Indeed, as $D_g$ and $D_f$ meet transversally
along $F$ (see \ref{Ci}), 
the curve $F$ is the zero divisor of the
restriction $g|_{G\times \C_z}$, which 
is a polynomial of
degree one in $z$, say, $az+b\in \C[G][z]$.
Note that $a$ and $b$ have no common zero on $G$ 
(otherwise the zeros of $g|_{G\times \C}=az+b$ 
would contain 
a vertical component). As $\pi|_{F}$ is surjective,
$a$ is nowhere zero, whence
$z=-b/a\in \C[G]$, and so $\pi|_{F}$ is an isomorphism.
This proves $\,(\alpha)$.

 $\,(\gamma)$ The group 
$G(\G^*)/H_{\G}=\tilde\pi_*(G(C^*)/H_C)$ being 
generated by the classes 
$$\tilde\pi([lk_{Q_j}(C^*_i)])\in 
G(\G^*_{\rm non-vert})\big/[H_{\G}\cap 
G(\G^*_{\rm non-vert})]\,$$
(where $Q_j$ runs over the set of punctures of the curve $C^*$)
we have
$$G(\G^*_{\rm non-vert})\big/[H_{\G}\cap G(\G^*_{\rm non-vert})] = G(\G^*)/H_{\G}\,.$$ 
Since in virtue of $(\beta)$, $G(\G^*_{\rm slant})=0$ we obtain
$$G(\G^*)/H_{\G}=G(\G^*_{\rm non-slant})/H_{\G}=
G(\G^*_{\rm horiz})\big/[H_{\G}\cap 
G(\G^*_{\rm horiz})]\,.$$
In particular, for each puncture $P$ of $\G^*_{\rm vert}$, 
we have
\query{(label=incl)}\be\la{incl}
[lk_P(\bar\G)]\in G(\G^*_{\rm horiz})\big/
[H_{\G}\cap G(\G^*_{\rm horiz})]\,.
\ee
It follows that on any vertical component of 
$\G$ there is only one puncture
at infinity. Indeed,
if there were a component $\G_i$ of $\GV$ 
with at least two punctures at infinity, 
say, $P_1$ and $P_2$, 
then any  free base as in (\ref{bas}) of the 
$\Z$-module $G(\G^*)/H_{\G}$ would contain 
at least one of the corresponding classes, say, 
$[lk_{P_1}(\bar\G)]$, which contradicts 
(\ref{incl}). 

A similar argument shows that the curve $\GV$ 
has no selfintersection. In particular, 
the components of $\GV$ are disjoint, and for 
each $i=k+1,\ldots,n$,
$\G_i$ is homeomorphic to 
$\G_i^{\rm norm}$. 

Furthermore, by \ref{mainprop}
$$\bar\pi_* : H_1({\ol C^*})\to H_1(\bar\G)
=H_1({\bar \Gamma}_{\rm non-vert})\oplus 
H_1({\bar \Gamma}_{\rm vert})$$
is an isomorphism, and 
$\bar\pi_*(H_1({\ol C^*}))\subseteq
H_1({\bar \Gamma}_{\rm non-vert})$. 
It follows that 
 $ H_1({\bar \Gamma}_{\rm vert})=0$,
whence the components of $\GV$ are rational. 
Finally, $\G_i\,\,\,(i=k+1,\ldots,n)$ is a rational curve 
with one place at infinity and without selfintersections,
therefore is homeomorphic to $\C$. 
 
As for any $i=h+2,\ldots,n$, $D_{h+1}\cap D_i=\emptyset$ 
and 
$D_i\simeq \G_i\times\C$ is simply connected, 
with the notation as in \ref{def} we have that
the restriction 
$f_{h+1}|_{D_i}$ does not vanish, and so is constant: 
$f_{h+1}|_{D_i}=:\lambda_i\in\C$. Thus
we obtain a decomposition as in $(\gamma)$, 
which completes the proof of $(\gamma)$. 

 $\,(\delta)$ As the matrix $B'$ in \ref{def} 
is unimodular and by $\,(\gamma)$ the vertical components are disjoint, 
arguing as in
\ref{def} we can easily see that the morphism $\pi:C\to\G$ 
maps any component
of $\CV$ into a component of $\GV$, and maps different 
components
of $\CV$ into different components of $\GV$.
Moreover, the columns of $B'$ are vectors of the standard basis
in $\R^{n-k}$. Hence up to a reordering we may assume that
$B'$ is the unit matrix. 
The slanted components $\G_{h+1},\dots,\G_k$
being isolated the last $k-h$ lines of 
the matrix $B$ are zero, 
which completes the proof of $(\delta)$.\qed\medskip

\query{(label=rmq)} \brem\la{rmq} {\rm Notice that 
for any vertical component $C_i$ of $C$ ($i=k+1,\ldots,n$),  
$P_i:=\pi(C_i)$ is a smooth point of the curve $\G_i$.
This follows from \ref{gone} as $m_{ii}=1$. }\erem

\query{(label=pracs)} \subsection{Preserving acyclicity: 
a criterion}\la{pracs} 
We have the following partial converse to \ref{3f}.

\query{(label=3g)}\no\bthm \label{3g} 
Let $X$ and $Y$ be irreducible smooth affine 3-folds 
which satisfy the condition (ii) of \ref{conv0}, 
as well as the following one 
(which replaces \ref{conv0}(iii)):

\begin{enumerate}
\item[(iii$'$)] $Y=Z\times \C$ 
is a cylinder over a smooth acyclic 
affine surface $Z$, and
$D=\G\times\C$ with $\G\subseteq Z$. 
\end{enumerate}

\no Suppose also that the conditions  
($\alpha$)-($\delta$) of Theorem \ref{3f} 
are fulfilled. \footnote{See \ref{UR} below.}
Then the induced homomorphisms
$$\s_*\,:\,H_*(X)\to H_*(Y)\quad {\rm and} 
\quad \si_*:\pi_1(X)\to\pi_1(Y)$$ 
are isomorphisms.
Henceforth, the 3-fold $X$ is acyclic; it is contractible 
if and only if $Y$ is so. \ethm

\no\proof The theorem immediately follows from \ref{pi1}, \ref{prdec} 
and \ref{newla} below.
Indeed, the assumption \ref{3f}($\delta$) allows to apply 
\ref{pi1} in order to get that $\si_*:\pi_1(X)\to\pi_1(Y)$
is an isomorphism. In turn, by \ref{newla} the assumptions of 
\ref{prdec} are fulfilled, and so by \ref{prdec}, 
$\s_*\,:\,H_*(X)\to H_*(Y)$ is an isomorphism as well. \qed\medskip

From \ref{3f} and \ref{3g} we obtain the following criterion 
for preserving the acyclicity.

\query{(label=c3g)}\no\bcor \label{c3g}
Let $X$ and $Y$ be smooth affine 3-folds satisfying 
the conditions \ref{conv0}($ii$)  and \ref{3g}($iii'$).
If $Y$ is acyclic resp., contractible 
then $X$ is so if and only if the conditions  
\ref{3f}($\alpha$)-($\delta$) are fulfilled. \ecor

\query{(label=UR)}\no\brem \label{UR} 
{\rm Assuming in \ref{3g} that the conditions 
\ref{3f}($\alpha$)-($\delta$)
are fulfilled we require implicitly  
that the things are as in \ref{def}-\ref{factor},
without supposing acyclicity as in \ref{conv0}(i).
In fact, \ref{def}-\ref{factor} refer only to the fact
that the multiplicity matrix $\M$ is unimodular, which is 
anyhow foreseen by \ref{3f}($\delta$). 
 
Actually the proofs of the lemmas below 
rely only on the conditions
(ii), (iii) and the following one.

\begin{enumerate}
\item[(iv)]  
There is a regular function $\varphi\in A=\C[Y]$ 
such that
for every $i=h+1,\ldots,n$ the restriction
$\varphi|_{D_i}$ is constant, 
as well as the restriction of $\varphi$
to any fiber of the morphism $\pi:D\to\G$,
and for any point $P\in \SG\backslash \GV$, both
$\Phi_P:=\varphi^*(\varphi(\pi^{-1}(P)))\subseteq Y$ and
$\Psi_P:=\si^{-1}(\Phi_P)\subseteq X$
are smooth, reduced surfaces with $\si^*(\Phi_P)=\Psi_P$.
\end{enumerate}

\no Under these assumptions $\s_*$ 
is an isomorphism in homology (even if we do not  
suppose as in \ref{3g} that $Y$ is acyclic). }\erem

  The next lemma shows that the conditions 
(ii) and (iii$'$) imply (iv).
 
\query{(label=oldla)}\blem \label{oldla} Under the assumptions
of \ref{3g}
there exists a (possibly empty) 
smooth, reduced curve $\G'$ on $Z$ such that 
 
\begin{enumerate}
\item[(a)] $\G'\cap\GNH=\emptyset$ and 
$\G'\supseteq S_{\G}^{(0)}\backslash
\GNH$,
\item[(b)] the surfaces $\Phi:=\G'\times\C\subseteq Y$ and 
$\Psi:=\si^{-1}(\Phi)\subseteq X$ are smooth, and
\item[(c)] $\s|_{\Psi} : H_*(\Psi)\to H_*(\Phi)$ is an isomorphism.
\end{enumerate}\elem

\no\proof If $h<n$ 
(i.e., $\GNH\neq\emptyset$) then we 
take for $\G'$ the union of the fibers
of the regular function $f_{h+1}|_{Z}$ 
through the points of $S_{\G}^{(0)}\backslash
\GNH$.
Indeed, since  
the Euler characteristics 
$e(\G_{h+1})=e(Z)=1$ and $\G_{h+1}=
(f_{h+1}|_{Z})^*(0)$, by \cite[6.2]{Za1} 
the fibers of $f_{h+1}|_{Z}$
are smooth, reduced and irreducible except 
for, possibly, $\G_{h+1}$. Thus 
$\G'$ is smooth and reduced.  
As $f_{h+1}$ is constant on any component $\G_i$ of $\GNH$
we have 
$\G'\cap\GNH=\emptyset$,
whence (a) is fulfilled. 

In the case where $h=n$ (i.e., $\GNH=\emptyset$)
the existence of a smooth, reduced curve $\G'$ satisfying (a) 
easily follows by Bertini's Theorem. 

Evidently, $\Phi=\G'\times\C$ is a smooth surface. 
Since $\s|_{X\backslash E} : X\backslash E \to Y\backslash D$ 
is an isomorphism,  
the surface $\Psi:=\si^{-1}(\Phi)$ is smooth if
it is smooth 
at the points $R\in \Psi\cap E$. 
Let $Q=(P,z_0):=\s (R) \in C$ 
with $P\in \G\cap\G'\subseteq \GH\backslash\GNH$.
As $\G'$ is smooth we can find local 
coordinates $(x,y)$ on $Z$ 
centered at $P$ such that (locally) $\G'=\{x=0\}$. 
Thus $Q=(P,z_0)=(0,0,z_0)$, $\,R=(Q,u_0)=(0,0,z_0,u_0)$,
and (locally) $$\Psi=\{x=0\}\cap X=\{f(0,y)u-g(0,y,z)=0\}\subseteq\C^3_{y,z,u}\,,$$ 
where by the condition \ref{3f}($\alpha$),
(locally) $g(0,y,z)=a(y)z+b(y)$ with $a(0)\neq 0$
(see the proof of \ref{3f}($\alpha$)). Therefore, 
$\p g/\p z (Q)=a(0)\neq 0$, whence the surface $\Psi$ 
is smooth at $R$. Thus (b) holds. 

To show (c) we need the following claim. 

\medskip

\no {\bf Claim.} $(\s|_{\Psi})^*
\big((D_{\rm horiz}\cap \Phi)^{\rm red}\big)=
( E_{\rm horiz}\cap \Psi)^{\rm red}.$

\medskip\no {\it Proof of the claim.} 
We will use the same local chart and the notation as above.
We have $(D_{\rm horiz}\cap \Phi)^{\rm red}=\{P\}\times\C_z=y^*(0)$
($y$ being regarded as a function on $\Phi$). 
Then locally, $(y\circ\s) : (y,z,u)\longmapsto y$,
whence the subspace
ker$\,d(y\circ\s)|_R=\{y=0\}$ 
does not contain the tangent space
$T_RX$ (indeed as $a(0)\neq 0$ the differential 
$d(fu-g)|_R=...-a(0)dz$
of the defining polynomial of $\Psi$ in $\C^3_{y,z,u}$
at the point $R$ is
not proportional to the differential
$dy$ of $y= y\circ\s$ at $R$).
It follows that $\,d(y\circ\s)|_{\Psi}$ does not vanish
at the points $R\in \Psi\cap E$, and so locally
$$( E_{\rm horiz}\cap \Psi)^{\rm red}=(y\circ\s)^*(0)=
(\s|_{\Psi})^*\big((D_{\rm horiz}\cap \Phi)^{\rm red}\big)\,,$$ 
as desired. \qed\medskip

Notice that $\Phi\cap D_{\rm horiz} = 
(\G'\cap\GH)\times\C\subseteq \Phi$
is a disjoint union of affine lines, whereas in virtue of
\ref{3f}($\alpha$), 
$\Phi\cap\CH$ is a finite set of points, one 
on every of those lines. Furthermore, 
$\Psi\cap  E_{\rm horiz}=(\Phi\cap\CH)\times\C$. 

The proof of (c) is based on \ref{prhom}. In the notation as in 
\ref{prhom} we let
${\hat X}:=\Psi,\,{\hat Y}:=\Phi,\,\hat E:=\Psi\cap E_{\rm horiz}$
and $\hat D:=\Phi\cap D_{\rm horiz}$. 
Then by (b) above, ${\hat X}, {\hat Y}, \hat E$ and $\hat D$ 
are smooth varieties,
$\si(\hat E)\subseteq\hat D$, 
$\si({\hat X}\backslash\hat E)={\hat Y}\backslash\hat D$,
and (by the above claim) $\s^*(\hat D)=\hat E$, 
that is, the condition (i)
of \ref{prhom} is fulfilled. 
As $\si|_{{\hat X}\backslash\hat E}:{\hat X}\backslash\hat E\to 
{\hat Y}\backslash\hat D$ is an isomorphism, 
taking into account the observations following the claim
it is easily seen that \ref{prhom}(ii) holds as well. 
Thus by \ref{prhom}, $\s_* : H_*({\hat X})\to H_*({\hat Y})$ 
is an isomorphism, which yields (c). \qed\medskip

\query{(label=newnot)}\no\bsit\la{newnot}  {\rm 
The proof of the next lemma relies on \ref{prdec},
where we let $\hE = E'\cup E''$ and $\hD=D'\cup D''$
with
$E':=E_{\rm horiz},\,\,D':=D_{\rm horiz}$,
$$E'':=E_{\rm non-horiz}\coprod \Psi\quad {\rm
and}\quad
D'':=D_{\rm non-horiz}\coprod \Phi\,,$$
$\Phi$ and $\Psi$ being the same as in \ref{oldla}(b).
}\esit

\query{(label=newla)}\blem \label{newla}  
Under the assumptions and the notation 
as in \ref{3g} and \ref{newnot}, the conditions 
($i'$) and ($ii'$) of \ref{prdec} are fulfilled. 
\elem

\no\proof In virtue of \ref{3f}($\gamma$),
for every $i=h+1,\ldots,n$ 
both curves $C_i$ and $\G_i$ are homeomorphic to $\C$, 
whence both surfaces $E_i=C_i\times\C$ and 
$D_i\simeq \G_i\times\C$ 
are homeomorphic to $\C^2$.
Therefore $E_{\rm non-horiz}$ and $D_{\rm non-horiz}$
are topological manifolds, and
$(\si|_{E_{\rm non-horiz}})_*:H_*(E_{\rm non-horiz})\to 
H_*(D_{\rm non-horiz})$ is an isomorphism.
As the surfaces $\Phi$ and $\Psi$ are smooth 
we have that $E''$ and $D''$ are topological 
manifolds as well, and in view of \ref{oldla}(c),
$(\s|_{E''})_*: H_*(E'')\to H_*(D'')$ is an isomorphism.

As by our construction 
$$\SG\subseteq \G'\cup \GNH\quad {\rm and}\quad
\SC\subseteq \pi^{-1}(\G')\cup \CNH\,,$$ 
the curves 
$\GH\backslash (\G'\cup \GNH)$ and 
$\CH\backslash (\pi^{-1}(\G')\cup \CNH)$
have no multibranch points. Therefore they are 
topological 
manifolds, as well as the surfaces
$$D'\backslash D''=[\GH\backslash 
(\G'\cup \GNH)]\times\C\,\,\, {\rm and}\,\,\,
E'\backslash E''=[\CH\backslash 
(\pi^{-1}(\G')\cup \CNH)]\times\C\,.$$
By \ref{3f}($\alpha$) the projection 
$$\pi : \CH\backslash (\pi^{-1}(\G')\cup \CNH)\to
\GH\backslash (\G'\cup \GNH)$$ is an isomorphism, 
whence $\si|_{E'\backslash E''}: E'\backslash E''\to
D'\backslash D''$ is so as well. Now the conditions 
(i$'$) and (ii$'$) of \ref{prdec} follow. \qed\medskip


\query{(label=SAM3)}
\section{ Simple affine modifications of $\C^{3}$ diffeomorphic to $\R^6$}\la{SAM3}

The main result of this section is \ref{tra},
which together with \ref{c3g} provides
a criterion for as when a simple affine 
modification $X\subseteq\C^4$ 
of $Y=\C^3$ is isomorphic to $\C^3$. 

\query{(label=subEXO)}
\subsection{Exotic simple modifications of $\C^{3}$}
\la{subEXO}
We keep the terminology and the notation of section \ref{S2}, 
and we adopt the following

\query{(label=conv3)}\no\bconv\la{conv3} {\rm   Hereafter
\begin{enumerate}
\item[(i)] $Y=\C^3$ with the coordinates $x,y,z$, and
\item[(ii)] $X$ is a smooth affine 3-fold in $\C^4$ 
diffeomorphic to $\R^6$, with equation of the form
\query{(label=polp)}
\be\la{polp} \,\quad p=f(x,y)u+g(x,y,z)=0\,,\ee
where $f\in\C[x,\,y]\sm\C,\,\, g\in\C[x,\,y,\,z]$
(thus by \ref{c3g}, 
the conditions \ref{3f}($\alpha$)-($\delta$) hold).
\end{enumerate}
}\econv

\query{(label=sam)}\no\bsit\la{sam} 
{\rm Note that the blowup morphism 
$\si: X\ni (x,y,z,u)\longmapsto (x,y,z)\in Y$
represents $X$ as a simple affine modification of $Y=\C^3$ 
along the cylindrical divisor 
$D=D_f^{\rm red}= \G\times\C$ (where $\G:=f^{-1}(0)\subset\C^2$)
with center $C=\{f=g=0\}\subset \C^3$
and with the morphism $\pi: C\to\G$ as in \ref{Ci} given by
$\pi: (x,y,z)\longmapsto (x,y)$. 
Hence the assumptions (i)-(iii) of \ref{conv0} are fulfilled. 
}\esit

\query{(label=fct)}\no\bsit\la{fct} {\rm 
As in \ref{factor} we factorize $f\in\C[x,y]$ 
into irreducible factors:
$f=\prod_{i=1}^n {f_i}^{a_i}$, 
and we write $g$ 
as a polynomial in $z$:
$$ g(x,y,z)=\sum_{j=0}^d b_j(x,y)z^j\quad 
{\rm with}\quad d:=\deg_zg\,.$$
We let $d_i:=\deg\,(\pi|_{C_{i}})$. 
Recall (\ref{def}, \ref{non-smth}) that 
an irreducible component $C_i$ 
(resp., $\G_{i}=f_i^{-1}(0)$) 
of the curve $C$ (resp., $\G$)
is vertical
(resp., horizontal resp., slanted) 
if and only if $d_i=0$
(resp., $d_i=1$  resp., $d_i\geq 2$). 
The following lemma is a simple observation, 
and so we omit the proof. 
}\esit

\query{(label=vhs)}\no\blem\la{vhs}  
For every  
$i=1,\dots,n$ we have 
$b_j\in (f_i)\,\,\forall j=d_i+1,\ldots,d$ 
and $b_{d_i}\not\in (f_i)$. Hence $p$ as in (\ref{polp})
admits a unique presentation
$$ p=fu+g_i+f_{i}h_i$$
with $g_i(x,y,z):=\sum_{j=0}^{d_i} b_j(x,y)z^j\in\C[x,y,z]$
and $h_i\in (z^{d_i+1})\subseteq\C[x,y,z]$. 
Furthermore, if $\fh\neq {\rm const}$
then $p$ can be written as
\query{(label=prp)}\be\la{prp} 
p=\fh\fnh u + b_0+b_1z+f_{\rm horiz}^{\rm red} h_0\ee
with $b_0,\,b_1\in\C[x,y]$ and $h_0\in (z^2)\subseteq\C[x,y,z]$, where
$b_1|_{\GH\backslash\GNH}$ has no zero.
\elem

\query{(label=EXO)}\no\bsit\la{EXO} {\rm 
Recall \cite{Za2} that
an {\it exotic $\C^3$} is a smooth affine 3-fold 
diffeomorphic to
$\R^6$ but non-isomorphic to $\C^3$. }\esit

The principal result of this subsection is the following theorem. 

\query{(label=tra)}\no\bthm\la{tra} 
If under the assumptions 
as in \ref{conv3}, at least one of the curves $\CNH$ and $\GNH$
is singular then $X$ is an
exotic $\C^3$.\ethm

The proof is done in \ref{exo1} and \ref{exo2} below.
For a converse result, see \ref{mthsec5}
in the next subsection. 

Notice that \ref{tra}
provides a regular way of constructing exotic $\C^3$-s
as hypersurfaces in $\C^4$. Let us give concrete examples. 

\query{(label=Exa)}\no\bsit\la{Exa} {\bf Examples.} {\rm
For the Russell cubic $3$-fold
$$X=\{x^2u+x+y^2+z^3=0\}\subset \C^4\,,$$
the curve 
$\G=\GS=\{x=0\}\subset \C^2$ is smooth and
isomorphic to $\C$,
whereas 
$C=\CS=\{x=y^2+z^3=0\}\subset \C^3$ (the center of modification)
is homeomorphic to $\C$ but singular. 
It is well known \cite{ML1, De, Za2} that $X$ represents an exotic
algebraic structure on $\C^3$. By \cite{KaML1} the  
Koras-Russell cubic $3$-fold
$$X=\{(x^2+y^3)u+x+z^2=0\}\subset \C^4\,$$ 
also is an exotic $\C^3$; here
both $\G=\GS=\{x^2+y^3=0\}\subset \C^2$ and 
$C=\CS=\{x^2+y^3=0=x+z^2\}\subset\C^3$ are
homeomorphic to $\C$ but singular.
}\esit

\query{(label=rmN)}\no\brem\la{rmN}{\rm
Generalizing a theorem of Sathaye \cite{Sat1},
in \cite[Thm. 7.2]{KaZa1}
it is proven that actually, every smooth acyclic surface in $\C^3$
with equation $p=f(x,y)u+g(x,y)=0$ (where $f,\,g\in \C^{[2]}$) 
is isomorphic to $\C^2$ and rectifiable.
Examples \ref{Exa} show that in general, 
this does not hold anymore in $\C^4$
without the additional
assumption of smoothness of $ \CNH$ and $\GNH$
(cf. \ref{tra} above). }\erem

The proof of \ref{tra} starts 
with the following proposition (cf. \ref{MLN} below).

\query{(label=exo1)}\no\bprop\la{exo1}
Let $X$ and $Y$ be as in \ref{conv3}. 
If the curve $\GNH$  
is singular then the Derksen invariant
${\rm Dk}\,(A')$ of the algebra $A':=\C[X]$ 
is non-trivial, whence $X\not\simeq \C^3$.
\eprop 

The proof is done in \ref{lemGNHsing}-\ref{l4}. 

\query{(label=lemGNHsing)}
\no\blem\la{lemGNHsing} 
Under the assumptions as in \ref{exo1},
choosing appropriate new coordinates 
in the $(x,y)$-plane
and rescaling the $z$-coordinate we
may write the polynomial $p$ in the form
\query{(label=form)}\be\la{forme}
p=(x^k-y^l)^m\fh(x,y)u+z^{e}+
g_0(x,y,z)+(x^k-y^l)h_0(x,y,z)\,,
\ee
where $k,l,{e}\ge 2,\,\,(k,\,l)=1,\,\,
\fh\in \C[x,y]$ and $\quad \fh(0,0)=1,\,\,
g_0,\,h_0\in\C[x,y,z]$, ${\rm deg}_z g_0<{e}$ 
and $z^e|h_0$. \elem

\no \proof In the notation as in 
\ref{def}-\ref{isol} we have
$h<n$. In virtue of the condition 
\ref{3f}($\gamma$)
we may suppose that the component 
$\G_{h+1}$ of $\GNH$ 
is a singular 
plane curve homeomorphic to $\C$. 
Hence by the Lin-Zaidenberg 
Theorem
\cite{LiZa}, choosing new coordinates 
in the $(x,y)$-plane 
we may assume
that $f_{h+1}=x^k-y^l$ with $k,\,l\ge 2$ 
and $(k,\,l)=1$. 
As the other fibers $x^k-y^l=c \,\,\,
(c\in\C\backslash \{0\})$ 
of the polynomial 
$f_{h+1}$ are not
homeomorphic to $\C$, in view of 
\ref{3f}($\gamma$) 
we have $h+1=n$, and so
$\fnh=f_n^{m}=(x^k-y^l)^{m}$ with $m:=a_n$. 
By
\ref{vhs} the polynomial (\ref{polp}) 
can be now written as follows:
$$
  p=(x^k-y^l)^m\fh(x,y)u+b_{e}(x,y)z^{e}+
g_0(x,y,z)+
  (x^k-y^l)h_0(x,y,z)\,
$$ with $b_e\in\C[x,y],\,\,g_0,\,h_0\in\C[x,y,z],
\,\,\deg_z g_0<e$ and $z^{e+1}|h_0$. 
 In virtue of \ref{3f}($\delta$) 
and \ref{Ci} 
the divisors $D_n=\G_n\times\C$ and 
$D_g^{\rm red}$ 
meet transversally at general points of 
$C_n=D_n\cap D_g^{\rm red}$. 
If $\G_n$ were
vertical (i.e., $g(x,y,z)\equiv b_0(x,y)\, 
\mod (x^k-y^l)$; see \ref{fct}, \ref{vhs}) 
then
the equation $b_0(t^l,t^k)=0$ would have 
a unique solution $t_0$ 
such that (in virtue of \ref{rmq})
$\pi(C_n)=(t_0^l,t_0^k)$ is a smooth point 
of $\G_n\subseteq\C^2_{x,y}$ 
i.e., $t_0\neq 0$. 
Thus
$b_0(t^l,t^k)=c(t-t_0)^r$ which is wrong as the
derivative of the left hand side vanishes at $t=0$. 
Therefore,
$\G_n=\GS$, whence $e=\deg \pi|_{C_{n}}\ge 2$ 
and $\pi|_{C_{n}}\,:\,C_{n}\to\G_{n}$
is a proper morphism (as each of these curves 
has only one puncture). 
It follows that the restriction
$b_{e}|_{\G_{n}}$ has no zero. As $\G_{n}$
is  simply connected
we have $b_{e}|_{\G_{n}}={\rm const}=:
b_{e}^0\in\C^*$, and so
$b_{e}(x,y)=b_{e}^0+(x^k-y^l)c_{e}(x,y)$ 
for a certain polynomial
$c_{e}\in\C^{[2]}$. Finally, rescaling 
the $z$-coordinate if
necessary, we  obtain the desired presentation. 
As by
\ref{3f}($\delta$), $\G_n=\GS$ 
is an isolated component
of $\G$, the restriction $\fh|_{\G_{n}}$ 
does not vanish, whence is constant, 
and we
may assume in addition that this constant 
is $\fh(0,0)=1$.  \qed\medskip

\query{(label=ONT)}\no \bsit\la{ONT}{\rm
We choose the weight degree function $d$
on the algebra $A$ with 
$$d_x=-lN,\,\,d_y=-kN,\,\,d_z=\sqrt{2}
\quad {\rm and}\quad d_u=klmN+e \sqrt{2}\,,$$ 
so that the polynomials
$$f_n=f_{h+1}=x^k-y^l\qquad {\rm and}\qquad
{\h p}:=(x^k-y^l)^m u+z^{e}$$
are $d$-quasihomogeneous.
Letting $N\in\N$ to be sufficiently large
(thus $d(h_0)<klN$, and hence $d((x^k-y^l)h_0)<0$) 
we may assume that
$\h p$ as above is the $d$-principal part 
of the polynomial $p$
as in (\ref{forme}).
}\esit

\query{(label=13)}\no\blem\label{l3} 
Let $q\in\C[x,y,z,u]\setminus\C$ be an
irreducible $d$-homogeneous polynomial 
with ${\rm deg}_z q < e$.
Then $q$ coincides (up to a constant factor) 
with one of the
following polynomials:
$$x,\,y,\,z,\,u,\,\lambda x^k+\mu y^l,\quad 
{\rm where}\quad\lambda,\,\mu \in\C^*\,.$$
\elem

\no\proof Letting $q=\sum_{i=0}^{e-1} a_i(x,y,u)z^i$ 
with
$ a_i\in\C[x,y,u]$ ($i=0,\ldots,e-1$) we claim 
that
there can be only one nonzero coefficient $a_i$.
Assuming on the contrary that $a_i,\,a_{i+j}\neq 0$ 
for some
$i,\,i+j\in\{0,\ldots,e-1\}$ with $1\le j\le e-1$ 
we would have
$d(a_iz^i)=d(a_{i+j}z^{i+j})$. Taking into account 
the equality
$d_y={k\over l}d_x$ we derive
$$
d(a_iz^i)-d(a_{i+j}z^{i+j})=jd_z=j\sqrt{2}=\alpha 
d_u+\beta d_x =
\alpha e \sqrt{2}+\beta' d_x$$ 
\query{(label=e1)}\be\la{e1}\Rightarrow 
(j-\alpha e)\sqrt{2}=\beta' d_x\in\Q\,\ee
for certain $\alpha\in\Z_{\ge 0}$ and 
$\beta,\,\beta'\in\Q$. 
But $\alpha e-j\in\Q\setminus\{0\}$, whence
(\ref{e1}) leads to a contradiction.

Therefore, the irreducible polynomial $q$ 
must coincide
(up to a constant factor) either with $z$ 
or with a polynomial
$a_0\in\C[x,y,u]$. Let
$q=a_0(x,y,u)=\sum_{i=0}^{t} c_i(x,y)u^i$ with
$c_i\in\C[x,y]$ ($i=0,\ldots,t$). The weights
$d_x$ (resp., $d_y$) and $d_u$ being independent 
over $\Q$ 
the same argument as above shows that at most 
one summand 
$c_i(x,y)u^i$ may be different from zero. 

Thus once again, the irreducible polynomial 
$q$ coincides
(up to a constant factor) either with $u$ or 
with a polynomial
$c_0\in\C[x,y]$. In the latter case, being 
$d$-homogeneous
the polynomial
$q$ must coincide (up to a constant factor) 
with one of the polynomials
$x,\,y,\,\lambda x^k+\mu y^l$ 
($\lambda,\,\mu \in\C^*$), 
as stated.
\qed\medskip

\query{(label=AA)}\no\bsit\la{AA} 
{\rm Let an irreducible polynomial 
$p\in\C^{[4]}$
be as in (\ref{forme}). Set as before
$X=p^{-1}(0)\subset\C^4$,
$\h X = {\h p}^{-1}(0)\subset\C^4$ and 
$\h {A'} =\C[\h X ]$
where $\h p$ is the $d$-principal part 
of $p$.
Denote by
$\rho_{\h {A'}}$ the canonical surjection
$$\rho_{\h {A'}}\,:\,\C^{[4]}\to 
\C^{[4]}/(\h p)=\h {A'}\,$$
and by $\h x=\rho_{\h {A'}}(x),\ldots, 
\h u = \rho_{\h {A'}}(u)\in \h {A'}$
the traces on $\h X$ of the
coordinate functions.

The polynomial $\h p$ being irreducible,
by \ref{KML2}(a)
$\h {A'}$ is just the graded algebra
associated with the filtered algebra $A'$,
a filtration being provided 
by the degree function $d$
on $A'$. 
}\esit

\query{(label=l4)}\no\blem\label{l4} 
In the notation as in \ref{AA}  
we have
$${\rm ML_{gr}}\,(\h {A'}) =
{\rm Dk_{gr}}\,(\h {A'}) = 
\C[\h x,\,\h y ]\,\subset
\h {A'}_{\le 0}.$$
Therefore (by \ref{DML})
$${\rm Dk}\,(A') \subset A'_0\neq A'\,,$$
so that
$A'\not\simeq \C^{[3]}$ and 
$X\not\simeq \C^{3}$,
which proves \ref{exo1}. \elem

\no\proof The Jacobian derivation
$$\p_0:={\p (\h p , \,x,\,y,\,*)\over 
\p (x,\,y,\,z,\,u)}$$
on the algebra $\h {A'}$ 
being homogeneous and locally 
nilpotent with ${\h {A'}}^{\p_0}=
{\rm ker}\,\p_0=\C[\h x,\,\h y ]$, 
it is sufficient to show that
$\h x,\,\h y\in {\rm ker}\,{\h \p}$
for any
$\h \p\in {\rm LND}_{\rm gr} (\h {A'})$. 
Indeed, by \ref{trdeg} both
$\C[\h x,\,\h y ]$ and  ${\rm ker}\,{\h \p}$
are algebraically closed subalgebras 
of $\h {A'}$ of
transcendence degree 2, thereby
they coincide
provided that 
$\C[\h x,\,\h y ]\subset {\rm ker}\,{\h \p}$.

The derivation $\h \p$ being homogeneous, 
the subalgebra
$ {\rm ker}\,{\h \p}\subset {\h A}$ is 
generated by homogeneous elements
(i.e., by the restrictions to $\h X$ of 
$d$-homogeneous polynomials on 
$\C^{4}$). 
Let $a=q|_{\h X}\in {\rm ker}\,{\h \p}$ 
(with a $d$-homogeneous $q\in\C^{[4]}$)
be nonconstant.
We may assume that deg$_z q<e$ 
(otherwise we replace
the polynomial $q$ by the rest of the 
Euclidean division of $q$ by
the $z$-monic $d$-homogeneous polynomial
$\h p = z^{e}+(x^k-y^l)^mu$).

The kernel ${\h A}^{\h \p}={\rm ker}\,{\h \p}$ 
being factorially closed
(see \ref{trdeg}(c)), the irreducible factors of 
$q$ restricted to $\h X$ belong to this kernel 
as well.
Therefore, ${\rm ker}\,{\h \p}$ 
is generated by the
traces of irreducible $d$-homogeneous 
polynomials
$q$ with  deg$_z q<e$.
By \ref{trdeg}(d) we have ${\h x},\,{\h y}\in 
{\rm ker}\,{\h \p}$
provided that
$\lambda {\h x}^k+\mu {\h y}^l\in 
{\rm ker}\,{\h \p}$
for some $\lambda,\,\mu \in\C^*$. 
Due to \ref{l3} and 
to the above argument on algebraic closeness, 
${\rm ker}\,{\h \p}$
coincides with the subalgebra  of $\h {A'}$ 
generated by
one of the following pairs:
$$({\h x},\, {\h y}),\quad ({\h x},\, {\h u}),
\quad ({\h y},\, {\h u}),\quad ({\h x},\,
{\h z}),\quad ({\h y},\,{\h z}),\quad 
({\h z},\,{\h u})\,.$$
If $\h z \in {\rm ker}\,{\h \p}$ then
by \ref{trdeg} and the equality
\query{(label=e2)}\be\label{e2} 
{\h z}^{e}+({\h x}^k+{\h y}^l)^m{\h u}=0\,
\ee 
also ${\h x},\,{\h y},\,{\h u}
\in {\rm ker}\,{\h \p}$, whence $\h \p=0$, 
a contradiction.
This eliminates the last three cases.

If ${\rm ker}\,{\h \p}=\C[{\h x},\,{\h u}]$ 
\footnote{
The case where ${\rm ker}\,{\h \p}=
\C[{\h y},\,{\h u}]$
can be eliminated by a similar argument.} 
then (\ref{e2}) yields
the equality
\query{(label=e3)}\be\la{e3}
lm \,{\rm deg}_{\h \p}\, {\h y} =
e \,{\rm deg}_{\h \p}\, {\h z}\,.\ee
On the other hand, as ${\rm ker}\,{\h \p}=
\C[{\h x},\,{\h u}]$, by \ref{KML1}
$\h \p$ is equivalent to 
the Jacobian derivation
$${\h \p}_1:= {\p({\h p},\,x,\,u,\,*)\over
\p(x,\,y,\,z,\,u)}\,.$$
It is easily seen that $${\h \p}_1\,(\h y) =
- ({\h p})'_z|_{\h X}
= - {\h z}^{e-1}\,,$$
and so as $\deg_{\p}=\deg_{\p_1}$
 (see \ref{JD}) we get
\query{(label=e4)}\be\label{e4} 
{\rm deg}_{\h \p}\, {\h y} =
(e-1){\rm deg}_{\h \p}\, {\h z}+1\,.\ee
From (\ref{e3}) and (\ref{e4}) we obtain:
$$\big\lbrack{(1-{1\over lm})e -1}\big\rbrack
{\rm deg}_{\h \p}\, {\h z} =-1\,.$$
But this is impossible because 
$l,\,e \ge 2,\,m\ge 1$, 
and so
$(1-{1\over lm})e\ge 1$.
This completes the proof of the lemma, 
as well as those of \ref{exo1}. 
\qed\medskip

Next we consider
the remaining possibility in \ref{tra}.

\query{(label=exo2)}\no\bprop\la{exo2}
Let $X$ and $Y$ be as in \ref{conv3}. 
Suppose that $\GNH$ is a smooth curve, 
whereas
there exists a singular component, say,
$C_{h+1}$ of the curve $\CS$.
Then once again, the Derksen invariant
${\rm Dk}\,(A')$  
is non-trivial, 
whence $X\not\simeq \C^3$.\eprop 

The proof is done in \ref{l5}-\ref{l7} 
below.

\query{(label=l5)}\no\blem\la{l5} 
Under the assumptions 
as in \ref{exo2}, after applying 
an appropriate 
tame automorphism of 
$\C^3_{x,y,z}$ the polynomial $p$ 
can be presented in the form
\query{(label=pres3)}\be\la{pres3} 
p=x^mh_0(x,y,z)u+y^k+z^l+
xh_1(x,y,z)\ee
with 
$k,\,l,\,m\ge 2,\,\,(k,\,l)=
1,\,\,h_0,\,h_1\in\C[x,y,z]$ 
and $h_0(0,y,z)={\rm const}=1$. \elem

\no\proof The curve $\G_{h+1}\subset\C^2_{x,y}$ 
(being smooth and
homeomorphic to $\C$) is isomorphic to $\C$.
Thus by the Abhyankar-Moh-Suzuki Theorem, 
choosing appropriate 
new coordinates
in the $(x,y)$-plane we may suppose 
that $f_{h+1}=x$. We factorize
$f=f_{h+1}^m {\tilde f} =x^m{\tilde f}$ with 
$m:=a_{h+1}$. As $C_{h+1}$ is
singular, by \ref{smooth} the gradient 
grad$f$ has to vanish 
at some point
of $\G_{h+1}=\{x=0\}$. The component 
$\G_{h+1}$ of $\G$
being isolated 
(see \ref{3f}($\beta$)) this implies 
that $m\ge 2$. 
Furthermore the polynomial
${\tilde f} (0,\,y)$ does not vanish,  
whence is a non-zero constant;
by rescaling the $x$-coordinate we may 
assume this constant being 1.

As the curve $C_{h+1}=\{x=0=
g(0,\,y,\,z)\}\subseteq\C^2_{y,z}$
is homeomorphic to $\C$ and singular, 
by the 
Lin-Zaidenberg Theorem (after
performing an appropriate 
automorphism of the plane
$\C^2_{y,z}$) we may suppose that 
$C_{h+1}$ is given by
$x=y^k+z^l=0$ with $k,\,l\ge 2$ and $(k,\,l)=1$.
In virtue of \ref{Ci}
and the condition \ref{3f}($\delta$) 
the divisors $D_{h+1}=\{x=0\}\subset 
\C^3_{x,y,z}$ and $D_g^{\rm red}$ meet
transversally at general points of 
the curve $C_{h+1}$, whence $g(0,\,y,\,z)=y^k+z^l$ 
(up to a constant factor which can be
put equal to $1$) i.e.,
$g=y^k+z^l+xh_1(x,y,z)$. 
Taking for $h_0(x,y,z)$ 
the polynomial obtained
from ${\tilde f}$ as a result of 
the latter coordinate change,
we obtain the desired presentation 
(\ref{pres3}) with $h_0(0,y,z)=
{\tilde f}(0,y)=1$. 
\qed\ms

\query{(label=DF2)}\no\bsit\la{DF2}{\rm
We consider the
weight degree function  $d$
on the algebra $A':=\C[X]$ with 
$$d_x=-1,\,\,d_y=d_z=0
\quad {\rm and}\quad d_u=m\,.$$
As the $d$-principal part $\h p=x^mu+y^k+z^l$ 
of the polynomial $p$
as in (\ref{pres3}) is irreducible, 
by \ref{KML2}
$\h {A'}:=\C[\h X]$ is just the graded 
algebra associated
with the filtration on $A'$ defined 
by the degree function $d$.
Notice that in $\h {A'}$ the 
following relation holds:
\be\label{rel} \h x^m\h u+
\h y^k+\h z^l=0\,.\ee}
\esit

With this notation we have 
the following lemma.

\query{(label=l6)}\no\blem\label{l6} 
${\rm Dk_{gr}}\,(\h {A'}) \subset 
\h {A'}_{\le 0}\,,$
and so (by \ref{DML})
${\rm Dk}\,(A') \subset A'_0\neq A'\,.$
Thereby
$A'\not\simeq \C^{[3]}$ and $X\not\simeq \C^{3}$,
which yields \ref{exo2}.\elem

\no\proof We must show that 
${\rm ker}\,{\h \p}\subset \h {A'}_{\le 0}$
whenever ${\h \p}\in {\rm LND_{gr}}(\h {A'})$. 
Assume that there
exist
${\h \p}\in {\rm LND_{gr}}(\h {A'})$ and 
$\h a \in {\rm ker}\,{\h \p}$ 
such that
$\h a \notin \h {A'}_{\le 0}$. Furthermore, 
by \ref{trdeg}(c)
we may suppose that this element $\h a$
is non-decomposable. 
We have $\h a = q|_{\h X}$ for some 
$d$-homogeneous polynomial
$q=\sum_{i,j} a_{ij}x^iu^j\in\C^{[4]}$ 
(with $a_{ij}\in\C[y,\,z]$).
Moreover, taking into account (\ref{rel}) 
we may suppose that
$i<m$ whenever $j>0$. \ms

\no {\bf Claim.} {\it In the above expression 
for $q$ there is 
only one non-zero monomial.}\ms

\no {\it Proof of the claim.}
Indeed, if there were two of them, say, if
$a_{i_1j_1}\neq 0$ and $a_{i_2j_2}\neq 0$ then
we would have
\be\label{re1}
d(q)=mj_1-i_1=mj_2-i_2\qquad \Longrightarrow
\qquad m(j_1-j_2)=i_1-i_2\,.\ee
Assuming that $i_1>i_2$, by (\ref{re1})
we get $j_1>j_2$, and vice versa.
Thus $j_1>0$, and then by our assumption
$i_1<m$. Hence also $i_1-i_2< m$, 
which contradicts
(\ref{re1}). This proves the claim. \qed\ms

Since the element
$\h a= q(\h x,\,\h y,\,\h z,\,\h u)=
a_{ij}(\h y,\,\h z) \h x^i \h u^j$
is supposed to be non-decomposable 
and of positive $d$-degree,
we have $\h a=\h u\in {\rm ker}\,{\h \p}$.
Thus $\h \p$ can be specified
to a locally nilpotent derivation
$\p_1$ of the algebra $B=\C[S]$, where
$S:=\{x^m+y^k+z^l=0\}=\{u=1\}\cap 
\h X \subset\C^3$ 
(see \ref{SPE}).

By \ref{trdeg}(a), 
tr.deg$\,({\rm ker}\,{\h \p})=2$, 
whence
there is a homogeneous $ \h \p$-constant $\h b$ 
such that
the elements $\h a=\h u,\,\h b\in \h {A'}$ 
are algebraically independent.
As above we obtain that either
$\h b=b(\h y,\,\h z)$
for some irreducible polynomial 
$b\in\C^{[2]}\setminus \C$,
or $\h b=\h x$. In the latter case by 
(\ref{rel}) we  have
$\h y^k+\h z^l\in {\rm ker}\,{\h \p}$, 
and thus by \ref{trdeg}(d)
also $\h y,\,\h z \in {\rm ker}\,{\h \p}$. 
Therefore 
$\h \p=0$, which is impossible. 
Finally, we conclude that
$b(\h y,\,\h z)\in {\rm ker}\,{\h \p}$ 
for a certain polynomial 
$b\in\C^{[2]}\setminus \C$, 
and so the restriction $b|_S$
is a $\p_1$-constant.

Now the proof can be completed 
by applying the next lemma
(cf. \cite{De, KaZa2, Za2}). \qed\ms

\query{(label=l7)}\no\blem\label{l7}
Let $B=\C[S]$ where
$S:=\{x^m+y^k+z^l=0\}\subset\C^3$ and 
$k,\,l,\,m\ge 2,\,\,\,
{\rm gcd} (k,l)=1$, and let
$\p_1\in {\rm LND}\,(B)$.
Then $b|_S\not\in {\rm ker}\, \p_1$ 
whenever 
$b\in\C[y,\,z]\setminus \C$. \elem

\no\proof We define the weight degree 
function $\h d$ on the algebra
$B$ by letting $$\h d_x=kl,\quad \h d_y=
lm,\quad \h d_z=km\,.$$
Actually $x^m+y^k+z^l$ is a 
$\h d$-homogeneous polynomial,
$B$ is a graded algebra, and we may 
consider
the associated homogeneous derivation 
${\h \p}_1\in 
{\rm LND_{gr}}\,(B)$.
Assuming that for a polynomial 
$b\in\C[y,\,z]\setminus \C$,
$b|_S\in {\rm ker}\,  \p_1$
we will get a $\h d$-homogeneous 
polynomial
$\h b_1\in\C[y,\,z]\setminus \C$ 
such that  
$\h b_1|_S\in {\rm ker}\, {\h {\p}}_1$ 
as well.
By \ref{trdeg}(c) an irreducible
$\h d$-homogeneous factor of the polynomial 
$\h b_1$
(this can be $y,\,z$ or $\lambda y^k+\mu z^l$, 
where 
$\lambda, \mu \in\C^*$)
restricts to $S$ as a ${\h \p}_1$-constant.
If it were $\lambda y^k+\mu z^l$ then by 
\ref{trdeg}(d)
both $y|_S$ and $z|_S$ would be 
${\h \p}_1$-constants,
whence ${\h \p}_1=0$, a contradiction.

Thus we may assume that, say, 
$y|_S \in {\rm ker}\,{\h \p}_1$.
As $(x^m+y^k+z^l)|_S=0$
we have 
$(x^m+z^l)|_S \in {\rm ker}\,{\h \p}_1$,
and (since $m,\,l\ge 2$) 
again by \ref{trdeg}(d)
we  obtain that $x|_S,\,z|_S \in 
{\rm ker}\,{\h \p}_1$,
which is impossible. 

This completes
the proof of \ref{l7}, \ref{l6} and \ref{exo2}. 
\qed\ms

\query{(label=MLN)}\no\brem\la{MLN} 
{\rm Choosing 
a weight degree function likewise in \ref{ONT}
and repeating the proof of Proposition 5.3 
in \cite{Ka2}  
it can be shown that under assumptions 
as 
in \ref{tra} (that is, as in \ref{exo1} or 
\ref{exo2}) one has
$\C[x]\subseteq {\rm ML}(A')$, 
whence the Makar-Limanov invariant 
${\rm ML}(A')$
is non-trivial as well.  
}\erem

\query{(label=ssec-c3)}
\subsection{Simple affine modifications 
of $\C^{3}$
isomorphic to  $\C^{3}$}\la{ssec-c3}
The following theorem is a central 
result of the paper. 
The expression 
`in appropriate $(x,y)$-coordinates' 
below 
means `after performing an appropriate 
automorphism 
$\alpha\in
\Aut\,\C[x,y]\subset\Aut\,\C[x,y,z,u]$'. 
(In fact,
the new $(x,y)$-coordinates are chosen 
in a way to get $\fnh\in\C[x]$; in particular, 
no coordinate change is necessary when 
$\fnh=$const.) 

\query{(label=mthsec5)}\no\bthm\la{mthsec5} 
For a hypersurface 
$X \subset \C^4$ given
by an equation
$$
  p=f(x,y)u+g(x,y,z)=0\qquad 
(f\in\C^{[2]}\backslash 
\{0\}\,,g\in\C^{[3]})\,,
$$
in appropriately chosen new 
$(x,y)$-coordinates
the following conditions {\rm (i)-(vi)} are equivalent.

\begin{enumerate}
\item[{\rm (i)}] $X\simeq\C^3$.
\item[{\rm (ii)}] Every fiber
$X_{\mu}=p^{-1}(\mu)\,\,\,(\mu\in \C)$ 
of the polynomial $p\in\C^{[4]}$
is isomorphic to $\C^3$.
\item[{\rm (iii)}]  
Every fiber of the regular function 
$x|_{X}\in A':=\C[X]$ is reduced and
isomorphic to $\C^2$.
\item[{\rm (iv)}] 
Every fiber of the morphism 
$\rho: \C^4\to\C^2,\,\,\,
(x,y,z,u)\longmapsto (x,\,p(x,y,z,u))$, 
is reduced and
isomorphic to $\C^2$.
\item[{\rm (v)}] The polynomial
$p$ is a residual $x$-variable 
\footnote{That is (\ref{spe})
for any $\lambda\in\C$, the specialization 
$p_{\lambda}\in\C[y,z,u]$ of $p$ 
is a variable.}.
\item[{\rm (vi)}] With $Y:=\C^3$ 
and $X\subseteq \C^4$ as above being 
smooth and irreducible,
the conditions \ref{3f}($\alpha$)-($\delta$) 
hold 
as well as the following one:
\item[($\varepsilon$)] The curves 
$\CNH$ and $\GNH$ are smooth.
\end{enumerate}
\ethm

\no\proof The implications 
(v)$\Rightarrow$(iv)$\Rightarrow$(iii)
 and (ii)$\Rightarrow$(i)
 are immediate; (i)$\Rightarrow$(vi)
follows from \ref{c3g} and \ref{tra}. 
Thus it suffices to
establish (vi)$\Rightarrow$(v) and 
(iii)$\Rightarrow$(ii).

Let us show (vi)$\Rightarrow$(v). If 
$\fnh\neq {\rm const}$ (that is, $h<n$) 
then
in virtue of our assumptions 
\ref{3f}($\gamma$) 
and 
(vi)($\varepsilon$), $\G_{h+1}\simeq\C$. 
By the 
Abhyankar-Moh-Suzuki Theorem, 
after performing an 
appropriate automorphism of the plane 
$\C^2_{x,y}$ 
we may suppose that $f_{h+1}=x$. Hence 
by \ref{3f}($\gamma$)
we have
$$\fnh=c\prod_{i=h+1}^n 
(x-\lambda_i)^{a_i}\in\C[x]\,.$$ 
Anyhow, $\fnh= {\rm const}$ or not, 
assuming that $\fh\neq {\rm const}$ 
by \ref{vhs} we may write
$$p_{\lambda}=\kappa f_{{\rm horiz},\lambda} u 
+ 
b_{0, \lambda} + b_{1, \lambda} z + 
f_{{\rm horiz},\lambda}^{\rm  red} 
h_{\lambda}$$ 
with 
$ f_{{\rm horiz},\lambda}, 
\,f_{{\rm horiz},\lambda}^{\rm  red},\,
b_{0, \lambda},\,  b_{1, \lambda}\in 
B:=\C[y],\,\,h_{\lambda}\in B[z]$
and $\kappa:=\fnh({\lambda})=c \prod_{i=h+1}^n 
(\lambda-\lambda_i)^{a_i}\in\C$. 
If $\lambda\in\C\backslash \{\lambda_{h+1},
\ldots,\lambda_n\}$ 
(i.e., $\kappa\neq 0$) then
in virtue of \ref{vhs}, $b_{1, \lambda}\in B$ 
is invertible 
$\mod\,f_{{\rm horiz},\lambda}$, whereas 
$f_{{\rm horiz},\lambda}^{\rm  red}$ is 
nilpotent $\mod\,f_{{\rm horiz},\lambda}$. 
Therefore by \ref{RuVe}, 
$p_{\lambda}\in B[z,u]$ is a $B$-variable; 
in particular, $p_{\lambda}\in \C[y,z,u]$ 
is a variable (whence for every $\mu\in\C$, 
the surface 
$\{p_{\lambda}=\mu\}\subseteq\C^3_{y,z,u}$
is reduced and isomorphic to $\C^2$). 

For every $i=h+1,\ldots,n$ we have 
$f_i=x-\lambda_i$ and
in virtue of \ref{3f}($\gamma$) and 
(vi)($\varepsilon$),
$\G_i\simeq\C\simeq C_i$. Moreover by 
\ref{smooth}(i)
the polynomial $g_{\lambda_i}\in \C[y,z]$ 
is irreducible
and $g_{\lambda_i}^{-1}(0)=C_i\simeq\C$, 
whence 
by the Abhyankar-Moh-Suzuki Theorem 
it is a variable
of the polynomial ring $\C[y,z]$. 
Thus $p_{\lambda_i}=g_{\lambda_i}\in\C[y,z,u]$ 
is a variable as well. 
Now (v) follows. 
In the remaining case where $\fh= {\rm const}$
the proof is similar (but simpler) 
and left to the reader. 

Next we prove (iii)$\Rightarrow$(ii).
We notice first of all that 
under the condition (iii), the hypersurface
$X$ is smooth and irreducible
(indeed, for every $\lambda\in\C$ we have
$${\rm grad}\, p_{\lambda}=
{\rm pr}\,(({\rm grad}\, p)_{\lambda})\qquad
\mbox{with}\qquad {\rm pr}: (x,y,z,u)
\longmapsto (y,z,u)\,{\rm ).}$$
Similarly, under the condition (iv) 
every fiber $X_{\mu}\,\,\,(\mu\in\C)$ 
is smooth and irreducible. 

Actually we establish 
(iii)$\Rightarrow$(i), 
which at the same time shows 
(iv)$\Rightarrow$(ii). 
Together with the implication 
(i)$\Rightarrow$(iv)
(which we already know) these yield 
(iii)$\Rightarrow$(ii),
as needed. 

To prove the implication 
(iii)$\Rightarrow$(i),
in virtue of \ref{miy}
it is enough to show that a 
smooth, irreducible $3$-fold 
$X$ satisfying (iii) is acyclic. 
Let $U=\C\backslash V$ be a Zariski 
open subset 
such that over $U$, the morphism 
$x|_X : X\to\C$
is a smooth fibration 
(with the fibers diffeomorphic to $\R^4$).
In the notation of \ref{prhom} we let
$\hX=X,\,\hY=Y=\C^3,\,\hE=
(x|_X)^{-1}(V)\subseteq \hX$
and $\hD=x^{-1}(V)\subseteq \hY$. 
It is easily seen that the induced 
homomorphisms
$$(x|_{\hX\backslash \hE})_*: 
H_*(\hX\backslash \hE)\to 
H_*(\C\backslash V)\,,$$
$$(x|_{\hY\backslash \hD})_*: 
H_*(\hY\backslash \hD)\to 
H_*(\C\backslash V)\,,$$ 
and hence also 
$$(\si|_{\hX\backslash \hE})_*: 
H_*(\hX\backslash \hE)\to 
H_*(\hY\backslash \hD)\,$$ are isomorphisms, 
as well as
$$(\si|_{\hE})_*: H_*(\hE)\to H_*(\hD)\,$$ 
(indeed, $H_0(\hE)\cong \Z^{{\rm card}\,V}
\cong H_0(\hD)$ and the senior homology are
trivial). Thus the conditions (i) and (ii) of 
\ref{prhom} are fulfilled, and so by \ref{prhom} 
$$\si_* : H_*(\hX)\to H_*(\hY)$$ 
is an isomorphism. 
Therefore, the $3$-fold $X=\hX$ is acyclic. 
This completes the proof.
\qed\medskip

\query{(label=FR)}\brem\la{FR} 
{\rm For another proof 
of the implication
(iii)$\Rightarrow$(i) which do not use \ref{miy}
see \ref{x-planes}, \ref{anpr} below. 
}\erem

\query{(label=NTM)} \no\bsit\la{NTM} {\rm
If the polynomial
$p\in\C^{[4]}$ is linear with respect to two 
(and not just one) variables 
we have a much stronger
result (see \ref{mthsec6} below). 
It is an analog in dimension $4$ 
of Theorem 7.2 in \cite{KaZa1} generalizing
Sathaye's Theorem \cite{Sat1}. For a polynomial
$\varphi\in\C[x,y]$ we let 
$\G_{\varphi}:=V(\varphi)\subseteq\C^2$. 
Also, we use below the variable $v$ instead of $z$
to emphasize the symmetry of the situation.
}\esit

\query{(label=mthsec6)}\no\bthm\la{mthsec6} 
For a hypersurface 
$X \subset \C^4$ given
by an equation
\query{(label=li2)}\be\la{li2}
  p=\tilde a(x,y)u+\tilde b(x,y)v+c(x,y)=
d(au+bv)+c=0
\ee
with $a,b,c,d\in\C[x,y],\,\,\, d\neq 0,\,\,(a,b)\neq (0,0)$ 
and ${\rm gcd} (a,b)=1$,
in appropriately chosen new 
$(x,y)$-coordinates
the following conditions {\rm (i$'$)-(v$'$)} 
are equivalent.

\begin{enumerate}
\item[{\rm (i$'$)}] 
The $3$-fold $X$ is irreducible, 
smooth and acyclic.
\item[{\rm (ii$'$)}] $X\simeq\C^3$.
\item[{\rm (iii$'$)}] The polynomial
$p\in\C^{[4]}$ is a residual $x$-variable.
\item[{\rm (iv$'$)}] $p\in\C[x]^{[3]}$ 
is an $x$-variable.
\item[{\rm (v$'$)}] $d\in\C[x]$, 
{\rm gcd}$(c,d)=1$, $\G_a\cap\G_b\subseteq\G_d$,
and for every root $x=x_i$ of 
$d\,\,\, (i=1,\ldots,\deg d)$ 
we have
\query{(label=fpl)}\be\la{fpl}
c_{x_i}(y)=c(x_i,y)\in\C^*y+\C\,.\ee  
\end{enumerate}
\ethm

\no\proof The strategy of the proof is as follows. 
In virtue of \ref{mthsec5} the conditions 
(ii$'$) and (iii$'$) are equivalent 
to each other and to the other conditions of 
\ref{mthsec5}. The equivalence 
(iii$'$)$\Longleftrightarrow$(iv$'$) will 
be established later on in \ref{partialres}(b).
By \ref{3f} we have 
(i$'$)$\Longrightarrow
$\ref{3f}($\alpha$)-($\delta$). 
We show below that (i$'$) also
implies the condition ($\varepsilon$) of
\ref{mthsec5}(vi). In virtue of \ref{mthsec5}
this yields the implications
(i$'$)$\Longrightarrow$\ref{mthsec5}(vi)
$\Longleftrightarrow$(ii$'$)$\Longrightarrow$(i$'$),
and so gives (i$'$)$\Longleftrightarrow$(ii$'$). 
Finally we show
(iii$'$)$\Longleftrightarrow$(v$'$), 
which concludes the proof. 

Interchanging (if necessary) the roles of $u$ and $v$ 
we may suppose in the sequel that
$a\neq 0$.

(i$'$)$\Longrightarrow$($\varepsilon$). 
We observe (see \ref{fct}) that for $p$ as 
in (\ref{li2}), $\CS=\GS=\emptyset$, 
whence $\CNH=\CV$ is smooth. 
And it was shown in the proof of \ref{lemGNHsing}
that if $\GNH$ is singular then $\GNH=\GS$. 
Henceforth
in our setting $\GNH$ is also smooth, that
shows ($\varepsilon$).

(iii$'$)$\Longrightarrow$(v$'$).
As $X$ is irreducible we have gcd$(c,d)=1$.  
Letting
$q:=p-c$ we
fix a point $P=(x_0,y_0)$ such that
$$q_P=d(x_0,y_0)(a(x_0,y_0)u+b(x_0,y_0)v)=0\,.$$
It is easily seen that $y-y_0$ divides 
$p_{x_0}-c(x_0,y_0)=q_{x_0}+(c_{x_0}-c(x_0,y_0))$. 
As $p_{x_0}\in\Var_{\C}[y,u,v]$
we obtain that $q_{x_0}=0$ and
$p_{x_0}-c(x_0,y_0)=c_{x_0}-c(x_0,y_0)=\kappa (y-y_0)$
with $\kappa\in\C^*$, as needed in (\ref{fpl}). 
As gcd$(a,b)=1$, from $q_{x_0}=0$ we get 
$d_{x_0}=0$. Thus for every root 
$(x_0,y_0)$ of $d=0$ we have $d_{x_0}=0$, which 
means that $d\in\C[x]$. Moreover, for every root 
$(x_0,y_0)$ of $a=b=0$ we also have  
$d_{x_0}=0$, which shows the inclusion 
$\G_a\cap\G_b\subseteq\G_d$ and gives (v$'$).

(v$'$)$\Longrightarrow$(iii$'$).
In view of (\ref{fpl}) 
for every root $x_i$ of $d$ 
($i=h+1,\ldots,n$),  
$p_{x_i}=c_{x_i}\in \Var_{\C} \C[y]\subseteq
\Var_{\C} \C[y,u,v]$. Let further
$\lambda\in\C$ and $d(\lambda)\neq 0$. 
As by (v$'$), $\G_a\cap\G_b\subseteq\G_d$
we have gcd$(a_{\lambda}, b_{\lambda})=1$, 
whence 
$r_{\lambda}a_{\lambda}-s_{\lambda}b_{\lambda}=1$
for certain polynomials 
$r_{\lambda},s_{\lambda}\in\C[y]$. 
As the matrix 
\begin{displaymath}
\mbox{$A_{\lambda}=$}
\left( \begin{array}{cc}
a_{\lambda} & b_{\lambda} \\
s_{\lambda} & r_{\lambda} \\
\end{array} \right)
\end{displaymath}
is invertible over the ring $\C[y]$, the mapping
$$\alpha : (u,v)\longmapsto (u_1,v_1):
=A_{\lambda}(u,v)$$
induces a $\C[y]$-automorphism of $\C[y][u,v]$
with $\alpha(p_{\lambda})=d(\lambda)u_1+c_{\lambda}
\in\Var_{\C} \C[y,u,v]$. Thus (iii$'$) follows.
\qed


\query{(label=VAR)}\section{Simple birational 
extensions of $\C^{[3]}$ as variables
in $\C^{[4]}$}\la{VAR}

\query{(label=ssecPART)}
\subsection{Partial positive results}\la{ssecPART}

We recall the problem stated in the Introduction. 

\query{(label=pro)} \no\bsit\la{pro} {\bf Problem.} 
{\it Is it true that 
if a hypersurface $X\subseteq \C^4$ 
with equation of the form
$$p:=f(x,\,y)u+g(x,\,y,\,z)=0\qquad 
({\rm where}\quad f\in \C^{[2]}\backslash \{0\}\quad 
{\rm and}\quad g\in \C^{[3]})$$
is isomorphic to $\C^3$ then necessarily
$p\in\Var_{\C} \C^{[4]}$ and moreover,
$p\in\Var_{B} B^{[3]}$ (that is, 
$p$ is an $x$-variable)
provided in addition that 
$\fnh\in B:=\C[x]$?
}\esit

The principal results of this subsection 
can be summarized as follows.

\query{(label=partialres)} \no\bthm 
\la{partialres} 
The answer to \ref{pro} is positive
in each of the following cases:

\begin{enumerate}
\item[{\rm (a)}] $f \in\C[x]$.
\item[{\rm (b)}] $\deg_z g\le 1$.
\item[{\rm (c)}] $f_{\rm vert}=1$.
\item[{\rm (d)}] $\fh=f_{\rm horiz}^{\rm red}$.
\end{enumerate}
\ethm

\no The proof is given in \ref{prtr} below.
From \ref{partialres}(c),(d) 
we obtain the following corollary.

\query{(cova)} \no\bcor\la{cova} If in 
\ref{pro}, $f=f_1^{a_1}$ 
with $f_1\in\C[x,y]$ irreducible then $p$ 
is a variable. \ecor

\query{(reva)} \no\bsit\la{reva} {\rm
By \ref{mthsec5}(v) we may suppose 
that $p$ as in \ref{pro} is a residual $x$-variable 
(see \ref{var}, \ref{spe}).
Thus \ref{pro} would be answered in positive 
if the following conjecture 
\footnote{Proposed by the second author.} 
were true.
}\esit

\query{(label=xvar)} \no\bconj\la{xvar} 
{\it If  
$p\in \C^{[n]}$ is a residual $x_1$-variable 
(that is,  
$$p_{\lambda}:=p(\lambda,x_2,\dots,x_n)
\in \Var_{\C} 
\C^{[n-1]}\qquad \forall\lambda\in \C)$$
then $p$ is an $x_1$-variable: 
$p\in\Var_{\C[x_1]} \C[x_1]^{[n-1]}$. 
}\econj

In \ref{f=0}-\ref{form} below 
we analyze the situation 
with the only assumption that $p=fu+g$
is a residual $x$-variable. 
First in \ref{f=0} we include 
the case where $f=0$,
then in \ref{specase} we deal with 
the special case 
where $f\in\C[x]\backslash \{0\}$, 
whereas the general case 
is treated in \ref{form}.

\query{(label=f=0)}\no\bprop\label{f=0}
Let $X$ be a smooth affine surface and 
$\varphi : X \to \C$ be a morphism with 
$\varphi^*(\lambda)\simeq \C\,\,\,
\forall \lambda\in\C$.
Then the following hold. 

\begin{enumerate}
\item[(a)] $X\simeq\C^2$. 
\item[(b)] Furthermore, 
if $X = {\rm spec}\,\bigl (\C^{[3]}/(g)\bigr )$
with $g\in \C^{[3]}=\C[x,y,z]$ 
and $\varphi=x|_X$ then 
$g$ is an $x$-variable. 
\item[(c)] Consequently, 
a residual $x$-variable $g\in\C[x,y,z]$ is
an $x$-variable. \footnote{Moreover, $g\in\C[x,y,z]$
is a residual $x$-variable of the ring $\C[x][y,z,u]$
iff it is an $x$-variable. Indeed, 
this follows from (c) and the cancellation for curves:
$\G\times\C\simeq\C^2\Longrightarrow\G\simeq\C$.}
\end{enumerate}
\eprop

\no\proof (a) The fact is well-known 
(e.g., cf. \cite[Thm. 2.2.1]{Miy3}); 
nevertheless we indicate the proof. 
First extending 
$\varphi : X \to \C$ to a 
projective morphism $\bar \varphi : V \to \PP^1$
of a smooth ruled surface $V$, 
we then pass to a smooth relatively minimal model 
of $(V, \bar \varphi)$ by
contracting successively the superfluous components
of the reducible fibers of $\bar \varphi$
different form the closures of the original fibers
of $\varphi$
(this is always possible, see e.g. 
\cite[Lemma 7]{Gi}, \cite[Lemma 3.5]{Za1} or 
\cite[Lemma 1.4.1(6)]{Miy3}). 
Thus we obtain (as a completion
 of the original family) 
a Hirzebruch surface $\pi : \Sigma_n\to\PP^1$ 
with a section $s: \PP^1\to \Sigma_n$ 
`at infinity',
so that $X=\Sigma_n\backslash 
(s(\PP^1)\cup F_{\infty})$. 
By means of elementary transformations over 
the point $\infty\in\PP^1$ we replace 
$\Sigma_n$ with $\Sigma_0=\PP^1\times\PP^1$ 
and $s(\PP^1)$ with a constant section,
say, $C_0$. 
This yields the desired isomorphism 
$X=\Sigma_0\backslash 
(C_0\cup F_{\infty})\simeq\C^2$.

(b) By (a) there is an isomorphism 
$$\gamma : \C[x,y,z]/(g)\to \C^{[2]}=\C[x,y]\,.$$
Letting $h:=\gamma (x)\in\C^{[2]}$, 
by our assumption we obtain:
$$\C^{[2]}/(h)=\C[x,y,z]/(g,x)=
\C[\varphi^*(0)]\cong \C^{[1]}\,.$$
Hence by the Abhyankar-Moh-Suzuki Theorem, 
up to an automorphism of $\C[x,y]$
we may suppose that $h=x$ i.e.,
$$\gamma : \C[x][y,z]/(g)\to \C[x]^{[1]}\,.$$ 
Now it follows from the generalized version 
of the Abhyankar-Moh-Suzuki Theorem
\cite[Thm. 2.6.2]{RuSat} that 
$g$ is an $x$-variable, as stated.

(c) immediately follows from (a) and (b). 
\qed\medskip

\query{(label=specase)}\no\bprop\label{specase}
If $p=f(x)u+g(x,y,z)\in\C^{[4]}$ (with 
$f\in\C[x]\backslash \{0\}$)
is a residual $x$-variable then 
actually it is an
$x$-variable.
\eprop

\no\proof As $u+g(x,y,z)$ is an $x$-variable, 
by \ref{corRuVe}(a) so is $p=fu+g$.\qed\medskip

\query{(label=form)}\no\bprop\label{form} 
If $p=f(x,y)u+g(x,y,z)\in \C[x][y,z,u]$ 
with $f\notin \C[x]$
is a residual $x$-variable then $p$ 
can be written as follows:
\query{(label=pform)}\be\la{pform}
 p= q[r\tilde f u+\tilde g_0 y^2
+\tilde g_1 z+\tilde f^{\rm red} \tilde g_2 z^2]
+a_0+a_1 y\,,
\ee
where 
\begin{enumerate}
\item[$\bullet$] 
$\quad q, r, a_0, a_1\in\C[x],\,\,\,\tilde f, 
\tilde g_0,
\tilde g_1\in\C[x,y],\,\,\,\tilde g_2\in\C[x,y,z]$,
\item[$\bullet$] $\quad \tilde f^{-1}(0)\cap 
\tilde g_1^{-1}(0)$ 
is a finite subset of $q^{-1}(0)\times \C_{y}$, and 
\item[$\bullet$] $\quad$for every root $x_0$ of $r$,  
$\,\,\,q(x_0)\tilde f(x_0,y)\in \C$. 
\end{enumerate}
\eprop

The proof is done in \ref{prpr} below. 

\query{(label=cform)} \bcor\la{cform}
If $p=f(x,y)u+g(x,y,z)\in \C[x][y,z,u]$ 
is a residual $x$-variable then it is a
$\C(x)$-variable of $\C(x)^{[3]}$.
\ecor

\no\proof As a $\C[x]$-variable is also a 
$\C(x)$-variable, by \ref{f=0} 
and \ref{specase}
we may suppose that $f\notin\C[x]$, and so 
\ref{form} applies. 
Let $p$ be presented as in (\ref{pform}).
As $\tilde f^{-1}(0)\cap \tilde g_1^{-1}(0)$ 
is finite, the polynomials 
$qr\tilde f,\, q\tilde g_1\in\C[x,y]$ regarded 
as elements of $B:=\C(x)[y]$ are coprime. 
Moreover the polynomial
$\tilde f^{\rm red} \tilde g_2 z^2\in B[z]$ 
is nilpotent $\mod (qr\tilde f)$, and so by
\ref{RuVe}(a), $p=fu+g\in\Var_B B[z,u]$. 
\qed\medskip

The proof of \ref{form} 
relies on the following lemma 
(cf. \cite{Sat1, KaZa1, Ve}). 

\query{(label=f(y)u+g(y,z))}\no\blem
\la{f(y)u+g(y,z)}
Let $p=f(y)u+g(y,z)$ be a variable of 
$\C[y,z,u]$. 

\begin{enumerate}
\item[{\rm (a)}] If $f\neq 0$ then $p$ 
is an $y$-variable, and it can be 
written as follows: 
$$
p=f(y)u+g_0(y)+g_1(y)z+f^{\rm red}(y)
\tilde g_2(y,z) z^2
$$
with ${\rm gcd} (f,g_1)=1$.
\item[{\rm (b)}] If $f=0$ 
then the coefficient of the highest 
order term in $z$ 
of $p$ (which {\it \`a priori} is a 
polynomial in $y$) is 
constant, unless $\deg_z g\le 0$. 
\end{enumerate}
\elem

\no \proof (a) The statement is evidently 
true if $f=$ const. 
Suppose that $f\notin\C$. By our assumption 
$\C^{[3]}/(p)\cong\C^{[2]}$. 
If $y_1\in\C$ is such that $f(y_1)=1$ 
then we have
$$
\C[y,z,u]/(f(y)u+g(y,z),y-y_1)\cong
\C[z,u]/(u+g_{y_1}(z))\cong\C^{[1]}. 
$$
By the Abhyankar-Moh-Suzuki Theorem, 
the same is true 
for every value of $y$.
In particular, for any root $y_0$ of $f$
we have
$$
 \C[y,z,u]/(f(y)u+g(y,z),y-y_0)\cong
\C[z,u]/(g_{y_0}(z))\cong \C^{[1]}\,.
$$
Therefore $\deg_z g_{y_0}(z)=1$, $p=fu+g$ 
has the desired form, and by \ref{RuVe}(a)
(with $B:=\C[y]$) 
it is an $y$-variable.

(b) In the course of the proof 
of \ref{corRuVe}
we have noticed that 
the polynomials of $\C[y,z]$ 
which are variables of $\C[y,z,u]$, 
are actually variables of $\C[y,z]$.
Thus $g\in \Var_{\C} \C[y,z]$, 
and so the statement (b) 
is well-known 
(e.g., see \cite{AM}).
\qed\medskip

\query{(label=prpr)}\no \bsit\la{prpr} 
{\it Proof of Proposition \ref{form}.} 
{\rm
It is easily seen that 
the polynomial $p=fu+g\in\C[x][y,z,u]$ 
admits a presentation
\query{(label=prst)}\be\la{prst}\qquad
p=q[\bar fu+\tilde g_0y^2+
\tilde g_1z+\bar g_2 z^2]
+a_0+a_1y\,
\ee
with $\quad a_0, a_1\in\C[x],\,\,\,\bar f, \tilde g_0,
\tilde g_1\in\C[x,y],\,\,\,\bar g_2
\in\C[x,y,z]$
and with $q\in\C[x]$ being a monic 
polynomial 
of maximal degree such that (\ref{prst}) 
holds.
Clearly,
$$q^{-1}(0)=L(p):=\{x_0\in\C\mid p_{x_0}
=f(x_0,y)u+g(x_0,y,z)\in \C y+\C\}$$
(possibly, this set is empty, 
and then we put $q:=1$).
As by our assumption, $f\notin\C[x]$ 
the subset  
$$
K=K(f):=\{\lambda\in \C\mid 
f_{\lambda}(y)\in \C\}
=q^{-1}(0)\cup \{\lambda\in \C\mid 
\bar f_{\lambda}(y)\in \C\}\subseteq  \C\,
$$
is finite.
We let $x_0\in\C\backslash K(f)$. 
The specialization
$p_{x_0}$ 
being a variable of $\C[y,z,u]$, 
by \ref{f(y)u+g(y,z)}(a) we have 
${\rm gcd} \bigl (\bar f(x_0,y),
\tilde g_1(x_0,y)\bigr )=1$, whence 
$\bar f^{-1}(0)\cap \tilde g_1^{-1}(0)
\subseteq  K\times \C_y$. 

Since by \ref{f(y)u+g(y,z)}(a),
$p_{x_0}\in\C[y,z,u]$ 
is a residual $y$-variable,
the equality 
$\bar f(x_0,y_0)=0$ implies that 
$\deg_z g_{x_0,y_0}(z)=1$, 
whence also
$\bar g_2(x_0,y_0,z)=0$.
It follows that  
$$\bar f^{-1}(0)\backslash 
(K\times \C_y)\subseteq 
\bar g_2^{-1}(0)\,.$$ 
It is easily seen that $\bar f\in\C[x,y]$ 
admits a factorization 
$\bar f(x,y)=r(x)\tilde f(x,y)$ 
such that $r^{-1}(0)\subseteq  K$ and 
$\tilde f^{-1}(0)\cap (K\times \C_y)$ 
is finite.
Then we have
$$\tilde f^{-1}(0)\backslash (K\times \C_y) 
\subseteq \bar g_2^{-1}(0)\,.$$
Passing to the Zariski closures 
we obtain  
$$\tilde f^{-1}(0) \subseteq 
\bar g_2^{-1}(0)\qquad {\rm 
i.e.,}\qquad \bar g_2=
\tilde f^{\rm red}\tilde g_2\,.$$
Now (\ref{prst}) yields a desired 
presentation 
(\ref{pform}). As 
$r^{-1}(0)\subseteq K(f)$, 
by the definition of $K(f)$ 
for every root
$x_0$ of $r$ we have $q(x_0)\tilde 
f(x_0,y)\in\C$,
and so it remains to check that 
$\tilde f^{-1}(0)\cap \tilde g_1^{-1}(0)$ 
is a finite subset of
 $q^{-1}(0)\times \C_{y}$. 

We already know that  
$\tilde f^{-1}(0)\cap \tilde g_1^{-1}(0)$ 
is a finite subset of 
$K\times \C_y$. Suppose that 
there is a point 
$$(x_0,y_0)\in \bigl[\tilde f^{-1}(0)\cap 
\tilde g_1^{-1}(0)\bigr] 
\backslash q^{-1}(0)\,.$$ 
Then 
$x_0\in K(f)\backslash q^{-1}(0)$, 
whence $\bar f(x_0,y)\in\C$, and so
$$\bar f(x_0,y)=\bar f(x_0,y_0)=
r(x_0)\tilde f(x_0,y_0)=0\,.$$ 
By our assumption,
$$
p_{x_0}= q(x_0)[\tilde g_0(x_0,y)y^2+
\tilde g_1(x_0,y)z+
\tilde f^{\rm red}(x_0,y)\tilde 
g_2(x_0,y,z)z^2]+
a_0(x_0)+a_1(x_0)y
$$
is a variable of $\C[y,z,u]$, 
and thus also of $\C[y,z]$.
Hence by \ref{f(y)u+g(y,z)}(b), 
$\tilde f^{\rm red}(x_0,y)=
\tilde f^{\rm red}(x_0,y_0)=0$, 
and then again by \ref{f(y)u+g(y,z)}(b),
$\tilde g_1(x_0,y)=\tilde g_1(x_0,y_0)=0$. 
Therefore, $p_{x_0}\in\C[y]$ 
is a variable, 
whence $q(x_0)=0$, a contradiction.\qed
}\esit

\query{(label=r(x))}\no\brem\label{r(x)}
In the notation of  \ref{form}, if
$p=fu+g=qr\tilde fu+g$ 
is a residual $x$-variable then  also 
$q\tilde fu+g$ is so. 
Moreover by \ref{corRuVe}(a), 
if $q\tilde fu+g$ 
is an $x$-variable then $p=fu+g$ 
is so as well.  
\erem

\query{label=var1}\no\blem\label{var1}
Let $p=fu+g\in\C[x][y,z,u]$ 
be a residual $x$-variable 
presented as in (\ref{pform}). 
If
$\tilde f^{-1}(0)\cap \tilde g_1^{-1}(0)
=\emptyset$ 
then $p$ is an $x$-variable.
\elem

\no\proof By \ref{r(x)} we may put  
in (\ref{pform}) $r=1$.
In virtue of \ref{RuVe}(a), the polynomial
$$
\tilde p:=\tilde f(x,y)u+\tilde g_0(x,y)y^2+
\tilde g_1(x,y)z+
\tilde f^{\rm red}(x,y)\tilde g_2(x,y,z)z^2
$$
is an ($x,y$)-variable i.e.,
there is an automorphism 
$\alpha\in\Aut_B B[z,u]$ 
(with $B:=\C[x,y]$) such that $\alpha(u)
=\tilde p$
(indeed by our assumption, $\tilde g_1$ 
is invertible $\mod \tilde f$ in $B$).
Furthermore by our assumption, 
for every root $x_0$ of
$q$ the specialization
$$p_{x_0}=f(x_0,y)u+g(x_0,y,z)=
a_0(x_0)+a_1(x_0)y\in\C[y,z,u]$$
is a variable. Consequently, 
$a_1(x_0)\neq 0$, 
and so ${\rm gcd} (q,a_1)=1$. 
It follows that
there is an affine automorphism 
$\beta\in \Aut_{\C[x]} \C[x][y,u]$ 
such that 
$$\beta(y)=q(x)u+a_0(x)+a_1(x)y\,.$$ 
Now $\alpha\beta\in\Aut_{\C[x]} 
\C[x][y,z,u]$
is such that
$$\alpha \beta(y)=q\tilde p+a_0(x)+a_1(x)y
=
fu+g=p\,,$$ as desired.
\qed\medskip

\query{label=deg1}\no\blem\label{deg1}
If $p=f(x,y)u+g(x,y,z)$ is a residual 
$x$-variable 
of degree at most 1 in
$z$ then it is an $x$-variable.
\elem

\no\proof It follows from \ref{form} that if 
$f\notin \C[x]$ then in (\ref{pform}) 
$\tilde g_1\neq 0$ 
(indeed, otherwise 
$\tilde f^{-1}(0)\cap\tilde g_1^{-1}(0)$ 
cannot be a finite subset of 
$q^{-1}(0)\times\C_y$). 
Thus if $\deg_z g\le 0$ then $f\in\C[x]$,
and so the statement follows from 
\ref{specase}.

If $\deg_z g=1$ then by \ref{form} $p=fu+g$ 
can be written as follows: 
$$
p=q(x)[r(x)\tilde f(x,y)u+\tilde 
g_0(x,y)y^2+\tilde 
g_1(x,y)z]+a_0(x)+a_1(x)y
$$
with $\tilde f^{-1}(0)\cap \tilde 
g_1^{-1}(0)$ 
being 
a finite subset of $q^{-1}(0)\times\C_y$. 
By \ref{r(x)} we may assume 
that $r=1$.

We proceed by induction on the 
intersection
index $n:=\tilde f^{*}(0)\cdot\tilde 
g_1^{*}(0)$.
If $n=0$ then the statement follows from 
\ref{var1}.
Suppose that $n\geq 1$.  

Let $x_0\in q^{-1}(0)$ be a coordinate of
an intersection point 
$(x_0,y_0)\in \tilde f^{-1}(0)\cap\tilde
g_1^{-1}(0)$. 
Up to a translation we may suppose that 
$x_0=0$.
As $$\min \{\deg_y \tilde g_1(0,y), 
\deg_y \tilde f(0,y)\}>0$$
we may write
$$
  \tilde g_1(0,y)z+\tilde f(0,y)u=
d(y)(a(y)z+b(y)u)\,,
$$
where ${\rm gcd} (a,b)=1$. 
Hence there is a linear automorphism 
$\gamma\in \Aut_B B[z,u]$
(with $B:=\C[x,y]$) such that 
$\gamma(a(y)z+b(y)u)=z$. 
It is not difficult to verify that
applying $\gamma$ amounts to: 
$$
      \mbox{Jac}(\gamma)
  \left [
    \begin{array}{c}
      \tilde g_1 \\
      \tilde f
    \end{array}
  \right ]
   = \left [
\begin{array}{c}
\hat g_1 \\
      x\hat f
    \end{array}
  \right ]
$$ for certain $\hat f, \hat g \in B$ 
(where Jac($\gamma)$ is the Jacobi 
matrix of $\gamma$), 
and hence 
$$
  \gamma(p)=\gamma(fu+g)=q(x)[x\hat f(x,y)u+
\tilde g_0(x,y)y^2+\hat g_1(x,y)z]+a_0(x)+
a_1(x)y\,.
$$
In view of
the equality of ideals: 
$$(\hat g_1, x\hat  f)=(\tilde g_1,
\tilde f)\,,$$
for the intersection indices we have:
$$\hat g_1^*(0)\cdot (x\hat f)^*(0)=
\tilde g_1^*(0)\cdot\tilde f^*(0)=n\,,$$ 
therefore 
$$\hat g_1^*(0)\cdot \hat f^*(0)\leq 
n-1\,.$$ 
Now in virtue of  \ref{r(x)}, 
$$
\hat p:=q(x)[\hat f(x,y)u+\tilde g_0(x,y)
y^2+
\hat g_1(x,y)z]+a_0(x)
+a_1(x)y
$$
is a residual $x$-variable.  
Hence by the inductive hypothesis, 
it is also an $x$-variable. Again by  
\ref{r(x)}, 
so are  
$\gamma(p)$ and $p=fu+g$ as well. 
The induction step is done, 
so the proof is completed.
\qed\medskip

\query{(label=prtr)}\no\bsit\la{prtr} 
{\it Proof of Theorem \ref{partialres}.}  
{\rm
By \ref{mthsec5}(v),(vi) we may suppose that 
$p=fu+g$ is a residual $x$-variable and
$\fnh \in \C[x]$ (cf. \ref{reva}). 
Thus we must show that
actually, under each of the assumptions (a)-(d) 
$p$ is an $x$-variable.

(a) (resp., (b)) immediately follows from 
\ref{specase} (resp., \ref{deg1}).

(c) We write the polynomial $p=fu+g$ in the form
(\ref{pform}). The assumption $f_{\rm vert}=1$ 
implies that $q^{-1}(0)=\emptyset$, and hence 
$\tilde f^{-1}(0)\cap \tilde g_1^{-1}(0)=
\emptyset$.
Now the conclusion follows from \ref{var1}.

(d) As by our assumption, 
$\fh=f_{{\rm horiz}}^{\rm red}$ 
we have $\tilde f=\tilde f^{\rm red}$. 
In view of \ref{r(x)} we may suppose that
$r=1$ in (\ref{pform}), and so $p$ can be 
written as follows: 
$$
p=fu+g=q[\tilde f u+
\tilde g_0 y^2+\tilde g_1 z+
\tilde f\tilde g_2 z^2]+
a_0+a_1 y\,.
$$
The $\C[x,y,z]$-automorphism 
$u\longmapsto u-z^2 \tilde g_2(x,y,z)$
transforms $p$ into a residual $x$-variable 
of degree at most $1$ in $z$. 
Now the conclusion follows from (b).
\qed
}\esit 
 
The following result extends our list of 
variables $p=fu+g$ beyond those 
provided in \ref{partialres}.

\query{(label=adpr)}\no\bprop \la{adpr}
Let $p=fu+g\in\C^{[4]}$ as in (\ref{pform}) 
be a residual $x$-variable.
If
$\tilde g_1\in\C[x]$ and 
$\tilde g_2z^2=\hat g_2(x,y,\tilde g_1(x)z, u)
\in\C[x][y, \tilde g_1 z]$ 
then $p$ is an $x$-variable.
\eprop

\no\proof By \ref{r(x)} 
we may suppose that 
$r=1$ in (\ref{pform}), and so 
$p=\bar p (x,y,\tilde g_1(x) z,u)$,
where the polynomial
\query{(label=nf)}\be\la{nf} 
\bar p:=q[\tilde f u+\tilde g_0 y^2
+ z +\tilde f^{\rm red} 
\hat g_2]+a_0+a_1 y\,\ee
still is a residual $x$-variable. 
Indeed by our assumption
for any
$\lambda\in\C\backslash\tilde g_1^{-1}(0)$,
$$\bar p_{\lambda}=
p_{\lambda}(y,\,z/\tilde g_1(\lambda),\, u)
\in\Var_{\C} \C^{[3]}\,.$$
For $x_0\in g_1^{-1}(0)$, either $q(x_0)=0$ 
and then by the assumption,
$\bar p_{x_0}=a_0(x_0)+a_1(x_0) y=
p_{x_0}\in\Var_{\C} \C^{[3]}$,
or $q(x_0)\neq 0$ and then again,
in virtue of \ref{RuVe}(a),
$\bar p_{x_0}\in \Var_{B} B[z,u]$ 
with $B:=\C[y]$. 
It follows from \ref{var1} that the
polynomial $\bar p$ as in (\ref{nf}) is
an $x$-variable. Hence by \ref{corRuVe}(a)
(with $z$ playing the role of $u$) so is 
$p=\bar p (x,y,\tilde g_1(x) z,u)$ as well.
\qed\medskip

\query{(label=x5)}\no\bexa\la{x5}
{\rm The polynomial
$$p:=y+x(xz+y(yu+x^2z^2))\in\C^{[4]}$$
is an $x$-variable. Indeed 
$p=\bar p (x,y,x z,u)$, where 
by \ref{RuVe}(a),
$\bar p :=y+x(z+y(yu+z^2))$ is a residual 
$x$-variable (and hence by \ref{var1} is
an $x$-variable). 
Therefore, $p$ is a
residual $x$-variable as well, 
and so \ref{adpr}
applies. 
}\eexa

\query{(label=QN}\no \bsit\la{QN}{\rm
To conclude, we would like to rise the 
following 
question. It concerns a (probably 
simplest) example of a residual 
$x$-variable of $\C^{[4]}$ 
which does not fit any of the assumptions
\ref{partialres}(a)-(d) or \ref{adpr}.} 
\medskip

\no {\bf Question:} 
{\it Is the polynomial 
$p=y+x(xz+y(yu+z^2))\in\C^{[4]}$ 
a variable?}\esit

\query{(label=1stably)}
\subsection{Simple modifications of 
$\C^3$ rectifiable in $\C^5$ }
\la{1stably}
The
following definitions are inspired by 
\cite{Sat1}.

\query{(label=hypl)}\no\bsit\la{hypl}{\rm
Let $B$ be a commutative ring. 
We say that $p\in B^{[n]}$ is a
{\it $B$-hyperplane} 
if $B^{[n]}/(p)\cong B^{[n-1]}$.
We say that $p$ defines a {\it $B$-hyperplane 
fibration} 
if $p-\lambda$ is a
$B$-hyperplane for every $\lambda\in B$.

If $B=\C$ or $B=\C[x]$ we simply say 
{\it hyperplane} resp.,  
{\it $x$-hyperplane} instead of {\it 
$B$-hyperplane}.

Notice that an $x$-hyperplane 
$p(x,y_1,\cdots,y_n)\in \C[x][y_1,\cdots,y_n]$ 
becomes
a hyperplane $p_{x_0}(y_1,\cdots,y_n)$ 
for every fixed $x=x_0\in \C$. 
A polynomial $p$ with the latter property is 
called a  
{\it residual $x$-hyperplane}.
}\esit

\query{(label=general)}\no\bsit\la{general}
{\rm Observe \cite{Sat1, KaZa2, Ve} that a 
polynomial $p=f(x,y)u+g(x,y,z)\in\C^{[4]}$ 
is a residual $x$-plane if and only if it is a 
residual $x$-variable, and henceforth, 
if and only if 
it defines a residual $x$-plane fibration. 
}\esit

\query{(label=1stab)}\no\bsit\la{1stab}{\rm
We say that two polynomials 
$p, \,q\in B[y_1,\cdots,y_n]$
are {\it 1-stably equivalent}  
if $$\gamma\bigl ((p,v)\bigr )=(q,v)$$
for an automorphism  $\gamma\in\Aut_B 
B[y_1,\cdots,y_n,v]$.
}\esit

\query{(label=fi)}\no\bsit\la{fi}{\rm
For instance, with this terminology 
\ref{isom}(b) 
says that 
$y+ap(y)$ with $a\in B,\,\,p\in B[y]$ 
is 1-stably equivalent to $y+p(ay)$, 
and moreover if $a^k|p(0)$ then 
$y+ap(y)$ is 1-stably 
equivalent to $y+{p(a^{k+1}y)}/{a^k}$.
}\esit

\query{(label=x-planes)}\bthm\label{x-planes}
Let $p=f(x,y)u+g(x,y,z)\in \C^{[4]}$ 
be a residual $x$-plane. 
Then it is 1-stably equivalent to an 
$x$-variable.
Consequently, $p$ defines an $x$-plane 
fibration,
and each fiber 
$X_{\lambda}:=\{p=\lambda\}\,\,\,(\lambda\in\C)$
of $p$ is a $3$-fold in $\C^4$ 
isomorphic to $\C^3$ which 
can be rectified in $\C^5$.    
\ethm

\query{(label=anpr)}\no\brem\la{anpr}
Notice that \ref{x-planes} provides yet 
another proof 
of the implication 
(iii)$\Rightarrow$(i) of \ref{mthsec5} 
which does not need
\ref{miy}. Indeed, it proves (v)$\Rightarrow$(i) 
while
(iii)$\Rightarrow$(v) is observed in \ref{general} 
above. 
\erem

In the proof of \ref{x-planes} 
(see \ref{prto} below)
we use the following lemma.

\query{(label=Cy+C)}\blem\la{Cy+C}
If $p\in \C[x,y,z,u]\setminus\C[x,y]$ 
is a residual $x$-plane 
then it is 1-stably equivalent 
(under an automorphism $\gamma\in\Aut_B B[y,v]$
with $B:=\C[x,z,u]$) 
to a polynomial 
$\bar p\in\C[x][y,z,u]$ which is still 
a residual $x$-plane and
such that for every 
$x_0\in \C,\;\;\;\bar p_{x_0}(y,z,u)\notin 
\C^*y+\C$. 
\elem

\no\proof As $p\notin\C[x,y]$, the subset 
$$L(p):=\{x_0\in \C\mid p_{x_0}\in\C^*y+\C\}=
:\{x_1,\cdots,x_n\}\subseteq\C$$
is finite. The proof proceeds by
a descending
induction on $n=$card$\,L(p)$.

If card$\,L(p)=0$ there is nothing to prove.
Suppose that $n>0$, and let $x_n=0$. 
Then there exists an affine automorphism 
$\tau\in\Aut_{\C[x]} \C[x][y]\subseteq\Aut_{D} 
D[y]$ 
(where $D:=\C[x,z,u,v]$) such that 
$\tau(p)=y+x\hat p$, 
where $\hat p\in \C[x,y,z,u]$ and 
$\hat p(x,0,0,0)=0$. 
Indeed, up to an affine automorphism
of $\C[y]$ we may assume that 
$p(0,y,z,u)=y$ i.e., 
$p=y+x\tilde p$ with $\tilde p\in\C^{[4]}$,
and then we apply the shift $y\longmapsto 
y-x\tilde p(x,0,0,0)$ to
obtain $\tau(p)=y+x\hat p$ with
$\hat p:=\tilde p-\tilde p(x,0,0,0)$. 
 
The polynomial 
$\hat p\in\C^{[4]}$ admits a
presentation:
$$
\hat p=yq_1(x,y,z,u)+x^k(xq_2(x,z,u)+
q_3(z,u))\,.
$$
As by our assumption
for every fixed $x\in\C$ the polynomial
$\tau(p)=y+x\hat p\not\in\C[x,y]$ is irreducible, 
we have
$xq_2+q_3\neq 0$, 
and thereby $q_3\neq 0$. Furthermore,
$$0=\hat p(x,0,0,0)=
x^k\bigl(xq_2(x,0,0)+q_3(0,0)\bigr)=q_3(0,0)\,,$$
whence $q_3\notin\C$.
Now by \ref{isom}(b) (cf. \ref{fi} above) 
there exists an automorphism 
$\alpha\in\Aut_B B[y,v]$
(with $B=\C[x,z,u]$) such that 
$$
\alpha\bigl ((y+x\hat p,v)\bigr )=(y+
\frac{\hat p(x,x^{k+1}y,z,u)}{x^k},v)\,.
$$ 
Letting $\gamma_n:=\alpha\tau$ we obtain:
$$
  \gamma_n\bigl ((p,v)\bigr )=(p_n,v)\,,
$$
where 
$p_n:=y+{\hat p(x,x^{k+1}y,z,u)}/{x^k}$.
We have:
\begin{eqnarray*}
  p_n & = & y+\frac{x^{k+1}yq_1(x,x^{k+1}y,z,u)+
x^k(xq_2(x,z,u)+q_3(z,u))}{x^k}\\
               &=  & 
y+xyq_1(x,x^{k+1}y,z,u)+xq_2(x,z,u)+q_3(z,u)\,.
\end{eqnarray*}
Hence $p_n(0,y,z,u)=y+q_3(z,u)\in\Var_{\C} 
\C[y,z,u]$ 
with $q_3\notin\C$, and so
$ p_n(0,y,z,u)\notin \C^*y+\C$. 
In other words, 
$$0=x_n \notin L(p_n):=
\{\lambda\in \C\mid (p_n)_{\lambda}\in
\C^*y+\C\}\subseteq\C\,.$$ 
Now the induction step  
and the proof of the lemma are completed
by the following \smallskip

\no {\bf Claim.} {\rm (a)} $\quad 
L(p_n)= \{x_1,\cdots,x_{n-1}\}=L(p)
\backslash\{x_n\}\,.$

\begin{enumerate}\item[(b)] Furthermore 
for any $x_0\neq 0$, 
$(p_n)_{x_0}\in\C[y,z,u]$ is a plane.
\end{enumerate}
\medskip

\no {\it Proof of the claim.} 
\footnote{Alternatively,  
(b) also follows from existence of 
a $\C[x]$-isomorphism  
$\C[x]^{[3]}/(p)\cong \C[x]^{[3]}/(\bar p)$, 
see \ref{isomm}(b).}
It is enough to show that for every $x_0\neq 0$ 
there exists 
an affine automorphism 
$\beta_{x_0}\in \Aut_{\C} \C[y]
\subseteq\Aut_{\C[z,u]} \C[z,u][y]$ 
such that 
\query{(label=ap)}\be\la{ap}
\beta_{x_0}(p_{x_0})=x_0^{k+1}(p_n)_{x_0}\,.\ee
Indeed, we have:
\query{(label=pa)}\be\la{pa} 
x_0^{k+1} (p_n)_{x_0} = x_0^{k+1}y+
x_0\hat p_{x_0}(x_0^{k+1}y,z,u)\,.\ee
Consider the linear automorphism 
$\si_{x_0} : y\longmapsto x_0^{k+1}y$ 
extended to $\C[z,u][y]$ in a natural way.
From (\ref{pa}) we get:
$$x_0^{k+1} (p_n)_{x_0} = 
\si_{x_0} \bigl (y+x_0\hat p_{x_0}(y,z,u)\bigr )=
\si_{x_0} 
\bigl ((\tau(p))_{x_0}\bigr )=
\si_{x_0}\bigl ((\tau_{x_0}(p_{x_0})\bigr )\,,$$
where $\tau_{x_0}$ denotes the specialization 
of $\tau$ at $x_0$. Thus we obtain (\ref{ap})
with 
$\beta_{x_0}:=\si_{x_0} \tau_{x_0}\in\Aut_{\C} 
\C[y]$
extended to $\C[z,u][y]$. 
Hence $p_{x_0}\in \C^*y+\C$ if and only 
if so is
$(p_n)_{x_0}$. 
\qed\medskip

\query{(label=prto)}\bsit\la{prto}
{\it Proof of Theorem \ref{x-planes}.} {\rm
If $f=0$ then in virtue of \ref{f=0}(c)
and the Abhyankar-Moh-Suzuki Theorem,
$p$ itself is an $x$-variable. 

If $f\neq 0$ then 
by \ref{Cy+C}, $p=fu+g$ is 1-stably equivalent 
to a polynomial 
$\bar p\in\C^{[4]}$ such that
\query{(label=L)}\be\la{L}
L(\bar p):=\{\lambda\in\C\mid \bar p_{\lambda}\in 
\C^*y+\C\}=\emptyset\,.\ee 
It is easily seen that the 
polynomial $\bar p$ as constructed in the proof of
\ref{Cy+C} still has the form 
$\bar p  =\bar f u +\bar g$ with 
$\bar f\in\C[x,y]$ and $\bar g\in \C[x,y,z]$.
It follows from (\ref{L})
that ${\bar f}_{\rm vert}=1$. 
As $\bar p$ is also a residual $x$-plane, by 
\ref{partialres}(c) it is 
an $x$-variable. 
Consequently,
\begin{eqnarray*}
\C[x][y,z,u]/(p)& \cong & \C[x][y,z,u,v]/(p,v)
 \cong  \C[x][y,z,u,v]/(\bar p,v)\\
&\cong &\C[x][y,z,u,v]/(u,v)  \cong  \C[x]^{[2]}\,.
\end{eqnarray*} Hence $p$ is an $x$-plane. 
In view of \ref{general} the same arguments work 
for every `fiber' $p-\lambda$ of $p=fu+g$ 
(with $\lambda\in\C[x]$). Therefore $p=fu+g$ 
defines an 
$x$-plane fibration.
\qed
}\esit

\query{(label=SW)}\subsection{On Sathaye-Wright's 
Theorem}
\la{SW}
Recall that the Sathaye Theorem on linear planes
\cite{Sat1}
cited in the introduction was generalized by
D. Wright \cite{Wr} for the embeddings
$\C^2\hookrightarrow\C^3$ of the form 
$p=f(x,y)u^n+g(x,y)=0$,
and further generalized in \cite{KaZa1} for
the acyclic surfaces in $\C^3$ of this type. 
Here we prove the following theorem 
(based on the latter result).

\query{(label=Wr)}\no\bthm\la{Wr}
If $p=f(x,y,z)u^n+g(x,y,z)\in\C^{[4]}$ 
(with $n\ge 2$) 
is a 
residual $x$-variable then it is an 
$x$-variable.\ethm

\no\proof By \cite[Thm. 7.2]{KaZa1}, 
the polynomials $f$ and $g$ 
satisfy the following condition:
\query{(label=con)}\be\la{con}
\forall \lambda\in\C\quad \exists 
\alpha_{\lambda}\in
\Aut_{\C} \C[y,z]
\quad \mbox{and}\quad \exists 
\varphi_{\lambda}\in 
\C[y]
\qquad \mbox{such that}\ee $$  
f_{\lambda}=\alpha_{\lambda}(\varphi_{\lambda})
\quad \mbox{and}\quad 
g_{\lambda}=\alpha_{\lambda}(z)\,.$$
Hence by \ref{f=0}(c), 
$g$ is an $x$-variable of $\C[x,y,z]$,
and so we may suppose that $g=z$ and 
$\alpha_{\lambda}\in\Aut_{\C[z]} \C[z][y]$. 
Then (\ref{con}) yields: 
\query{(label=ccon)}\be\la{ccon}
\forall \lambda\in\C\quad 
\exists q_{\lambda}\in\C[z]
\quad \mbox{and}\quad \exists 
\varphi_{\lambda}\in \C[y]
\qquad \mbox{such that}\ee 
$$  
f_{\lambda}=\alpha_{\lambda} (\varphi_{\lambda})=
\varphi_{\lambda} (\alpha_{\lambda} (y))=
\varphi_{\lambda} (y+q_{\lambda}) \,.
$$
Eventually, we prove below (by induction 
on $\deg_y f$) 
that
for any polynomial $f\in\C^{[3]}$
which satisfies (\ref{ccon}), we have:
\query{(label=icon)}\be\la{icon} 
\quad f=\Phi(x,a(x)y+zQ(x,z))\ee
with $a\in\C[x],\,\, Q\in\C[x,z]$
and $\Phi\in \C[x,y]$.
Therefore $p=P(x,a(x)y,z,u)$ 
is a residual $x$-variable 
with $P:=\Phi(x,y')u+z\in
\Var_{\C[x]}\C[x]^{[3]}$ 
(where $y':=y+zQ(x,z)$).  
Thus by \ref{corRuVe}(a), 
$p$ is an $x$-variable, as
stated.

If $\deg_y f\le 0$ then (\ref{ccon}) 
implies that
$f(x,y,z)=f(x,0,0)\in\C[x]$ and (\ref{icon})
follows. 
Assume further that $\deg_y f\ge 1$. 
By (\ref{ccon}) we have: 
$$\partial_y f_{\lambda}(y,z) = 
(\varphi_{\lambda})'(y+q_{\lambda}(z))\,.$$
It follows that the polynomial 
$\partial_y f$ with 
$\deg_y \partial_y f= \deg_y f-1$
also satisfies (\ref{ccon}), whence by the 
inductive hypothesis,  
$$ 
\partial_y f=
\tilde \Phi(x,\tilde a(x)y+z\tilde Q(x,z))
$$
with $\tilde a\in\C[x],\,\, 
\tilde Q\in\C[x,z]$ and $\tilde \Phi\in \C[x,y]$.
Thereby, 
$$\tilde af=\Psi (x,\tilde a(x)y
+z\tilde Q(x,z))+zR(x,z)\qquad 
\mbox{with}\quad\Psi\in\C[x,y]
\quad\mbox{and}\quad R\in\C[x,z]\,,$$
where $\partial_y \Psi=\tilde \Phi$. 
If $\deg_y f=1=\deg_y \Psi$ then 
$$
  \tilde a f =\Psi_0(x)+\Psi_1(x)[\tilde a(x)y
+z\tilde Q(x,z)]+zR(x,z)\,,
$$
and so 
$\tilde a$ divides $\Psi_0(x)
+\Psi_1(x)z\tilde Q(x,z)+zR(x,z)$ 
i.e.,
$$\Psi_0(x)+\Psi_1(x)z\tilde Q(x,z)+zR(x,z)
=\tilde a(zQ(x,z)+\Phi_0)\,.$$ It follows that
$$
   f = \Psi_1(x)y+zQ(x,z)+\Phi_0
$$ 
has the desired form (\ref{icon}).
 
Thus we may assume that $\deg_y f\ge 2$.
By (\ref{ccon}) for every $\lambda\in\C$ 
we obtain:
\query{(label=econ)}\be\la{econ}
\tilde a(\lambda )f_{\lambda}=
\Psi (\lambda,\tilde a(\lambda)y+
z\tilde Q(\lambda,z))+zR(\lambda,z)=\tilde a
(\lambda)
\varphi_{\lambda}(y+q_{\lambda}(z))\,,\ee
and so
$$\Psi (\lambda,\tilde a(\lambda)y-
\tilde a(\lambda)q_{\lambda}(z)+z
\tilde Q(\lambda,z))+zR(\lambda,z)=\tilde a
(\lambda)
\varphi_{\lambda}(y)\,.$$
As $\deg_y \Psi=\deg_y f\ge 2$ it follows that 
for every $\lambda\in\C$ such that 
$\deg_y \Psi(\lambda,y)\geq 0$ we have
$\deg_z (z\tilde Q(\lambda,z)-
\tilde a(\lambda) q_{\lambda}(z))\le 0$,
and so $\deg_z zR(\lambda,z)\le 0$ as well i.e., 
$R(\lambda,z)=0$ for almost every 
$\lambda\in \C$. 
Hence  $R=0$, and from (\ref{econ}) 

we obtain:
\query{(label=NN)}\be\la{NN}
\tilde af=\Psi (x,\tilde a(x)y+z\tilde Q(x,z))\,.
\ee
If $\tilde Q=0$ then 
$f=\Psi (x,\tilde a(x)y)/\tilde a \in\C[x,y]$, 
as 
desired. If $\tilde Q\neq 0$ then (up to
factorizing $\tilde a(x)y+z\tilde Q(x,z)$ 
and  
changing appropriately $\Psi (x,y)$) 
we may assume that 
gcd$\,(\tilde a,\tilde Q(x,z))=1$. 
Putting in (\ref{NN}) $y=0$ gives:
$$
\tilde a(x)f(x,0,z)=\Psi (x,z\tilde Q(x,z))\,.
$$
From this equality and the assumption that
gcd$\,(\tilde a,\tilde Q(x,z))=1$ 
it is not difficult to deduce that 
$\tilde a$ divides 
$\Psi$ in $\C[x,y]$. Thus from (\ref{NN}) 
with $\Phi:=\Psi/\tilde a$ 
we get  
$$
  f=\Phi (x,\tilde a(x)y+z\tilde Q(x,z))\,,
$$
as needed.
\qed

\end{document}